\documentclass[11pt,a4paper]{article}
\usepackage{amssymb,amsmath}

\numberwithin{equation}{section}

\usepackage{times,theorem,latexsym,color}

\newcommand{\BOX}{\ensuremath\Box}

\newtheorem{theorem}{Theorem }[section]

{\theorembodyfont{\rmfamily}}
{\theorembodyfont{\rmfamily}}
{\theorembodyfont{\rmfamily}}
\newtheorem{lemma}[theorem]{Lemma}
\newtheorem{proposition}[theorem]{Proposition}
{\theorembodyfont{\rmfamily}\newtheorem{remark}[theorem]{Remark}}
{\theorembodyfont{\rmfamily}}

\newcommand{\N}{\mathbb{N}}

\newcommand{\C}{\mathbb{C}}
\newcommand{\R}{\mathbb{R}}
\newcommand{\Z}{\mathbb{Z}}
\newcommand{\dd}{\,{\rm d}}
\newcommand{\Ai}{{\rm Ai}}

{\hspace*{.1pt}\hspace*{\fill}\BOX\vskip\baselineskip}

{\vskip\baselineskip\noindent\textbf{Proof of {#1}:}}%
{\hspace*{.1pt}\hspace*{\fill}\BOX\vskip\baselineskip}
{\vskip\baselineskip\noindent\textbf{Proof of Theorem \protect\ref{#1}:}}%
{\hspace*{.1pt}\hspace*{\fill}\BOX\vskip\baselineskip}
{\vskip\baselineskip\noindent\textbf{Proof of Lemma~\protect\ref{#1}:}}%
{\hspace*{.1pt}\hspace*{\fill}\BOX\vskip\baselineskip}
\newenvironment{proofprop}[1]%
{\vskip\baselineskip\noindent\textbf{Proof of Proposition~\protect\ref{#1}:}}%
{\hspace*{.1pt}\hspace*{\fill}\BOX\vskip\baselineskip}
{\vskip\baselineskip\noindent\textbf{Proof of Theorems \protect\ref{#1} --
\protect\ref{#2}:}}%
{\hspace*{.1pt}\hspace*{\fill}\BOX\vskip\baselineskip}

\begin{document}

\title{On stationary two-dimensional flows around a fast rotating disk}

\author{Isabelle Gallagher\thanks{IMJ-PRG, Universit\'e Paris-Diderot and DMA, Ecole Normale Sup\'erieure de Paris, France. E-mail: \texttt{gallagher@math.ens.fr}}  \and Mitsuo Higaki\thanks{Department of Mathematics, Kyoto University, Japan. E-mail: \texttt{mhigaki@math.kyoto-u.ac.jp}} \and Yasunori Maekawa\thanks{Department of Mathematics, Kyoto University, Japan. E-mail: \texttt{maekawa@math.kyoto-u.ac.jp}}}


\maketitle

\noindent {\bf Abstract} 
We study the two-dimensional stationary Navier-Stokes equations 
describing   flows around a rotating disk. 
The   existence of unique solutions is established for any rotating speed, and qualitative effects of a large rotation are described precisely by exhibiting a boundary layer structure and an axisymmetrization of the flow.
\vspace{0.3cm}

\noindent {\bf Keywords}\, Navier-Stokes equations $\cdot$ two-dimensional exterior flows $\cdot$ flows around a rotating obstacle $\cdot$ boundary layer

\vspace{0.3cm}

\noindent {\bf Mathematics Subject Classification (2000)}\, 35B35 $\cdot$ 35Q30 $\cdot$ 76D05 $\cdot$ 76D17


\section{Introduction}\label{intro}

Understanding the structure of   flows generated by the motion of a rigid body is an important subject in fluid dynamics. The typical motions of a rigid body (called the obstacle below) are translation and rotation, and it is well known that the flow generated by an obstacle translating with   constant speed possesses a wake structure behind it, and that this structure is well described by the Oseen approximation; see Galdi \cite{Ga} for   mathematical results on this topic. Flows around a rotating obstacle have been   studied mathematically   mainly in  three space dimensions,
while there are   only   few mathematical results for   two dimensions.
Moreover, in the two dimensional case most   known mathematical results are restricted to the case when the Reynolds number is sufficiently small.
In general, the motion of the obstacle leads to a drastic change of the decay structure of the flow,
and in the two dimensional case it yields a significant localizing effect that enables one to construct  corresponding steady state solutions when the motion of the obstacle is slow enough.
Although on the one hand a faster motion of the obstacle might give a stronger localizing and stabilizing effect, on the other hand   fast motion produces a rapid flow and creates a strong shear near the boundary that can be a source of   instability. As a result,   rigorous  analysis becomes quite difficult for the nonlinear problem in general.
Hence it is useful to study the problem under a simple geometrical setting and to understand a typical fluid structure that describes these two competitive mechanisms; localizing and stabilizing effects on the one hand, and the presence of a rapid flow and the boundary layer created by the fast motion of the obstacle on the other hand.

In this paper we study   two dimensional flows around a rotating obstacle   assuming that the obstacle is a unit disk centered at the origin, in the case when the rotation speed is sufficiently fast and   the Reynolds number is high.
Note that in a three dimensional setting, these flows are considered as a model for   two dimensional flows around a rotating infinite cylinder with a uniform cross section which is a unit disk.

Let us formulate our problem mathematically. 
For notational convenience we first state the problem when the obstacle has a general shape.
Let us consider the following Navier-Stokes equations for viscous incompressible flows in two dimensions:

\begin{equation}\label{NS}
  \left\{
\begin{aligned}
 \partial_t w -\Delta w + w\cdot \nabla w + \nabla \phi & \,=\, g \,,  \quad {\rm div}\, w \,=\, 0\,,   \qquad t>0\,,~y \in \Omega (t)\,, \\
  w & \,=\, \alpha y^\bot \,, \qquad \qquad \qquad ~ t>0\,,~ y \in \partial \Omega (t)\,. \\
\end{aligned}\right.
\end{equation}

\noindent Here $w=w(y,t) = (w_1(y_1,y_2,t), w_2 (y_1,y_2,t))^\top$ and $\phi=\phi (y,t)$ are respectively the unknown velocity field and pressure field, and $g(y,t) = (g_1(y,t), g_{2}(y,t))^\top$ is an external force. 
The vector $y^\bot$ is defined as $y^\bot=(-y_2,y_1)^\top$. Here $M^\top$ denotes the transpose of a matrix $M$. 
 The domain $\Omega(t)$ is given by

\begin{align}
\begin{split}
\Omega(t) & \, := \, \big \{ y\in \R^2~|~ y = O (\alpha t) x\,, \ \ x:= (x_1,x_2)^\top \in \Omega \big \}\,,\\
O (\alpha t ) & \, := \, 
\begin{pmatrix}
\cos \alpha t & -\sin \alpha t\\
\sin \alpha t & \cos \alpha t
\end{pmatrix}\,,
\end{split}
\end{align}

\noindent 
where $\Omega$ is an exterior domain in $\R^2$ and its complement $\Omega^c:=\R^2\setminus \Omega$ describes the obstacle at initial time, while the real number $\alpha$ represents the rotation speed of the obstacle.
The condition $w (t,y)=\alpha y^\bot$ on the boundary  $\partial \Omega (t)$  represents the noslip boundary condition. 
Taking into account the rotation of the obstacle, we introduce the following change of variables and unknowns:
\begin{align*}
y = O (\alpha t) x\,, \quad u (x,t) = O (\alpha t)^\top w (y,t)\,, \quad p (x,t) = \phi (y,t)\,, \quad f (x,t) = O (\alpha t)^\top g (y,t)\,.
\end{align*}
Then \eqref{NS} is equivalent to the equations in the time-independent domain $\Omega$:
\begin{equation*}
  \left\{
\begin{aligned}
 \partial_t u -\Delta u - \alpha ( x^\bot \cdot \nabla u - u^\bot ) + \nabla p & \,=\, 
  - u\cdot\nabla u  + f \,,  \quad {\rm div}\, u \,=\, 0\,,   \quad t>0\,,~x \in \Omega\, , \\
  u & \,=\, \alpha x^\bot \,, \qquad \qquad t>0\,,~ x \in \partial \Omega\,. \\
\end{aligned}\right.
\end{equation*}
We are interested in   stationary solutions to this system when the obstacle $\Omega^c$ is a unit disk centered at the origin, and therefore we assume that $f$ is independent of $t$ and consider the elliptic system
\begin{equation}\tag{NS$_{\alpha}$}\label{NS_alpha.intro}
  \left\{
\begin{aligned}
  -\Delta u - \alpha ( x^\bot \cdot \nabla u - u^\bot ) + \nabla p & \,=\, 
  - u\cdot\nabla u + f \,,   \qquad x\in \Omega \,,\\
{\rm div}\, u &  \,=\, 0\,,   \qquad \qquad  \quad  \quad \ \ x \in \Omega\,, \\
  u & \,=\, \alpha x^\bot \,, \qquad \qquad  \quad ~ x \in \partial \Omega\,, \\
\end{aligned}\right.
\end{equation}
with $\Omega := \{x\in \R^2~|~|x|>1\}$. We note that in the original coordinates the stationary solution to \eqref{NS_alpha.intro} gives a specific time periodic flow with a periodicity $ \frac{2\pi}{|\alpha|}$.
Due to the symmetry of the domain  there is an explicit stationary solution to \eqref{NS_alpha.intro} when $f=0$:
\begin{align}
\big ( \alpha U , \alpha^2 \nabla P \big)  ~~~~~{\rm with}~~~~~ U (x) := \frac{x^\bot}{|x|^2}\,,~~~~~ P (x) := -\frac{1}{2 |x|^2}\,\cdotp  \label{Sol}
\end{align}
Thus it is natural to consider an expansion around this explicit solution. By using the identity $u\cdot \nabla u = \frac12 \nabla |u|^2 + u^\bot {\rm rot} \, u$ with ${\rm rot}\, u :=\partial_1 u_2 - \partial_2 u_1$ and ${\rm rot}\, U =0$ for $x\ne 0$, 
the equations for~$v:= u- \alpha U$ are written as 
\begin{equation}\tag{$\widetilde {\rm NS}_{\alpha}$}\label{tildeNS_alpha}
  \left\{
\begin{aligned}
  -\Delta v - \alpha ( x^\bot \cdot \nabla v - v^\bot ) + \nabla q + \alpha U^\bot {\rm rot}\, v & \,=\, 
  - v^\bot {\rm rot}\, v   + f \,,   \quad x\in \Omega \,,\\
{\rm div}\, v &  \,=\, 0\,,   \qquad \quad  x \in \Omega\,, \\
  v & \,=\, 0 \,, \qquad \quad  x \in \partial \Omega\,. \\
\end{aligned}\right.
\end{equation}
The goal of this paper is to show   the  existence and uniqueness of solutions to \eqref{tildeNS_alpha} for arbitrary~$\alpha \in \R\setminus\{0\}$ under a suitable condition on the given external force $f$ in terms of regularity and summability, and   moreover, we shall give a detailed qualitative analysis for the fast rotation case $|\alpha|\gg 1$.

Before stating our results let us recall known mathematical results related to this problem.
Flows around a rotating obstacle have been studied mainly for the three dimensional case.
Borchers \cite{Bo} proved the existence of global weak solutions for the nonstationary problem, 
and  the   existence and uniqueness of local in time regular solutions is shown by Hishida~\cite{H1} and Geissert, Heck, and Hieber \cite{GHH}. Global strong solutions for small data are obtained by Galdi and Silvestre \cite{GSi}. The spectrum of the linear operator related to this problem is studied by Hishida \cite{H2}, and Farwig and Neustupa \cite{FN}. The existence of stationary solutions to the three dimensional problem is proved in \cite{Bo}, Galdi~\cite{G1}, Silvestre \cite{Si}, and Farwig and Hishida \cite{FH0}. In particular, \cite{G1} constructs   stationary flows with  a decay of order  $O(|x|^{-1})$, while~\cite{FH0} discusses  the weak $L^{3}$ framework.
The asymptotic profiles of these stationary flows at spatial infinity are studied by Farwig and Hishida \cite{FH1,FH2} and Farwig, Galdi, and Kyed \cite{FGK}, where it is proved that the asymptotic profiles are described by the Landau solutions, which are stationary self-similar solutions to the Navier-Stokes equations in 
$\R^3\setminus\{0\}$. The stability of   small stationary solutions has been well studied in the three-dimensional case. Indeed,  the global $L^2$ stability is proved in \cite{GSi} and   local $L^3$ stability is obtained by Hishida and Shibata \cite{HShi}.

All the results mentioned above are for the three dimensional case.
So far there have been only   few results in the two dimensional case. 
Hishida \cite{H3} revealed the asymptotic behavior of the stationary Stokes flow around a rotating obstacle, and showed that the rotation of the obstacle leads to the resolution of the Stokes paradox as in the case of the translation of the obstacle.
On the other hand, the linear result in \cite{H3} was not sufficient to solve the nonlinear problem due to the nature of the singular perturbation limit $\alpha \rightarrow 0$ which is specific to the two dimensional case.
Recently the linear result of \cite{H3} was extended by Higaki, Maekawa, and Nakahara \cite{HMN}, where the  existence and uniqueness of solutions to~\eqref{NS_alpha.intro} decaying at the scale-critical order $O(|x|^{-1})$ was also proved when the rotation speed~$\alpha$ is sufficiently small and the external force $f$ is of a divergence form $f={\rm div}\, F$ for some~$F$ which is small in a scale critical norm. Moreover, the leading profile at spatial infinity was shown to be~$C \frac{x^\bot}{|x|^2}$ for some constant $C$ under an additional decay condition on~$F$ such as~$F=O(|x|^{-2-r})$ with~$r>0$, which is compatible with the work of \cite{H3} for the Stokes case. 
The stability of these stationary solutions, though they are small in a scale-critical norm, is a difficult problem and is still largely open in the two dimensional case. The only known result is by Maekawa~\cite{Ma} for a specific case, which shows the local~$L^2$ stability of the explicit solution~\eqref{Sol} in the original frame~\eqref{NS} when $\alpha$ is small enough and~$\Omega$ is the exterior to the unit disk as assumed in this paper. 
Few results are known   in the case of nonsmall~$\alpha$. As a related work in this direction, Hillairet and Wittwer \cite{HW} established the existence of stationary solutions to \eqref{NS} for large $|\alpha|$ when $\Omega(t)=\{y\in \R^2~|~|y|>1\}$ as in this paper, although  in \cite{HW} the boundary condition is $\alpha y^\bot + b$ instead of $\alpha y^\bot$, with some time-independent given data $b$ and the external force $f$ is assumed to be zero. In \cite{HW} the stationary solution is constructed around \eqref{Sol} when $|\alpha|$ is large and $b$ is small. Our problem is in fact essentially different from the one discussed in \cite{HW}. Indeed, the stationary solution to \eqref{NS_alpha.intro} is a time periodic solution in the original frame \eqref{NS}, and therefore, the result in~\cite{HW}  is not applicable to our problem and vice versa.
The reader is also referred to~\cite{HMN2} for a model problem of \eqref{NS_alpha.intro} formulated in the whole plane $\R^2$, where the  existence and uniqueness of stationary solutions is proved when $\alpha\ne 0$. However, the argument in \cite{HMN2} relies on the absence of   physical boundary, and the key mechanism originating from the boundary is not analyzed.

Let us now return to   \eqref{tildeNS_alpha}.
The novelty of the results in this paper is the followings:

\medskip

\noindent (1) Existence and uniqueness of solutions to \eqref{tildeNS_alpha} for arbitrary $\alpha\in \R\setminus\{0\}$.

\noindent (2) Relaxed summability condition on  $f$ and on the class of solutions, allowing slow spatial decay with respect to  scaling.

\noindent (3) Qualitative analysis of solutions in the fast rotation case $|\alpha|\gg 1$.

\medskip

\noindent As for (1), the result is new compared with \cite{HMN} in which the stationary solutions are obtained only for nonzero but small $|\alpha|$, though there is no restriction on the shape of the obstacle in \cite{HMN}.
The reason why we can construct solutions for all  nonzero $\alpha$ is a remarkable coercive estimate for the term $-\alpha ( x^\bot \cdot \nabla v - v^\bot ) + \alpha U^\bot {\rm rot}\, v$ in   polar coordinates; see \eqref{proof.apriori.middle.intro.1} below. As for (2), we note that the given data $f$ and the class of solutions in \cite{HMN} are in a scale critical space.
A typical behavior for $f$ assumed in \cite{HMN} is that $f={\rm div}\, F$ with $F(x)=O(|x|^{-2})$, and then the solution $v$ satisfies the estimate $|v(x)|\leq C|x|^{-1}$ for~$|x|\gg 1$. In this paper the summability condition on $f$ is much weaker than this scaling, see \eqref{def.Y} below. Moreover, the radial part of the solution constructed in this paper only behaves  like $o(1)$ as $|x|\rightarrow \infty$ in general, which is considerably  slow, while the nonradial part of the solution belongs to $L^2(\Omega)$ which is just in the scale critical regime.
The point (3) is important both physically and mathematically.   Understanding the fluid structure around the fast rotating obstacle up to the boundary is one of the main subjects of this paper, and we   show the appearance of a boundary layer as well as an axisymmetrization mechanism due to the fast rotation of the obstacle.

Let us state our functional setting.
Due to the symmetry of the domain it is natural to introduce the relevant function spaces in terms of   polar coordinates.
As usual, we set
\begin{align*}
& x_1 \, = \, r \cos \theta\,,~~~~~ x_2 \, = \, r\sin \theta\,,~~~~~~~~~~ r \, = \, |x| \geq 1\,,~~~\theta\in [0,2\pi)\,,\\
& {\bf e}_r \, = \, \frac{x}{|x|}\,, \qquad {\bf e}_\theta \, = \, \frac{x^\bot }{|x|} \, =\, \partial_\theta {\bf e}_r\,,
\end{align*}
and 
\begin{align*}
& v \, = \, v_r \,  {\bf e}_r \, + \,  v_\theta \,  {\bf e}_\theta\,, \qquad \quad v_r \, = \, v\cdot {\bf e}_r\,, \qquad  v_\theta \, = \, v\cdot {\bf e}_\theta\,.
\end{align*}
Next, for each $n\in \Z$, we denote by $\mathcal{P}_n$ the projection on the Fourier mode $n$ with respect to the angular variable $\theta$:
\begin{align}
\mathcal{P}_n v \,: = \, v_{r,n} e^{i n \theta} {\bf e}_r + v_{\theta,n} e^{i n \theta} {\bf e}_\theta\,, \label{def.P_n}
\end{align}
where 
\begin{align*}
v_{r,n} (r) & \,: = \, \frac{1}{2\pi} \int_0^{2\pi} v_r (r \cos \theta, r\sin\theta) e^{-i n \theta} \dd \theta\,,\\
v_{\theta,n} (r) & \, := \,  \frac{1}{2\pi} \int_0^{2\pi} v_\theta (r \cos \theta, r\sin\theta) e^{-i n \theta} \dd \theta\,.
\end{align*}
We also set  for $m\in \N\cup \{0\}$,
\begin{align}
\mathcal{Q}_m v := \sum_{|n|=m+1}^\infty \mathcal{P}_n v\,.\label{def.Q_m}
\end{align}
For notational convenience we will often write $v_n$ for $\mathcal{P}_n v$. 
Each $\mathcal{P}_n$ is an orthogonal projection in $L^2 (\Omega)^2$, and the space $L^2_\sigma (\Omega):=\overline{\{f\in C_0^\infty (\Omega)^2~|~{\rm div}\, f=0\}}^{L^2(\Omega)^2}$ is invariant under the action of $\mathcal{P}_n$.
Note that $v_0:=\mathcal{P}_0 v$ is the radial part of $v$, and thus, $\mathcal{Q}_0 v$ is the nonradial part of $v$. We will set $\mathcal{P}_n L^2 (\Omega)^2 := \{f\in L^2 (\Omega)^2~|~f=\mathcal{P}_n f\}$, and   similar notation will be used for $L^2_\sigma (\Omega)$ and $\mathcal{Q}_0$. 
A vector field $f$ in $\Omega$ is formally identified with the pair~$(\mathcal{P}_0f, \mathcal{Q}_0 f)$.
Then, for the class of   external forces we introduce the product space 
\begin{align}
Y := \mathcal{P}_0 L^1 (\Omega)^2 \times \mathcal{Q}_0 L^2 (\Omega)^2\,.\label{def.Y}
\end{align}
For the class of solutions we set 
\begin{align}
X := \mathcal{P}_0 W_0^{1,\infty} (\Omega)^2 \times \mathcal{Q}_0 W^{1,2}_0 (\Omega)^2\,.\label{def.X}
\end{align}
Here $W^{1,r}_0 (\Omega):=\{f\in W^{1,r} (\Omega)~|~f=0 ~ \text{on}~\partial\Omega\}$ for $1< r\leq \infty$.
Our first result is stated as follows.
\begin{theorem}\label{thm.main.1} There exists $\gamma>0$ such that the following statements hold.

\noindent {\rm (i)} Let $0<|\alpha|<1$. Then for any external force $f=(\mathcal{P}_0 f, \mathcal{Q}_0 f)\in Y$ satisfying 
\begin{align}
\|(\mathcal{P}_0 f)_{\theta}\|_{L^1(\Omega)} \leq \gamma |\alpha|\,, \qquad 
\| \mathcal{Q}_0 f\|_{L^2 (\Omega)}\leq \gamma |\alpha|^{2} \,,
\end{align}
there exists a solution $v\in X\cap L^\infty (\Omega)^2 \cap W^{2,1}_{loc} (\overline{\Omega})^2$ to \eqref{tildeNS_alpha} with a suitable pressure~$q \in W^{1,1}_{loc}(\overline{\Omega})$, and $v$ satisfies 
\begin{align}
\|\mathcal{P}_0 v\|_{L^\infty (\Omega)} 
+ \big\| \nabla \mathcal{P}_0 v\big\|_{L^\infty (\Omega)} 
& \leq
C \|(\mathcal{P}_0 f)_{\theta}\|_{L^1(\Omega)}
+ \frac{C}{|\alpha|^{\frac32}} \|\mathcal{Q}_0 f\|_{L^2(\Omega)}^{2}
\,,\label{est.thm.main.1.1} \\
\| \mathcal{Q}_0 v \|_{L^2 (\Omega)} 
& \leq \frac{C}{|\alpha|} \|\mathcal{Q}_{0} f\|_{L^{2} (\Omega)}
\,,\label{est.thm.main.1.2}\\
\sum_{|n|\geq 1} \| \mathcal{P}_n v \|_{L^\infty (\Omega)} 
&  \leq \frac{C}{|\alpha|^\frac34} \|\mathcal{Q}_{0} f\|_{L^{2} (\Omega)}
\,,\label{est.thm.main.1.2'}\\
\| \nabla \mathcal{Q}_0 v \|_{L^2 (\Omega)} 
& \leq \frac{C}{|\alpha|^{\frac12}} \|\mathcal{Q}_{0} f\|_{L^{2} (\Omega)}
\,.\label{est.thm.main.1.3}
\end{align}
This solution is unique in a suitable closed convex set in $X$ (see Subsection \ref{subsec.nonlinear.general.alpha} for the precise description).

\noindent {\rm (ii)} Let $|\alpha|\geq 1$. Then for any external force $f=(\mathcal{P}_0 f, \mathcal{Q}_0 f)\in Y$ satisfying 
\begin{align}
\| (\mathcal{P}_0 f)_{\theta}\|_{L^1(\Omega)} \leq \gamma \,, \qquad 
\| \mathcal{Q}_0 f\|_{L^2 (\Omega)}\leq \gamma\,,
\end{align}
there exists a solution $v\in X\cap L^\infty (\Omega)^2 \cap W^{2,1}_{loc} (\overline{\Omega})^2$ to \eqref{tildeNS_alpha} with a suitable pressure $q\in W^{1,1}_{loc}(\overline{\Omega})$, and $v$ satisfies 
\begin{align}
\| \mathcal{P}_0 v\|_{L^\infty (\Omega)} 
+  \big\| \nabla \mathcal{P}_0 v\big\|_{L^\infty (\Omega)} 
& \leq
C \|(\mathcal{P}_0 f)_{\theta}\|_{L^1(\Omega)}
+ \frac{C}{|\alpha|^{\frac12}} \|\mathcal{Q}_0 f\|_{L^2(\Omega)}^{2}
\,,\label{est.thm.main.1.4}\\
\| \mathcal{Q}_0 v \|_{L^2 (\Omega)} 
& \leq 
\frac{C}{|\alpha|^{\frac12}} \|\mathcal{Q}_{0} f\|_{L^{2} (\Omega)}
\,,\label{est.thm.main.1.5}\\
  \sum_{|n|\geq 1} \| \mathcal{P}_n v \|_{L^\infty (\Omega)} 
&  \leq \frac{C}{|\alpha|^\frac14} \|\mathcal{Q}_{0} f\|_{L^{2} (\Omega)}
\,,\label{est.thm.main.1.5'}\\
\| \nabla \mathcal{Q}_0 v \|_{L^2 (\Omega)} 
& \leq C \|\mathcal{Q}_{0} f\|_{L^{2} (\Omega)}
\,.\label{est.thm.main.1.6}
\end{align}
This solution is unique in a suitable subset of $X$.
\end{theorem}
Note that the summability of $f$ assumed in Theorem \ref{thm.main.1} is much weaker than the scale-critical one.  For the radial part $\mathcal{P}_0 v$ we can show $\lim_{|x|\rightarrow \infty} |\mathcal{P}_0 v (x)|=0$ but there is no rate in general under the assumptions of Theorem \ref{thm.main.1}.
Theorem \ref{thm.main.1} already exhibits the axisymmetrizing effect of  the fast rotating obstacle in $L^2$ and $L^\infty$, which will be further extended in Theorems \ref{thm.main.2} and \ref{thm.main.3} below.
The proof of Theorem \ref{thm.main.1} consists in two ingredients: the analysis of the linearized problem \eqref{S_alpha} (defined and studied in Section~\ref{sec.linear}), and the estimate of the interaction between the radial part and the nonradial part in the nonlinear problem (see Section~\ref{sec.nonlinear}).
The linear result used in Theorem \ref{thm.main.1} is stated in Proposition \ref{prop.all.alpha}, and the proof is based on an energy method. Although the proof of the linear result is not so difficult, there is a key observation for the term $-\alpha ( x^\bot \cdot \nabla v - v^\bot ) + \alpha U^\bot {\rm rot}\, v$.
Indeed, for the linearized problem \eqref{S_alpha} the energy computation for $v_n=\mathcal{P}_n v$ with $n\ne 0$ gives the key identity 
\begin{align}\label{proof.apriori.middle.intro.1}
\begin{split}
& \alpha n  \bigg ( \| v_{r,n} \|_{L^2(\Omega)}^2 - (1-\frac{2}{n^2})\big \| \frac{v_{r,n}}{r}\big\|_{L^2 (\Omega)}^2 +  \| v_{\theta,n} \|_{L^2 (\Omega)}^2 - \big\| \frac{v_{\theta,n}}{r}\big\|_{L^2 (\Omega)}^2 \bigg ) \\
& \qquad = - \Im \langle f_n\,, v_n \rangle_{L^2(\Omega)}\,.
\end{split}
\end{align}
Here $f_n$ denotes $\mathcal{P}_n f$ and the norm $\| g\|_{L^2 (\Omega)}$ for the function $g:[1,\infty)\rightarrow \C$ is defined as~$(2\pi)^\frac12\|g\|_{L^2((1,\infty); r  \dd r)}$.
The key point here is that the bracket in \eqref{proof.apriori.middle.intro.1} is nonnegative and provides a bound for $\|\frac{\sqrt{|x|^2-1}}{|x|} v_n \|_{L^2 (\Omega)}^2$ since $\Omega=\{|x|>1\}$. 
Then by combining with an interpolation inequality of the form 
\begin{align}\label{proof.apriori.middle.intro.2}
\| g \|_{L^2(\Omega)} \leq C \| \partial_r g \|_{L^2 (\Omega)}^\frac13 \big\| \frac{\sqrt{r^2-1}}{r} g \big\|_{L^2 (\Omega)}^\frac23 + C\big \| \frac{\sqrt{r^2-1}}{r} g\big \|_{L^2 (\Omega)}
\end{align}
for any scalar function $g\in W^{1,2}((1,\infty); r\dd r)$ and the dissipation from the Laplacian in the energy computation, we can close the energy estimate for all $\alpha\ne 0$. The proof of \eqref{proof.apriori.middle.intro.2} is given in Appendix~\ref{appendix.interpolation}.
In solving the nonlinear problem the key observation is that the product of the radial parts in the nonlinear term can always be written in a gradient form and thus regarded as a pressure term, which yields the identity
\begin{align}\label{proof.nonlinear.intro.1}
v^\bot {\rm rot}\ v=  v_0^\bot {\rm rot}\, \mathcal{Q}_0 v  + (\mathcal{Q}_0 v) ^\bot {\rm rot}\, v_0 + (\mathcal{Q}_0 v)^\bot {\rm rot}\, \mathcal{Q}_0 v +\nabla \tilde q 
\end{align}
for a suitable $\tilde q$. Since $\mathcal{P}_0 \big (v_0^\bot {\rm rot}\, \mathcal{Q}_0 v  + (\mathcal{Q}_0 v) ^\bot {\rm rot}\, v_0 \big ) = 0$ as long as $\mathcal{Q}_0 v\in W^{1,2}_0 (\Omega)^2$ the radial part of the velocity in the right-hand side of \eqref{proof.nonlinear.intro.1} (neglecting $\nabla \tilde q$) belongs to~$L^1 (\Omega)^2$, which is the same summability as the space $Y$. This is a brief explanation for the reason why we can close the nonlinear estimate and solve \eqref{tildeNS_alpha} in $X$ for a source $f\in Y$.
 
Our second result is focused on the fast rotation case $|\alpha|\gg 1$.
 In this regime there are three fundamental mechanisms in our system:
 
\medskip

\noindent  (I) an axisymmetrization due to the fast rotation of the obstacle, 
  
\noindent  (II) the presence of a boundary layer for the nonradial part of the flow due to the noslip boundary condition, 
  
\noindent  (III) the diffusion in high angular frequencies due to the viscosity.

\medskip

 (I) and (II) are potentially in a competitive relation,  for the noslip boundary condition and the boundary layer can suppress the effect of the fast rotation to some extent. 

(II) and (III) are also competitive. Indeed, it is natural that if the viscosity is strong enough then the boundary layer  is diffused and is no longer observable.
The important task here is to determine the regime of   angular frequencies in which the boundary layer appears, and to estimate the thickness of the boundary layer. We   show that the boundary layer appears in the regime $1\leq |n|\ll O(|\alpha|^\frac12)$, and the thickness of the boundary layer is~$(2|\alpha n|)^{-\frac13}$ for each $n$ in this regime. In constructing the boundary layer the term~$\alpha U^\bot {\rm rot}\,v$ plays a crucial role as well as the term $-\alpha ( x^\bot \cdot \nabla v - v^\bot )$. 
In fact, if we drop the term $\alpha U^\bot {\rm rot}\, v$ as a model problem then the thickness of the boundary layer arising from the rotation term $-\alpha ( x^\bot \cdot \nabla v - v^\bot )$ is $|\alpha n|^{-\frac12}$, and the leading boundary layer profile is simply described by   exponential functions.
The term $\alpha U^\bot {\rm rot}\, v$ leads to a significant change both in the thickness and in the profile of the boundary layer,
and we need to introduce the Airy function to describe the profile of the boundary layer associated with the term $-\alpha ( x^\bot \cdot \nabla v - v^\bot )+ \alpha U^\bot {\rm rot}\, v$.

By performing the boundary layer analysis we can improve the result  stated in (ii) of Theorem \ref{thm.main.1} in the regime $|\alpha|\gg 1$, which is briefly described as follows.
\begin{theorem}\label{thm.main.2} There exists $\gamma>0$ such that the following statement holds. For all sufficiently large $|\alpha| \ge 1$ and for any external force  $f=(\mathcal{P}_0 f, \mathcal{Q}_0 f)\in Y$ satisfying 
\begin{align}
\| (\mathcal{P}_0 f)_{\theta}\|_{L^1(\Omega)} \leq \gamma |\alpha|^{\frac13}\,, \qquad \| \mathcal{Q}_0 f\|_{L^2 (\Omega)}\leq \gamma |\alpha|^{\frac13} \,,\label{assumption.thm.main.2}
\end{align}
there exists a solution $v\in X\cap L^\infty (\Omega)^2 \cap W^{2,1}_{loc} (\overline{\Omega})^2$ to \eqref{tildeNS_alpha} with a suitable pressure $q\in W^{1,1}_{loc}(\overline{\Omega})$, and $v$ satisfies 
\begin{align}
\| \mathcal{P}_0 v\|_{L^\infty (\Omega)} 
+ \| \nabla \mathcal{P}_0 v\|_{L^\infty (\Omega)}
& \leq 
C \|(\mathcal{P}_0 f)_{\theta}\|_{L^1(\Omega)}
+ \frac{C}{|\alpha|} \|\mathcal{Q}_0 f\|_{L^2(\Omega)}^{2}
\,,\label{est.thm.main.2.1}\\
\| \mathcal{Q}_0 v \|_{L^2 (\Omega)} & \leq \frac{C}{|\alpha|^\frac23} \| \mathcal{Q}_0 f \|_{L^2 (\Omega)} \,,\label{est.thm.main.2.2}\\
 \sum_{|n|\geq 1} \| \mathcal{P}_n v \|_{L^\infty (\Omega)} & \leq \frac{C\big(\log |\alpha| \big)^\frac12}{|\alpha|^\frac12} \| \mathcal{Q}_0 f \|_{L^2 (\Omega)} \,,\label{est.thm.main.2.2'}\\
\| \nabla \mathcal{Q}_0 v \|_{L^2 (\Omega)} & \leq \frac{C}{|\alpha|^{\frac13}} \| \mathcal{Q}_0 f\|_{L^2(\Omega)}\label{est.thm.main.2.3}\,.
\end{align}
This solution is unique in a suitable subset of $X$.
\end{theorem}
By fixing the external force $f$ we can state 
Theorem \ref{thm.main.2} 
in a different but more convenient way to understand the qualitative behavior of solutions in the fast rotation limit.
\begin{theorem}\label{thm.main.3}  For any $f=(\mathcal{P}_0 f, \mathcal{Q}_0 f)\in Y$ there is $\alpha_0=\alpha_0 (\|f\|_Y)\geq 1$ such that the following statements hold. If $|\alpha|\geq \alpha_0$ then
there exists a solution $v=v^{(\alpha)}\in X\cap W^{2,1}_{loc} (\overline{\Omega})^2$ to \eqref{tildeNS_alpha} with a suitable pressure $q^{(\alpha)} \in W^{1,1}_{loc}(\overline{\Omega})$, and $v^{(\alpha)}$ satisfies 
\begin{align}\label{est.thm.main.3.1}
\| v^{(\alpha)} -  v_0^{{\rm linear}}\|_{L^\infty (\Omega)} & \leq \frac{C (\log |\alpha| )^\frac12}{|\alpha|^\frac12}\,\cdotp
\end{align}
Here $v_0^{{\rm linear}}$ is the solution to the linearized problem \eqref{S_alpha} defined in page~{\rm\pageref{Salpha}} with $f$ replaced by $\mathcal{P}_0f$ which is in fact independent of $\alpha$, and $C$ depends only on $\|f\|_Y$.
Moreover, there exists $\kappa>0$ independent of $\alpha$ and $f$ such that if $1\leq |n| \leq \kappa |\alpha|^\frac12$ then $v_n^{(\alpha)}=\mathcal{P}_n v^{(\alpha)}$ is written in the form
\begin{align}\label{decom.thm.main.3} 
v_n^{(\alpha)}=v_n^{(\alpha),{\rm slip}} +v_n^{(\alpha),{\rm slow}} + v_{n,{\rm BL}}^{(\alpha)} + \widetilde v_n^{(\alpha)}.
\end{align}
Here $v_n^{(\alpha),{\rm slip}}$ satisfies $v_{n,r}^{(\alpha),{\rm slip}}={\rm rot}\, v_n^{(\alpha),{\rm slip}}=0$ on $\partial\Omega$,  $v_n^{(\alpha),{\rm slow}}$ is irrotational in $\Omega$, and~$v_{n,{\rm BL}}^{(\alpha)}$ possesses a boundary layer structure with   the boundary layer thickness $|2\alpha n|^{-\frac13}$. Finally the following estimates hold:
$$
\|v_n^{(\alpha),{\rm slip}}\|_{L^2 (\Omega)} +\|v_n^{(\alpha),{\rm slow}}\|_{L^2 (\Omega)} +\|v_{n,{\rm BL}}^{(\alpha)}\|_{L^2 (\Omega)} \leq \frac{C}{|\alpha n|^\frac23}\,,
$$
while $\widetilde v_n^{(\alpha)}$ is a remainder which satisfies 
$$\| \widetilde v_n^{(\alpha)} \|_{L^2 (\Omega)} \leq \frac{C}{|\alpha n| }\,\cdotp$$
Here the constant $C$ depends only on $\|f\|_Y$.
\end{theorem}

By going back to \eqref{NS_alpha.intro}, Theorems \ref{thm.main.2} and \ref{thm.main.3} show that there exists a unique solution~$u=u^{(\alpha)}$ which satisfies 
\begin{align}\label{expansion.u}  
\| u^{(\alpha)} - \alpha U - v_0^{{\rm linear}} \|_{L^\infty (\Omega)} \leq \frac{C (\log |\alpha| )^\frac12}{|\alpha|^\frac12}\,, \qquad |\alpha|\gg 1\,.
\end{align}
The expansion \eqref{expansion.u} verifies the axisymmetrizing effect (measured in $L^\infty$) due to the fast rotation. The logarithmic factor  $(\log |\alpha| )^\frac12$ is simply due to the regularity of $f$, and if $f$ has more regularity such as $\sum_{n\ne 0} \| \mathcal{P}_n f\|_{L^2 (\Omega)}^s<\infty$ for some $s<2$ then the factor $(\log |\alpha|)^\frac12$ in  \eqref{est.thm.main.3.1} and \eqref{expansion.u} can be dropped.
Moreover, the power $|\alpha|^{-\frac12}$ can be also improved by assuming enough regularity of $f$. For example, if $\mathcal{Q}_0f\in W^{1,2}_0 (\Omega)^2$ in addition, then  $|\alpha|^{-\frac12}$ is replaced by $|\alpha|^{-\frac34}$, though we do not go into the detail on this point. 
The new ingredient of the proof of Theorems \ref{thm.main.2} and \ref{thm.main.3} is 
 stated in Proposition \ref{decompositionvelocity} and consists in  
  refined estimates for the linearized problem \eqref{S_alpha}.
The nonlinear problem is handled exactly in the same manner as in the proof of Theorem \ref{thm.main.1}.
For \eqref{S_alpha} we observe that in  polar coordinates  the  angular mode~$n$ of the streamfunction satisfies the ODE in $r \in (1,\infty)$
\begin{equation}\label{eq.fast.slow.intro}
\begin{aligned}
& \bigg(\frac{\dd^2}{\dd r^2} + \frac1r \frac{\dd}{\dd r}  -\frac{n^2}{r^2}+  i \alpha n \big (  1-\frac1{r^2}\big ) \bigg ) \big ( \frac{\dd^2 }{\dd r^2} + \frac1r\frac{\dd}{\dd r} - \frac{n^2}{r^2} \big ) \psi_n  = 0\,,
\end{aligned}
\end{equation}
with the  boundary condition $\psi_n (1)=\frac{\dd \psi_n}{\dd r} (1) =0$ when $|n|\geq 1$.
The thickness of the boundary layer originating from the fast rotation is determined by the balance between $ \frac{\dd^2}{\dd r^2}$ and $i\alpha n (1-\frac{1}{r^2})\approx 2i\alpha n (r-1)$ near $r=1$ as long as the dissipation $-\frac{n^2}{r^2}\approx -n^2$ is moderate. This implies that the thickness  is $|2\alpha n|^{-\frac13}$. Then the regime of $n$ where the dissipation is relatively moderate is estimated from the condition $n^2\ll \frac{\dd ^2}{\dd r^2} \approx O(|\alpha n|^\frac23)$, which leads to $|n|\ll O(|\alpha|^\frac12)$. From this observation we employ the boundary layer analysis when the angular frequency $n$ satisfies $1\leq |n|\ll O (|\alpha|^\frac12)$, while we just apply Proposition \ref{prop.all.alpha} in the regime $|n|\geq O(|\alpha|^\frac12)$ where the boundary layer due to the fast rotation is no longer present.
In the regime $1\leq |n|\ll O(|\alpha|^\frac12)$ we first consider \eqref{S_alpha} with $f=\mathcal{P}_n f$ but under the slip boundary condition $v_{r,n}={\rm rot}\, v_n =0$ on $\partial\Omega$. The estimate of this {\it slip} solution is obtained by the energy method for the vorticity equations thanks to the boundary condition~${\rm rot}\, v_n=0$ on $\partial\Omega$. The key point is that the term $U^\bot {\rm rot}$ in the velocity equation becomes $U\cdot \nabla$ in the vorticity equation which is antisymmetric because of ${\rm div}\, U=0$, and thus, it is easy to apply the energy method for the vorticity   under the slip boundary condition.
The {\it noslip} solution is then obtained by correcting the boundary condition. To this end we construct the boundary layer solution called the fast mode.
The leading profile of the boundary layer is given by a suitable integral of the Airy function.
In order to recover the noslip boundary condition, the fact that~$\displaystyle \int_0^\infty \Ai  (s) \dd s \ne 0$ is crucial and it plays the role of a nondegeneracy condition in our construction of the solution.
This construction gives a formula as in \eqref{decom.thm.main.3} for the solution to \eqref{S_alpha}. 
Compared with Proposition \ref{prop.all.alpha}, which is based only on an energy computation for the velocity field, the estimate of the $n$ mode $v_n$ is drastically improved for $1\leq |n|\ll O(|\alpha|^\frac12)$ thanks to the boundary layer analysis.  
On the other hand, in the regime $|n|\geq O(|\alpha|^\frac12)$,  Proposition \ref{prop.all.alpha} for the {\it noslip} solution already gives the same decay estimates as in Proposition \ref{prop.mS_alpha} for the {\it slip} solution,  as expected.

This paper is organized as follows. In Section \ref{sec.pre} we recall some basic facts on operators and   vector fields   in   polar coordinates. In Section \ref{sec.linear} the linearized problem \eqref{S_alpha} is studied. This section is the core of the paper. In Subsection \ref{subsec.general} we prove the linear estimates which are valid for all $\alpha\in \R\setminus\{0\}$. These are summarized in Proposition \ref{prop.all.alpha}. 
Subsection \ref{subsec.fast} is devoted to the linear analysis for the case $|\alpha|\gg 1$, and the main result of this section is Proposition~\ref{decompositionvelocity}. The nonlinear problem is discussed in Section \ref{sec.nonlinear}.  Some basics on the Airy function and the proof of the interpolation inequality \eqref{proof.apriori.middle.intro.2} are given in the appendix.

\section{Preliminaries}\label{sec.pre}

In this preliminary section we state basic results on some differential operators and   vector fields in   polar coordinates.

\subsection{Operators in   polar coordinates}
The following formulas will be used frequently:
\begin{align}
{\rm div}\,  v & \, = \, \partial_1 v_1 + \partial_2 v_2 \, = \,  \frac1r \partial_r ( r v_r ) + \frac1r \partial_\theta v_\theta\,, \label{polar.div}\\
{\rm rot}\, v & \, = \, \partial_1v_2 - \partial_2 v_1 \, = \,  \frac1r \partial_r ( r v_\theta ) - \frac1r \partial_\theta v_r \,,\label{polar.rot} \\
| \nabla v |^2 & \, = \,  | \partial_r v_r |^2 + |\partial_r v_\theta |^2 + \frac{1}{r^2} \big ( | \partial_\theta v_r - v_\theta |^2 + |v_r + \partial_\theta v_\theta |^2 \big )\,,\label{polar.grad}
\end{align} 
and 
\begin{align}
\begin{split}\label{polar.laplace}
-\Delta v &  \, = \, \bigg ( - \partial_r \big  ( \frac1r \partial_r (r v_r ) \big )  - \frac{1}{r^2} \partial_\theta^2 v_r + \frac{2}{r^2} \partial_\theta v_\theta \bigg ) {\bf e}_r \\
& ~~~~~ +  \bigg ( - \partial_r \big ( \frac1r \partial_r (r v_\theta ) \big ) - \frac{1}{r^2} \partial_\theta^2 v_\theta -  \frac{2}{r^2} \partial_\theta v_r \bigg ) {\bf e}_\theta \,,
\end{split}
\end{align}
\begin{align*}
\begin{split}
{\bf e}_r \cdot \nabla v & \, = \, (\partial_r v_r ) \, {\bf e}_r \, + \, (\partial_r v_\theta) \,  {\bf e}_\theta\,,\\
{\bf e}_\theta \cdot \nabla v & \, = \, \frac{\partial_\theta v_r - v_\theta}{r}\,  {\bf e}_r  \, + \, \frac{\partial_\theta v_\theta + v_r}{r} \, {\bf e}_\theta\,. 
\end{split}
\end{align*}
In particular, we have 
\begin{align}\label{polar.e_theta}
x^\bot \cdot \nabla v - v^\bot & = |x| \big ( {\bf e}_\theta \cdot \nabla v \big ) - \big ( v_r {\bf e}_r^\bot + v_\theta {\bf e}_\theta^\bot \big ) \nonumber \\
& = (\partial_\theta v_r - v_\theta )\,  {\bf e}_r  \, + \, (\partial_\theta v_\theta + v_r) \, {\bf e}_\theta - \big ( v_r {\bf e}_r^\bot + v_\theta {\bf e}_\theta^\bot \big ) \nonumber \\
& = \partial_\theta v_r \,  {\bf e}_r + \partial_\theta v_\theta \, {\bf e}_\theta\,.
\end{align} 

\vspace{0.3cm}

\noindent From \eqref{polar.grad} and the definition of $\mathcal{P}_n$ in \eqref{def.P_n} it follows that
for $n\in \N \cup\{0\}$ and for $v$ in~$W^{1,2}(\Omega)^2$,
\begin{align*}
\| \nabla v \|_{L^2 (\Omega)}^2 & = \sum_{n\in \Z} \| \nabla \mathcal{P}_n v \|_{L^2 (\Omega)}^2\,, 
\\
\begin{split}
|\nabla \mathcal{P}_n v |^2 & = |\partial_r v_{r,n}|^2 + \frac{1+n^2}{r^2} |v_{r,n}|^2 \\
& \qquad + |\partial _r v_{\theta,n}|^2 + \frac{1+n^2}{r^2}|v_{\theta,n}|^2  - \frac{4 n}{r^2} \Im  ( v_{\theta,n} \overline{v_{r,n}} )\,. 
\end{split}
\end{align*}
In particular, we have 
\begin{align}
| \nabla \mathcal{P}_n v |^2 \geq  |\partial_r v_{r,n}|^2 + \frac{(|n|-1)^2}{r^2} |v_{r,n}|^2 + |\partial _r v_{\theta,n}|^2 + \frac{(|n|-1)^2}{r^2}|v_{\theta,n}|^2 \,,\label{polar.grad.n'} 
\end{align}
and thus, from the definition of $\mathcal{Q}_m$ in \eqref{def.Q_m},
\begin{align*}
\| \nabla \mathcal{Q}_m v\|_{L^2(\Omega)}^2 \geq  \| \partial_r (\mathcal{Q}_m v)_r \|_{L^2(\Omega)}^2 +  \| \partial_r (\mathcal{Q}_m v)_\theta \|_{L^2(\Omega)}^2 + m^2 \| \frac{v}{|x|} \|_{L^2(\Omega)}^2\,.
\end{align*}

\subsection{The Biot-Savart law in polar coordinates}

For a given scalar field $\omega$ in $\Omega$,
the streamfunction $\psi$ is formally defined as the solution to the Poisson equation:  $-\Delta \psi = \omega$ in $\Omega$,~ $\psi =0$ on $\partial \Omega$. For $n\in \Z$ and $\omega\in L^2(\Omega)$ we set 
\begin{align}
\begin{split}
\mathcal{P}_n \omega & := \frac{1}{2\pi} \int_0^{2\pi} \omega (r \cos s, r\sin s) e^{- i n s} \dd s  \, e^{i n \theta} \,,\\
\omega_n  & := \big ( \mathcal{P}_n \omega \big ) e^{-i n \theta}\,.\label{def.w_n}
\end{split}
\end{align}
By using  the Laplace operator in   polar coordinates, the Poisson equation for the Fourier mode $n$
is given by
\begin{align}
H_n \psi_n \, := \, - \psi_n^{''}  - \frac{1}{r} \psi_n^{'} + \frac{n^2}{r^2} \psi_n = \omega_n\,,~~~ r>1\,,~~~~~~~ \psi_n (1) =0\,.\label{eq.stream}
\end{align}
Let $|n|\geq 1$. Then the solution $\psi_n= \psi_n[\omega_n]$ to the ordinary differential equation \eqref{eq.stream} decaying at spatial infinity is formally given as
\begin{align*}
\begin{split}
\psi_n[\omega_n] (r) & =  \frac{1}{2 |n|}  \bigg (- \frac{d_n [\omega_n] }{r^{|n|} } + \frac{1}{r^{|n|}} \int_1^r s^{1+|n|} \omega_n (s) \dd s   + r^{|n|}\int_r^\infty s^{1-|n|} \omega_n (s) \dd s \bigg ) \,,\\
d_n [\omega_n]  & :=   \int_1^\infty s^{1-|n|} \omega_n (s) \dd s\,.
\end{split}
\end{align*}
The Biot-Savart law $V_n[\omega_n]$ is then written as
\begin{align}\label{def.V_n}
\begin{split}
&V_n [\omega_n] \, := \, V_{r,n} [\omega_n] \, e^{i n \theta}  \,  {\bf e}_r \, + \, V_{\theta,n} [\omega_n]  \, e^{i n\theta}\, {\bf e}_\theta\,,\\
& V_{r,n} [\omega_n] \, := \, \frac{i n}{r} \psi_n [\omega_n] \,,~~~~~~~~ V_{\theta,n} [\omega_n] \, := \, - \frac{\dd   }{\dd r} \psi_n [\omega_n] \,.
\end{split}
\end{align}
The velocity $V_n[\omega_n]$ is well defined at least when $r^{1-|n|} \omega_n\in L^1 ((1,\infty))$, and it is straightforward to see that
\begin{align}\label{V_n}
\begin{split}
&{\rm div}\, V_n [\omega_n] \, =\, 0\,,~~~~~{\rm rot}\, V_n [\omega_n] \, = \, \omega_n \,  e^{i n \theta}~~~~~~~{\rm in}~~\Omega\,,\\
& {\bf e}_r \cdot V_n [\omega_n] \, =\,  0 ~~~~~{\rm on}~\partial\Omega\,.
\end{split}
\end{align}
 The condition $r^{1-|n|} \omega_n\in L^1 ((1,\infty))$ is automatically satisfied when $\omega\in L^2 (\Omega)$ and~$|n|\geq 2$. When $|n|=1$ the integral in the definition of $\psi_n[\omega_n]$ does not always converge absolutely for general $\omega\in L^2 (\Omega)$. However, it is well-defined if $\omega={\rm rot}\, u$ for some $u\in W^{1,2}(\Omega)^2$, for  one can apply the integration by parts that ensures the convergence of  $\displaystyle \lim_{N\rightarrow \infty} \int_r^N \omega_n \dd r$ even when $|n|=1$. 
As a result, we can check that 
for any solenoidal vector field $v$ in~$ L^2_\sigma (\Omega) \cap W^{1,2} (\Omega)^2$, the $n$ mode $v_n = \mathcal{P}_n v$ is expressed in terms of its vorticity $\omega_n$ by the formula \eqref{def.V_n} when $|n|\geq 1$.

\section{Analysis of the linearized system}\label{sec.linear}

The linearized system around~$\alpha U$ for \eqref{tildeNS_alpha} is  \label{Salpha}
\begin{equation}\tag{${\rm S}_{\alpha}$}\label{S_alpha}
  \left\{
\begin{aligned}
  -\Delta v - \alpha ( x^\bot \cdot \nabla v - v^\bot ) + \nabla q + \alpha U^\bot {\rm rot}\, v & \,=\, 
  f \,,  \qquad \qquad x\in \Omega \,,\\
{\rm div}\, v &  \,=\, 0\,,   \qquad \qquad  x \in \Omega\,, \\
  v & \,=\, 0 \,, \qquad \qquad  x \in \partial \Omega\,. \\
\end{aligned}\right.
\end{equation}
In this section we study \eqref{S_alpha} for $\alpha\in \R\setminus\{0\}$. 

\subsection{General estimate}\label{subsec.general}

In this subsection we establish estimates on solutions to \eqref{S_alpha} that are valid for all $\alpha\ne 0$.  For convenience we set for scalar functions $g, h:[1,\infty) \rightarrow \C$,
\begin{align*}
\langle g, h \rangle_{L^2 (\Omega)} : = 2\pi \int_1^\infty g (r) \overline{h(r)} r\dd r\,, \qquad \| g\|_{L^2(\Omega)}^2: = 2\pi \int_1^\infty |g(r)|^2 r \dd r\,.
\end{align*}
Before going into   details let us give a remark on the verification of the energy argument.
Let us assume that $f\in L^2(\Omega)^2$ and $v\in W_{0}^{1,2}(\Omega)^2 \cap W^{2,2}_{loc} (\overline{\Omega})^2$ is a solution to \eqref{S_alpha} for some $q\in W^{1,2}_{loc}(\overline{\Omega})$.
We have to be careful  when applying the energy argument due to the presence of the term $x^\bot \cdot \nabla v$ in the first equation of \eqref{S_alpha}, for this term involves a linearly growing coefficient, and therefore it is not clear whether the inner product $\langle x^\bot \cdot \nabla v, v\rangle_{L^2(\Omega)}$ makes sense or not. A similar difficulty appears in taking the inner product $\langle \nabla q, v\rangle_{L^2(\Omega)}$, since we are assuming only $q\in W^{1,2}_{loc} (\overline{\Omega})$. 
The most convenient way to overcome this difficulty is to consider the equation for $v_n := \mathcal{P}_n v$, which is identified with $(v_{r,n},v_{\theta,n})$. 
Note that $v_{r,n}, v_{\theta,n}\in W^{1,2}_0 ((1,\infty); r \dd r)\cap W^{2,2}_{loc} ([1,\infty))$ if $v\in W_{0}^{1,2}(\Omega)^2 \cap W^{2,2}_{loc} (\overline{\Omega})^2$.
Let us denote by $\omega_n=\omega_n (r)$ the $n$ mode of the vorticity of $v$ in the polar coordinates, i.e., $\omega_n (r)  = ({\rm rot} \, \mathcal{P}_n v )e^{-i n \theta}$. Similarly, we set $q_n = q_n (r) = (\mathcal{P}_n q)e^{-in\theta}$, where the projection $\mathcal{P}_n$ for the scalar $q$ is defined as
$$\mathcal{P}_n q:=\frac{1}{2\pi} \int_0^{2\pi} q(r\cos s, r\sin s) e^{-in s} \dd s \, e^{in\theta}\,.$$
Then, from \eqref{polar.laplace} and \eqref{polar.e_theta}, $v_{r,n}$ and $v_{\theta,n}$ obey the following equations:
\begin{align}\label{eq.v_rn}
\begin{split}
 - \partial_r \big  ( \frac1r \partial_r (r v_{r,n} ) \big )  + \frac{n^2}{r^2}  v_{r,n} + \frac{2in}{r^2} v_{\theta,n}  -i \alpha n v_{r,n} - \alpha \frac{\omega_n}{r} + \partial_r q_n = f_{r,n}\,,
\end{split}
\end{align} 
\begin{align}\label{eq.v_thetan}
\begin{split}
 - \partial_r \big  ( \frac1r \partial_r (r v_{\theta,n} ) \big )  + \frac{n^2}{r^2}  v_{\theta,n} - \frac{2in}{r^2} v_{r,n}  -i \alpha n v_{\theta,n}  + i n q_n = f_{\theta,n}\,,
\end{split}
\end{align} 
together with the divergence free condition $\partial_r (r v_{r,n}) + in v_{\theta,n}=0$ and the noslip boundary condition $v_{r,n}(1)=v_{\theta,n}(1)=0$. Then the key observation   is that the factor $-i\alpha n$ is now regarded as a resolvent parameter, and by setting $\lambda:= -i\alpha n$, the above system is  equivalent to 
\begin{align}\label{eq.reduced.stokes}
(\lambda - \Delta ) v_n  + \nabla \mathcal{P}_n q + \alpha U^\bot {\rm rot}\, v_n = f_n\,, \qquad x\in  \Omega\,,
\end{align}
with ${\rm div}\, v_n=0$ and $v_n|_{\partial\Omega}=0$, where $f_n = \mathcal{P}_n f\in \mathcal{P}_n L^2 (\Omega)^2$.
Indeed,   system \eqref{eq.reduced.stokes} in   polar coordinates is exactly \eqref{eq.v_rn} and \eqref{eq.v_thetan}. The key point is that there is no term involving a linearly growing coefficient in \eqref{eq.reduced.stokes}, and therefore we can apply the standard regularity theory of the Stokes resolvent system with a resolvent parameter~$\lambda$. 
Let us assume that~$n\ne 0$. Then $\lambda\ne 0$ since we are assuming that $\alpha\ne 0$.
If $v_n \in L^2_\sigma (\Omega)\cap W^{1,2}_{0} (\Omega)^2 \cap W^{2,2}_{loc}(\overline{\Omega})^2$ and~$\mathcal{P}_n q \in W^{1,2}_{loc} (\overline{\Omega})$ is a solution to \eqref{eq.reduced.stokes}, then $(v_n, \mathcal{P}_n q)$ is a weak solution to the Stokes system with   source term $f_n - \alpha U^\bot {\rm rot}\, v_n$ which clearly belongs to $L^2 (\Omega)^2$, and thus, the regularity theory of the Stokes system implies that $v_n \in W^{2,2}(\Omega)^2$ and $\nabla \mathcal{P}_n q\in L^2 (\Omega)^2$.
In this way we can recover the summability of $\nabla^2 v_n, \nabla \mathcal{P}_n q\in L^2 (\Omega)^2$.
Then, by going back to the system \eqref{eq.v_rn} and \eqref{eq.v_thetan}, we also find that $q_n \in L^2 ((1,\infty); r\dd r)$ from \eqref{eq.v_thetan}, for all the other terms in \eqref{eq.v_thetan} belong to $L^2([1,\infty); r\dd r)$. As a summary, for any solution~$v\in L^2_\sigma (\Omega)\cap W^{1,2}_{0} (\Omega)^2 \cap W^{2,2}_{loc}(\overline{\Omega})^2$ and $q\in W^{1,2}_{loc}(\overline{\Omega})$ to \eqref{S_alpha}, we can rigorously verify the energy computation for the system \eqref{eq.v_rn}-\eqref{eq.v_thetan} in each $n$ mode $(v_{r,n},v_{\theta,n})$ with $n \ne 0$. 
The estimate for the $0$ mode is handled in a different way from the energy method, and is discussed in Subsection \ref{subsubsec.0} below.

\medskip
\noindent 
Our main result in this subsection is stated as follows. 
Let us recall that the projection~$\mathcal{Q}_0$ is defined as $\mathcal{Q}_0: = \sum_{n\ne 0} \mathcal{P}_n$.

\begin{proposition}\label{prop.all.alpha} Let $\alpha\in \R\setminus\{0\}$.

\noindent {\rm (i)} For any external force $f_0 \in \mathcal{P}_0 L^1 (\Omega)^2$ the system \eqref{S_alpha} admits a unique solution $(v_0,\nabla q)$ with $v_0\in  \mathcal{P}_0 L^\infty (\Omega)^2 \cap W_0^{1,\infty} (\Omega)^2$, $\nabla^2 v_0\in L^1_{loc} (\overline{\Omega})^{2\times 2}$, $q\in W^{1,1}_{loc}(\overline{\Omega})$. Moreover, $v_0 =v_{\theta,0} {\bf e}_\theta$ and 
\begin{align}\label{est.prop.all.alpha.0}
\| v_0 \|_{L^\infty (\Omega)} + \big\| |x|\nabla v_0 \big\|_{L^\infty (\Omega)} & \leq C \| f_{\theta,0} \|_{L^1 (\Omega)}\,.
\end{align}
Here $f_{\theta,0} := f_0 \cdot {\bf e}_\theta$ and $C$ is independent of $\alpha$. 

\noindent {\rm (ii)} For any external force $f\in \mathcal{Q}_0 L^2(\Omega)^2$ the system \eqref{S_alpha} admits a unique solution $(v, \nabla q)$ with $v\in \mathcal{Q}_0 L^2_\sigma (\Omega)\cap W^{1,2}_{0} (\Omega)^2 \cap W^{2,2}_{loc}(\overline{\Omega})^2$ and $q\in W^{1,2}_{loc} (\overline{\Omega})$. 
Moreover, $v_n = \mathcal{P}_n v$ satisfies the following estimates: if $1\leq |n| < 1+ \sqrt{2|\alpha|}$ then 
\begin{align}
\| v_{n} \|_{L^2 (\Omega)} & \leq \frac{C}{|\alpha|^\frac12} \big ( \frac{1}{|n|} + \frac{1}{|\alpha|^\frac12} \big ) \| f_n \|_{L^2 (\Omega)}\,, \label{est.prop.all.alpha.1}\\
\big\| \frac{\sqrt{|x|^2-1}}{|x|} v_{n} \big\|_{L^2 (\Omega)}  & \leq \frac{C}{|\alpha n|^\frac12|\alpha|^\frac14} \big ( \frac{1}{|n|} + \frac{1}{|\alpha|^\frac12} \big )^\frac12  \| f_n \|_{L^2(\Omega)}\,, \label{est.prop.all.alpha.2} \\
\| v_{n} \|_{L^\infty (\Omega)} & \leq \frac{C}{|\alpha|^\frac14} \big ( \frac{1}{|n|} + \frac{1}{|\alpha|^\frac12} \big ) \| f_n \|_{L^2 (\Omega)}\,, \label{est.prop.all.alpha.1'}\\
\| \nabla v_n  \|_{L^2 (\Omega)} & \leq C \big ( \frac{1}{|n|} + \frac{1}{|\alpha|^\frac12} \big ) \| f_n \|_{L^2(\Omega)}\,,\label{est.prop.all.alpha.3} 
\end{align}
while if $|n|\geq  1+ \sqrt{2|\alpha|}$ then 
\begin{align}
\| v_n \|_{L^2 (\Omega)} & \leq \frac{C}{|n|} \big ( \frac{1}{|n|}  + \frac{1}{|\alpha|} \big ) \| f_n \|_{L^2 (\Omega)}\,, \label{est.prop.all.alpha.4} \\
\big\| \frac{\sqrt{|x|^2-1}}{|x|} v_{n} \big\|_{L^2 (\Omega)} &\leq \frac{C}{|\alpha n|^\frac12|n|^\frac12}   \big ( \frac{1}{|n|} + \frac{1}{|\alpha|} \big )^\frac12 \| f_n \|_{L^2(\Omega)}\,,\label{est.prop.all.alpha.5} \\
\| v_n \|_{L^\infty (\Omega)} & \leq \frac{C}{|n|^\frac34} \big ( \frac{1}{|n|}  + \frac{1}{|\alpha|} \big )^\frac34 \| f_n \|_{L^2 (\Omega)}\,, \label{est.prop.all.alpha.4'} \\
\| \nabla v_n \|_{L^2 (\Omega)}  & \leq  \frac{C}{|n|^\frac12} \big ( \frac{1}{|n|} + \frac{1}{|\alpha|} \big )^\frac12 \| f_n \|_{L^2(\Omega)} \,.\label{est.prop.all.alpha.6} 
\end{align}
Finally if~$|\alpha|\gg1$ and~$|n| = O(|\alpha|^\frac12)$ then
\begin{align}
\| v_n \|_{L^2 (\Omega)} & \leq \frac{C}{|\alpha|} \| f_n \|_{L^2 (\Omega)}\,, \label{est.prop.large.alpha.1} \\
\| v_n \|_{L^\infty (\Omega)} & \leq \frac{C}{|\alpha|^\frac34}   \| f_n \|_{L^2 (\Omega)}\,, \label{est.prop.large.alpha.2} \\
\| \nabla v_n \|_{L^2 (\Omega)}  & \leq  \frac{C}{|\alpha|^\frac12}   \| f_n \|_{L^2(\Omega)} \,.\label{est.prop.large.alpha.3} 
\end{align}
Here $f_n:=\mathcal{P}_n f$ and $C$ is independent of $n$ and $\alpha$.
\end{proposition}

\subsubsection{Structure and estimate of the  $0$ mode}\label{subsubsec.0}

Firstly we observe that if $v_0$ satisfies $v_0 \in \mathcal{P}_0 L^\infty (\Omega)^2 \cap W^{1,\infty}_0 (\Omega)^2$ then 
the divergence-free condition in   polar coordinates \eqref{polar.div} implies that 
\begin{align*}
\frac{\dd (r v_{r,0})}{\dd r} =0\,,
\end{align*}
and thus, $v_{r,0} =\displaystyle \frac{C}{r}$ with some constant $C$. Then the noslip boundary condition leads to~$C=0$, and therefore, $v_{r,0}=0$.
So it suffices to consider the angular part $v_{\theta,0}$. 
From~\eqref{polar.e_theta} we have 
\begin{align*}
x^\bot \cdot \nabla v_0^\bot - v_0^\bot =0\,,
\end{align*}
and we also note that the term $U^\bot {\rm rot}\, v_0$ in \eqref{S_alpha} with $U^\bot =\displaystyle-\frac{x}{|x|^2}$ is always written in a gradient form, so can be absorbed in a pressure term.

Collecting these remarks, we see that any solution $(v_0,\nabla q)$ to \eqref{S_alpha} with $f_0 \in \mathcal{P}_0 L^1(\Omega)^2$ satisfying $v_0\in \mathcal{P}_0 L^\infty (\Omega)^2 \cap W^{1,\infty}_0 (\Omega)^2$,~$\nabla^2 v_0 \in L^1_{loc}(\overline{\Omega})^{2\times 2}$, $q\in W^{1,1}_{loc} (\overline{\Omega})$ must be written as~$ v_0=v_{0,\theta} {\bf e}_\theta$, where $v_{\theta,0}=v_{\theta,0}(r)$ obeys from \eqref{eq.v_thetan} the ordinary differential equation
\begin{align}\label{proof.prop.all.alpha.0-1}
-\frac{\dd^2 v_{\theta,0} }{\dd r^2} - \frac1r \frac{\dd v_{\theta,0} }{\dd r}  + \frac{v_{\theta,0} }{r^2}  = f_{\theta,0}\,, \quad r>1\,, \qquad v_{\theta,0}(1)=0\,.
\end{align}
The bounded solution to \eqref{proof.prop.all.alpha.0-1} is written as 
\begin{align}\label{proof.prop.all.alpha.0-2}
v_{\theta,0} (r) = \frac12 \bigg ( - \frac{1}{r} \int_1^\infty f_{\theta,0} \dd s + \frac1r \int_1^r s^2 f_{\theta,0} \dd s + r \int_r^\infty f_{\theta,0} \dd s \bigg )\,.
\end{align}
We note that 
\begin{align*}
\|f_{\theta,0}\|_{L^1 (\Omega)} = 2\pi \int_1^\infty |f_{\theta,0} | \, s \dd s\,.
\end{align*}
Thus we see from \eqref{proof.prop.all.alpha.0-2} that
\begin{align*}
\| v_{\theta,0} \|_{L^\infty(\Omega)} + \| r \frac{\dd v_{\theta,0}}{\dd r} \|_{L^\infty (\Omega)}  \leq C \| f_{\theta,0} \|_{L^1 (\Omega)}\,, 
\end{align*}
which implies \eqref{est.prop.all.alpha.0}. 

\subsubsection{A priori estimate of the $n$ mode with $|n|\geq 1$}\label{subsubsec.middle}

Let $v$ denote a solution to \eqref{S_alpha} satisfying $v\in \mathcal{Q} L^2_\sigma (\Omega)\cap  W^{1,2}_0 (\Omega)^2\cap W^{2,2}_{loc} (\overline{\Omega})^2$ for $f\in \mathcal{Q}_0 L^2 (\Omega)^2$ with some $q\in W^{1,2}_{loc} (\overline{\Omega})$. Then, as we have already seen in the beginning of Section \ref{subsec.general}, for each~$n\ne 0$, the $n$ mode $v_n =\mathcal{P}_n v$ belongs in addition to $W^{2,2}(\Omega)^2$ and we also have that~$\mathcal{P}_n q$ belongs to~$ W^{1,2} (\Omega)$. Hence the energy computation below for $(v_{r,n},v_{\theta,n})$ to the system~\eqref{eq.v_rn}-\eqref{eq.v_thetan} is rigorously verified. With this important remark in mind we multiply both sides of~\eqref{eq.v_rn} by $r\bar{v}_{r,n}$ and of \eqref{eq.v_thetan} by $r\bar{v}_{\theta,n}$ and integrate over $[1,\infty)$, which results in the following identities:
\begin{equation} \label{proof.apriori.middle.1}
\| \nabla v_n \|_{L^2(\Omega)}^2   =  -\alpha \Re \langle U^\bot {\rm rot}\, v_n\,, v_n\rangle_{L^2(\Omega)} +  \Re \langle f_n\,, v_n\rangle _{L^2(\Omega)}\,,\end{equation}
\begin{equation}\label{proof.apriori.middle.2}
-\alpha n  \big ( \| v_{r,n} \|_{L^2 (\Omega)}^2 + \| v_{\theta,n} \|_{L^2(\Omega)}^2 \big )  + \alpha \Im \langle U^\bot {\rm rot}\, v_n\,, v_n \rangle_{L^2(\Omega)}      =  \Im \langle f_n \,, v_n\rangle _{L^2 (\Omega)}\,.
\end{equation}
Note that $U^\bot (x) =\displaystyle -\frac{x}{|x|^2}$
thus we see from ${\rm rot}\, v_n = - \Delta \psi_n$ by   definition of the streamfunction $\psi_n$ for the $n$ mode with $n\ne 0$, 
\begin{align*}
\Re \langle U^\bot {\rm rot}\, v_n\,, v_n\rangle_{L^2(\Omega)}  & =  \Re \langle \frac{1}{|x|} \Delta \psi_n \,, v_{r,n} \rangle_{L^2(\Omega)}\nonumber \\
& = 2\pi n \,  \Im \int_1^\infty (\partial_r^2 \psi_n + \frac1r \partial_r \psi_n  - \frac{n^2}{r^2}  \psi_n ) \overline{\psi_n}  \frac{\dd r}{r}\\
& = 4 \pi n \,  \Im \int_1^\infty \frac1r \partial_r \psi_n  \, \overline{\psi_n}  \frac{\dd r}{r}\,\cdotp
\end{align*}
This gives the bound 
\begin{align}\label{proof.apriori.middle.3}
\Big | \Re \langle U^\bot {\rm rot}\, v_n\,, v_n\rangle_{L^2(\Omega)} \Big | \leq 2 \big \| \frac{v_{r,n}}{r} \big\|_{L^2(\Omega)}\big \| \frac{v_{\theta,n}}{r} \big\|_{L^2 (\Omega)}\,.
\end{align}
Therefore, \eqref{proof.apriori.middle.1} and \eqref{proof.apriori.middle.3} with the lower bound \eqref{polar.grad.n'} imply
\begin{align}\label{proof.apriori.middle.4}
\begin{split}
& \| \partial_r v_{r,n} \|_{L^2(\Omega)}^2 + \| \partial_r v_{\theta,n} \|_{L^2(\Omega)}^2 + \big ( (|n|-1)^2-|\alpha|  \big ) \big ( \| \frac{v_{r,n}}{r} \|_{L^2(\Omega)}^2 +  \| \frac{v_{\theta,n}}{r} \|_{L^2 (\Omega)}^2 \big )\\
& \qquad \qquad  \leq| \Re \langle f_n\,, v_n \rangle_{L^2 (\Omega)}|\,.
\end{split}
\end{align}
Next we study identity \eqref{proof.apriori.middle.2}. We see  that
\begin{align*}
\Im \langle U^\bot {\rm rot}\, v_n\,, v_n \rangle_{L^2(\Omega)}  & =  \Im \langle \frac{1}{|x|} \Delta \psi _n\,, v_{r,n} \rangle_{L^2(\Omega)}\nonumber \\
& = -  2\pi n \,  \Re \int_1^\infty (\partial_r^2 \psi_n + \frac1r \partial_r \psi_n  - \frac{n^2}{r^2}  \psi_n ) \overline{\psi_n}  \frac{\dd r}{r}\,\cdotp
\end{align*}
Integrations by parts yield
\begin{align*} 
& \Re \int_1^\infty (\partial_r^2 \psi_n + \frac1r \partial_r \psi_n  - \frac{n^2}{r^2}  \psi_n ) \overline{\psi_n}  \frac{\dd r}{r} \nonumber \\
& = - \int_1^\infty |\partial_r \psi_n |^2 \frac{\dd r}{r} + 2 \Re \int_1^\infty \frac1r \partial_r \psi_n \overline{\psi_n} \frac{\dd r}{r} - \int_1^\infty \frac{n^2|\psi_n|^2}{r^2} \frac{\dd r}{r} \nonumber \\
& = - \int_1^\infty |\partial_r \psi_n |^2 \frac{\dd r}{r} -  (n^2-2) \int_1^\infty \frac{|\psi_n|^2}{r^2} \frac{\dd r}{r}\,\cdotp
\end{align*}
Hence,
\begin{align}\label{proof.apriori.middle.5}
\Im \langle U^\bot {\rm rot}\, v_n\,, v_n\rangle_{L^2(\Omega)} =  n \big ( \| \frac{v_{\theta,n}}{r}\|_{L^2 (\Omega)}^2 + (1-\frac{2}{n^2}) \| \frac{v_{r,n}}{r}\|_{L^2 (\Omega)}^2 \big )\,.
\end{align}
Hence, \eqref{proof.apriori.middle.2} and \eqref{proof.apriori.middle.5} give
\begin{equation}\label{proof.apriori.middle.6}\begin{aligned}
\alpha n  \bigg ( \| v_{r,n} \|_{L^2(\Omega)}^2 - (1-\frac{2}{n^2}) \| \frac{v_{r,n}}{r}\|_{L^2 (\Omega)}^2 +  \| v_{\theta,n} \|_{L^2 (\Omega)}^2 &- \| \frac{v_{\theta,n}}{r}\|_{L^2 (\Omega)}^2 \bigg ) \\
&= - \Im \langle f_n\,, v_n \rangle_{L^2(\Omega)}\,,
\end{aligned}
\end{equation}
which in particular leads to 
\begin{align}
\big \| \frac{\sqrt{r^2-1}}{r} v_{r,n} \big\|_{L^2 (\Omega)}^2 +\big \| \frac{\sqrt{r^2-1}}{r} v_{\theta,n} \big\|_{L^2(\Omega)}^2 & \leq \frac{1}{|\alpha n|}\big |\Im \langle f_n\,, v_n \rangle_{L^2(\Omega)}\big|\,, \label{proof.apriori.middle.7}\\
\big\|  \frac{v_{r,n}}{r}\big\|_{L^2 (\Omega)}^2& \leq\big |\frac{n}{2\alpha}\big| \, \big |\Im \langle f_n\,, v_n \rangle_{L^2(\Omega)}\big|\,.\label{proof.apriori.middle.7'}
\end{align}
Then, from \eqref{proof.apriori.middle.6} and \eqref{proof.apriori.middle.7'} we gather
\begin{align}\label{proof.apriori.middle.8}
\| v_{r,n} \|_{L^2 (\Omega)}^2 & \leq \| \frac{v_{r,n}}{r} \|_{L^2(\Omega)}^2 + \frac{1}{|\alpha n|}  \big |\Im \langle f_n\,, v_n \rangle_{L^2(\Omega)} \big |\nonumber \\
& \leq \big (\big |\frac{n}{2\alpha}\big| + \frac{1}{|\alpha n|} \big ) \big |\Im \langle f_n\,, v_n \rangle_{L^2(\Omega)} \big | \nonumber \\
& \leq \big|\frac{3n}{2\alpha}\big |\, \big |\Im \langle f_n\,, v_n \rangle_{L^2(\Omega)} \big | \,, \qquad |n|\geq 3\,,
\end{align}
while for $|n|=1,2$, a slightly finer estimate is available from \eqref{proof.apriori.middle.6}:
\begin{align}\label{proof.apriori.middle.8}
\| v_{r,n} \|_{L^2 (\Omega)}^2 & \leq \frac{1}{|\alpha n|} \big |\Im \langle f_n\,, v_n \rangle_{L^2(\Omega)} \big | \,, \qquad |n|=1,2\,.
\end{align}
To obtain the estimate of $\| v_{n} \|_{L^2 (\Omega)}$ we first observe that the following interpolation inequality holds for any scalar function $g\in W^{1,2} ((1,\infty); r\dd r)$:
\begin{align}\label{proof.apriori.middle.9}
\| g \|_{L^2(\Omega)} \leq C \| \partial_r g \|_{L^2 (\Omega)}^\frac13 \| \frac{\sqrt{r^2-1}}{r} g \|_{L^2 (\Omega)}^\frac23 + C \| \frac{\sqrt{r^2-1}}{r} g \|_{L^2 (\Omega)}\,. 
\end{align}
See Appendix \ref{appendix.interpolation} for the proof of \eqref{proof.apriori.middle.9}.
 From \eqref{proof.apriori.middle.9} and  \eqref{proof.apriori.middle.7} we have 
\begin{align}
|\alpha |\big ( \| v_{r,n} \|_{L^2 (\Omega)}^2 + \| v_{\theta,n} \|_{L^2(\Omega)}^2 \big ) 
& \leq C |\alpha|^\frac32 \big ( \big\| \frac{\sqrt{r^2-1}}{r} v_{r,n}\big \|_{L^2 (\Omega)}^2 +\big \| \frac{\sqrt{r^2-1}}{r} v_{\theta,n}\big \|_{L^2(\Omega)}^2 \big )  \nonumber \\
& \quad + \frac14 \big ( \| \partial_r v_{r,n} \|_{L^2(\Omega)}^2 + \| \partial_r v_{\theta,n} \|_{L^2 (\Omega)}^2 \big )\nonumber \\
& \qquad + C|\alpha|  \big ( \big\| \frac{\sqrt{r^2-1}}{r} v_{r,n}\big \|_{L^2 (\Omega)}^2 +\big \| \frac{\sqrt{r^2-1}}{r} v_{\theta,n}\big \|_{L^2(\Omega)}^2 \big ) \nonumber\\ 
& \leq  C(\frac{|\alpha|^\frac12}{|n|} + \frac{1}{|n|} )  \| f_n \|_{L^2(\Omega)} \| v_n \|_{L^2 (\Omega)} \nonumber \\
& \quad  + \frac14 \big ( \| \partial_r v_{r,n} \|_{L^2(\Omega)}^2 + \| \partial_r v_{\theta,n} \|_{L^2 (\Omega)}^2 \big )\,,\nonumber 
\end{align}
and thus,
\begin{align}\label{proof.apriori.middle.11}
\begin{split}
& |\alpha |\big ( \| v_{r,n} \|_{L^2 (\Omega)}^2 + \| v_{\theta,n} \|_{L^2(\Omega)}^2 \big ) \\
& \leq \frac{C}{n^2} ( 1+\frac{1}{|\alpha|} ) \| f_n \|_{L^2 (\Omega)}^2 +  \frac12 \big ( \| \partial_r v_{r,n} \|_{L^2(\Omega)}^2 + \| \partial_r v_{\theta,n} \|_{L^2 (\Omega)}^2 \big )\,.
\end{split}
\end{align}
Hence \eqref{proof.apriori.middle.4} and \eqref{proof.apriori.middle.11} imply 
\begin{align}\label{proof.apriori.middle.12}
\| \partial_r v_{r,n} \|_{L^2 (\Omega)}^2 + \| \partial_r v_{\theta,n} \|_{L^2 (\Omega)}^2  & \leq \frac{C}{n^2} (1+ \frac{1}{|\alpha|} ) \| f_n \|_{L^2 (\Omega)}^2 + 2  |\Re \langle f_n\,, v_n \rangle_{L^2 (\Omega)}|\,.
\end{align} 
Then \eqref{proof.apriori.middle.11} and \eqref{proof.apriori.middle.12} yield
\begin{align*}
\| v_{r,n} \|_{L^2 (\Omega)}^2 + \| v_{\theta,n} \|_{L^2(\Omega)}^2  \leq \frac{C}{|\alpha| n^2} (1+\frac{1}{|\alpha|}) \| f_n \|_{L^2 (\Omega)}^2 + \frac{1}{|\alpha|} \| f_n \|_{L^2(\Omega)}  \| v_n \|_{L^2 (\Omega)}\,,
\end{align*}
that is,
\begin{align}\label{proof.apriori.middle.13}
\| v_{n} \|_{L^2 (\Omega)}^2 \leq C \big ( \frac{1}{|\alpha| n^2} + \frac{1}{|\alpha|^2} \big ) \| f_n \|_{L^2 (\Omega)}^2\,.
\end{align}
This proves \eqref{est.prop.all.alpha.1} and \eqref{est.prop.large.alpha.1}.
Note that the factor $\frac{1}{|\alpha| n^2}$ dominates $\frac{1}{\alpha^2}$ in the regime $|n|\leq O(|\alpha|^\frac12)$ and $|\alpha|\geq 1$. Once we have proved \eqref{proof.apriori.middle.13} the following estimates are immediately obtained from \eqref{proof.apriori.middle.7} and \eqref{proof.apriori.middle.12}:
\begin{align}
\big\| \frac{\sqrt{r^2-1}}{r} v_{r,n} \big\|_{L^2 (\Omega)}^2 +\big \| \frac{\sqrt{r^2-1}}{r} v_{\theta,n} \big\|_{L^2(\Omega)}^2 & \leq \frac{C}{|\alpha n|} \big ( \frac{1}{|\alpha|^\frac12 |n|} + \frac{1}{|\alpha|} \big ) \| f_n \|_{L^2(\Omega)}^2\,, \label{proof.apriori.middle.14} \\
\| \partial_r v_{r,n} \|_{L^2 (\Omega)}^2 + \| \partial_r v_{\theta,n} \|_{L^2 (\Omega)}^2 & \leq C \big ( \frac{1}{n^2} + \frac{1}{|\alpha|} \big ) \| f_n \|_{L^2(\Omega)}^2\,.\label{proof.apriori.middle.15}
\end{align}
The constant $C$ in \eqref{proof.apriori.middle.13}, \eqref{proof.apriori.middle.14}, and \eqref{proof.apriori.middle.15} is independent of $|n|\geq 1$ and~$\alpha$. 
Inequality in \eqref{proof.apriori.middle.14} proves \eqref{est.prop.all.alpha.2}.
 Moreover, \eqref{proof.apriori.middle.13} and \eqref{proof.apriori.middle.4} prove \eqref{est.prop.all.alpha.3} and~\eqref{est.prop.large.alpha.3} since 
 $$
 |\alpha|\| \frac{v_n}{|x|} \|_{L^2 (\Omega)}^2 \leq |\alpha | \| v_n \|_{L^2 (\Omega)}^2\, .
 $$
Finally  \eqref{est.prop.all.alpha.1'} and \eqref{est.prop.large.alpha.2}  are obtained 
from the interpolation inequality for functions in the space~$W_0^{1,2}((1,\infty); \dd r)$.

\subsubsection{A priori estimate of the $n$ mode with $|n|\gg O (\sqrt{1+|\alpha|})$}\label{subsubsec.high.n}
If $|n|\geq 1 +\sqrt{2|\alpha|}$ then \eqref{proof.apriori.middle.4} yields 
\begin{align}\label{proof.apriori.high.1}
\begin{split}
\| \partial_r v_{r,n} \|_{L^2(\Omega)}^2 + \| \partial_r v_{\theta,n} \|_{L^2(\Omega)}^2 + \frac{n^2}{8} \big ( \| \frac{v_{r,n}}{r} \|_{L^2(\Omega)}^2 +  \| \frac{v_{\theta,n}}{r} \|_{L^2 (\Omega)}^2 \big )\\
\leq \big| \Re \langle f_n\,, v_n \rangle_{L^2 (\Omega)} \big|\,.
\end{split}
\end{align}
Then \eqref{proof.apriori.middle.6}  and \eqref{proof.apriori.high.1} give
\begin{align*}
\| v_n \|_{L^2 (\Omega)}^2 & \leq \|\frac{v_n}{|x|} \|_{L^2 (\Omega)}^2 + \frac{1}{|\alpha n|} \| f_n \|_{L^2 (\Omega)} \| v_n \|_{L^2 (\Omega)} \\
& \leq \big ( \frac{8}{n^2}  + \frac{1}{|\alpha n|} \big ) \| f_n \|_{L^2 (\Omega)} \| v_n \|_{L^2 (\Omega)} \,,
\end{align*}
which shows 
\begin{align}\label{proof.apriori.high.2}
\| v_n \|_{L^2 (\Omega)}^2 \leq \big ( \frac{8}{n^2}  + \frac{1}{|\alpha n|} \big ) ^2 \| f_n \|_{L^2 (\Omega)}^2\,.
\end{align} 
Note that \eqref{proof.apriori.high.2} is better than \eqref{proof.apriori.middle.13} in the regime $|n|\gg 1+\sqrt{2 |\alpha|}$, while both are of the same order in the regime $|n|=O(\sqrt{1+|\alpha|})$.
The estimates \eqref{proof.apriori.high.1} and \eqref{proof.apriori.high.2} yield
\begin{align}\label{proof.apriori.high.3}
\| \nabla v_n \|_{L^2 (\Omega)}^2 \leq  C \big ( \frac{1}{n^2} + \frac{1}{|\alpha n|} \big ) \| f_n \|_{L^2(\Omega)}^2 \,,
\end{align}
while \eqref{proof.apriori.middle.7} and \eqref{proof.apriori.high.2} lead to 
\begin{align}\label{proof.apriori.high.4}
\| \frac{\sqrt{r^2-1}}{r} v_{r,n} \|_{L^2 (\Omega)}^2 + \| \frac{\sqrt{r^2-1}}{r} v_{\theta,n} \|_{L^2(\Omega)}^2 \leq \frac{1}{|\alpha n|}   \big ( \frac{8}{n^2} + \frac{1}{|\alpha n|} \big ) \| f_n \|_{L^2(\Omega)}^2\,.
\end{align}
Again, \eqref{proof.apriori.high.3} and \eqref{proof.apriori.high.4} are better than \eqref{proof.apriori.middle.15} and \eqref{proof.apriori.middle.14} in the regime  $|n|\gg 1 + \sqrt{2|\alpha|}$.
The estimates \eqref{proof.apriori.high.2}, \eqref{proof.apriori.high.3}, and \eqref{proof.apriori.high.4} show \eqref{est.prop.all.alpha.4}, \eqref{est.prop.all.alpha.6}, and \eqref{est.prop.all.alpha.5}. Then \eqref{est.prop.all.alpha.4'} follows by interpolation using \eqref{est.prop.all.alpha.4} and \eqref{est.prop.all.alpha.6}

\subsubsection{Proof of Proposition \ref{prop.all.alpha}}\label{subsubsec.all.alpha}

The statement (i) of Proposition \ref{prop.all.alpha} is proved in Subsection \ref{subsubsec.0}.
In particular, we have  that~$v_0=v_{\theta,0}{\bf e}_\theta$ and $v_{\theta,0}$ is given by the formula \eqref{proof.prop.all.alpha.0-2}.
It remains to prove (ii) of Proposition \ref{prop.all.alpha}.

\medskip

\noindent (Estimates and Uniqueness) We have already proved the a priori estimates of solutions in Subsections \ref{subsubsec.middle} and \ref{subsubsec.high.n}, which give \eqref{est.prop.all.alpha.1}-\eqref{est.prop.all.alpha.6}. The uniqueness of solutions directly follows from these a priori estimates.

\medskip
\noindent (Existence) By considering the Helmholtz-Leray projection which is bounded in $L^2$ and invariant under the action of $\mathcal{P}_n$, we may assume that $f_n$ belongs to $\mathcal{Q}_0 L^2_\sigma (\Omega)$, rather than~$\mathcal{Q}_0 L^2(\Omega)^2$.
To show the existence of solutions we consider the operator 
\begin{align*}
\mathbb{A}_\alpha := \mathbb{P} \Delta - \alpha \mathbb{P}U^\bot {\rm rot}
\end{align*}
in the $L^2$ framework, where $\mathbb{P}: L^2(\Omega)^2 \rightarrow L^2_\sigma (\Omega)$ is the Helmholtz-Leray projection and $\Delta$ is the Dirichlet Laplacian in $L^2(\Omega)$.  
The operator $\mathbb{P}\Delta$ is thus the standard Stokes operator in~$L^2_\sigma (\Omega)$.
Let us consider the operator $\mathbb{A}_\alpha$ in the invariant space $\mathcal{P}_n L^2_\sigma (\Omega)$, $n\ne 0$. Note that the spectrum of the Stokes operator $\mathbb{P}\Delta$ in $\mathcal{P}_n L^2_\sigma (\Omega)$ is included in the half real line $\overline{\R_-}=\{\lambda\leq 0\}$, while the operator $\mathbb{P}U^\bot {\rm rot}$ is relatively compact with respect to $\mathbb{P}\Delta$ in $\mathcal{P}_n L^2_\sigma (\Omega)$, for $U^\bot$ is smooth and decays at infinity. Thus, the difference between the spectrum of $\mathbb{A}_\alpha$ and the one of $\mathbb{P}\Delta$ consists only of discrete eigenvalues with finite multiplicities.
Then the scalar~$\lambda:=-i \alpha n$ must belong to the resolvent set of $\mathbb{A}_\alpha$ in $\mathcal{P}_n L^2_\sigma (\Omega)$, otherwise $-i\alpha n$ is an eigenvalue of $\mathbb{A}_\alpha$ in $\mathcal{P}_n L^2_\sigma (\Omega)$ but this cannot be true due to the a priori estimates on solutions of \eqref{eq.reduced.stokes} with $\lambda=-i\alpha n$ which we have   shown above.
Hence, when $\lambda=-i\alpha n$,  for any $f_n = \mathcal{P}_n f\in \mathcal{P}_n L^2_\sigma (\Omega)$ there exists a unique solution $v_n$  to \eqref{eq.reduced.stokes} belonging to the space~$\mathcal{P}_n L^2_\sigma (\Omega)\cap W^{1,2}_0 (\Omega)^2 \cap W^{2,2}(\Omega)^2$ with a suitable pressure field belonging to~$W^{1,2}_{loc} (\overline{\Omega})$. 
Then $v=\sum_{n\ne 0} v_n$ belongs to $\mathcal{Q}_0 L^2_\sigma (\Omega)\cap W^{1,2}_0 (\Omega)^2 \cap W^{2,2}(\Omega)^2$ and  solves \eqref{S_alpha} by   construction, for \eqref{eq.reduced.stokes} with $\lambda=-i\alpha n$ is equivalent to \eqref{eq.v_rn} and \eqref{eq.v_thetan} for each $n\ne 0$.
The proof of (ii) of Proposition~\ref{prop.all.alpha} is complete. \BOX

\subsection{Analysis in the fast rotation case $|\alpha| \gg 1$}\label{subsec.fast}

In this subsection we focus on the  behavior of solutions to \eqref{S_alpha} in the case~$|\alpha|\gg 1$.
Let us define the parameter
\begin{align}\label{def.beta}
\beta & := (-2i\alpha n)^\frac13  =
\begin{cases}
&  (2|\alpha n|)^{\frac13} c_- \quad \mbox{if} \quad \alpha n >0\,,\\
&  (2|\alpha n|)^{\frac13} c_+ \quad \mbox{if} \quad \alpha n <0 \,,
\end{cases}
\end{align}
where 
\begin{equation}\label{defc} 
c_\pm := \frac{\sqrt{3}\pm i}{2} \,\cdotp
\end{equation}

\noindent Our goal is to prove the following structure result on solutions to \eqref{S_alpha}. 
\begin{proposition}\label{decompositionvelocity}
There is a constant~$\kappa\in (0,1)$ such that  as long as~$1 \leq |n| \leq \kappa |\alpha|^\frac12$ and~$|\beta|  $ large enough (implying~$|\alpha|$ large enough), the~$n$ mode of the  velocity field~$v_n$ solving~\eqref{S_alpha} may be decomposed into four parts
$$
v_n =  v_n^{\rm slip} + v_n^{\rm slow} + v_{n,{ \rm BL}} + \widetilde v_n \, , 
$$
where
\begin{itemize}
\item the term $v_n^{\rm slip}$ satisfies the system~$({\rm mS}_{\alpha})$  and the estimates~{\rm(\ref{est.prop.mS_alpha.6})-(\ref{est.prop.mS_alpha.8})} of Proposition~{\rm\ref{prop.mS_alpha}}.

\item the  stream function of~$v_n^{\rm slow}$ is given by
$$
\psi_n^{\rm slow} (r) = a_n r^{-|n|}\,,
$$
where
$$
|a_n| \leq C |\alpha n|^{-\frac56} \|f_n\|_{L^2} \, .
$$ 
\item the boundary layer term~$v_{n,{ \rm BL}}$ is given by  
$$
v_{n,{ \rm BL},r} (r)= \frac{inb_n}{r}    G_{n,\alpha}  \big( |\beta| (r -1) \big) \, , \quad v_{n,{ \rm BL},\theta} =-|\beta| b_n   G_{n,\alpha}'  \big( |\beta| (r -1) \big)  
$$
with~$  G_{n,\alpha} $ a smooth function,   decaying exponentially at infinity, uniformly in~$n$ and~$\alpha$, and  where
$$
|b_n| \leq C |\alpha n|^{-\frac56} \|f_n\|_{L^2} \, .
$$ 

\item the  term  $\widetilde v_n$ is a remainder term in the sense that
\begin{equation}\label{est.decompositionvelocitytilde}
\begin{aligned}
\|\widetilde v_n \|_{L^2(\Omega)}  \le C  |\alpha n |^{-1} \|f_n\|_{L^2} \, , \quad  \|\widetilde v_n \|_{L^\infty (\Omega)}   \le C |\alpha n|^{-\frac56} \|f_n\|_{L^2}\,,  \\\mbox{and}\quad
\| \nabla\widetilde v_n \|_{L^2(\Omega) } \le C  |\alpha n|^{-\frac23}\|f_n\|_{L^2} \,.
\end{aligned}
\end{equation}
\end{itemize}
In particular, the following estimates hold for $v_n$:
\begin{equation}\label{est.decompositionvelocity}
\begin{aligned}
\| v_n \|_{L^2(\Omega)}   \le C |\alpha n|^{-\frac23}\|f_n\|_{L^2}  \,,\quad
\| v_n \|_{L^\infty (\Omega)}   \le C |\alpha n|^{-\frac12} \|f_n\|_{L^2}\,,  \\
\mbox{and}\quad \| \nabla v_n \|_{L^2(\Omega)} \le C  |\alpha n|^{-\frac13}\|f_n\|_{L^2} \,.
\end{aligned}
\end{equation}
The constant $C$ is independent of $n$, $\alpha$, and $f_n$. 
\end{proposition}

\begin{remark}\label{rem.decompositionvelocity}{\rm (1) The velocity fields~$v_n^{\rm slip}, v_n^{\rm slow}$ and~$v_{n,{\rm BL}}$ have the same decay order in the spaces~$L^2$,~$L^\infty$ and $H^1$, although their structure  is different. In particular contrary to the other velocity fields, the term~$v_{n, {\rm BL}}$ is negligible away from the boundary due to its rapid decay.
The key point is that the boundary layer analysis enables us to improve the decay order in the low and middle frequencies~$|n|\ll O(|\alpha|^\frac12)$, compared with the results of Proposition \ref{prop.all.alpha} which are based only on   energy computations. Indeed, it is not clear whether \eqref{est.decompositionvelocity} can be shown only from   energy computations without using the boundary layer analysis.

\noindent  (2) The vorticity of $v_n^{{\rm slow}}$ vanishes in $\Omega$, that is, $v_n^{{\rm slow}}$ is irrotational in $\Omega$. The function~$G_{n,\alpha}$ and its derivatives of finite order are uniformly bounded in $n$ and $\alpha$ which satisfy $|\alpha|\geq 1$ and~$1\leq |n|\leq |\alpha|^\frac12$.  The uniform decay estimate of $G_{n,\alpha} (\rho)$ for $\rho\gg 1$ is stated in~\eqref{nondegenerate.6} below.
The smallness of $\kappa$ in Proposition \ref{decompositionvelocity} is required only in obtaining some lower bound for the quantity $\Big | |\beta| G_{n,\alpha}' (0) + |n| G_{n,\alpha} (0)| \Big |$ which is essential to construct the solution satisfying the noslip boundary condition and possessing the structure stated in Proposition~\ref{decompositionvelocity}. 
}
\end{remark}
\begin{proofprop}{decompositionvelocity}
We first consider a linearized problem similar to \eqref{S_alpha} but with a different boundary condition,
such that the vorticity vanishes on the boundary:
\begin{equation}\tag{${\rm mS}_{\alpha}$}\label{mS_alpha}
  \left\{
\begin{aligned}
  -\Delta v - \alpha ( x^\bot \cdot \nabla v - v^\bot ) + \nabla q + \alpha U^\bot {\rm rot}\, v & = 
  f \,,  \qquad \qquad x\in \Omega \,,\\
{\rm div}\, v &  = 0\,,   \qquad \qquad  x \in \Omega\,, \\
  v_r  = {\rm rot}\, v & =  0 \,, \qquad \qquad  x \in \partial \Omega\,. \\
\end{aligned}\right.
\end{equation}
Note that the vorticity field $\omega$ then satisfies the following equations 
\begin{equation}\label{eqw}
\left\{
\begin{aligned}
& -\Delta \omega - \alpha x^\bot \cdot \nabla \omega + \alpha U \cdot \nabla \omega =  {\rm rot}\, f\,, \quad x\in \Omega \,,\\
& \omega|_{\partial\Omega} =0\,.
\end{aligned}\right.
\end{equation}
Let us prove the following proposition.
\begin{proposition}\label{prop.mS_alpha} For any $\alpha\in \R$ with $|\alpha|\geq 1$ and external force $f\in \mathcal{Q}_0 L^2 (\Omega)^2$ the system \eqref{mS_alpha} admits a unique solution $(v,\nabla q)$ with $v\in \mathcal{Q}_0 L^2_\sigma (\Omega) \cap  W^{1,2}(\Omega)^2 \cap  W^{2,2}_{loc} (\overline{\Omega})$ and $q\in W^{1,2}_{loc} (\overline{\Omega})$. Moreover, $\omega={\rm rot}\, v\in \mathcal{Q}_0 W^{1,2}_0 (\Omega)$ satisfies for $n\ne 0$,
\begin{align}
\| \mathcal{P}_n \omega \|_{L^2(\Omega)} & \leq \frac{C}{|\alpha n|^\frac13} \| f_n \|_{L^2(\Omega)}\,, \label{est.prop.mS_alpha.1} \\
\| \frac{\sqrt{|x|^2-1}}{|x|} \mathcal{P}_n \omega \|_{L^2 (\Omega)} & \leq \frac{1}{|\alpha n|^\frac12} \| f_n \|_{L^2 (\Omega)}\,,\label{est.prop.mS_alpha.2}\\
\| \nabla \mathcal{P}_n \omega \|_{L^2(\Omega)} & \leq \| f_n \|_{L^2 (\Omega)} \,,\label{est.prop.mS_alpha.3}
\end{align}
while  $v_n = \mathcal{P}_n v$ satisfies   for $1\leq |n|\leq 2|\alpha|^\frac12$,
\begin{align}
\| v_{n} \|_{L^2(\Omega)} & \leq \frac{C}{|\alpha n|^\frac23} \| f_n \|_{L^2 (\Omega)}\,,\label{est.prop.mS_alpha.6}\\
\| v_{n} \|_{L^\infty(\Omega)} & \leq \frac{C}{|\alpha n|^\frac12} \| f_n \|_{L^2 (\Omega)}\,,\label{est.prop.mS_alpha.7}\\
\| \frac{\sqrt{|x|^2-1}}{|x|} v_n \|_{L^2 (\Omega)} & \leq \frac{C}{|\alpha n|^\frac56} \| f_n \|_{L^2 (\Omega)}\,.\label{est.prop.mS_alpha.8}
\end{align}
Here $f_n:=\mathcal{P}_n f$ and $C$ is independent of $n$ and $\alpha$. 
\end{proposition}

\begin{remark}\label{rem.prop.mS_alpha}{\rm Note that 
\begin{align}
\| \nabla v_n \|_{L^2(\Omega)} \leq C \big ( \| \mathcal{P}_n \omega \|_{L^2 (\Omega)} + \| \frac{v_n}{|x|} \|_{L^2 (\Omega)} \big ) \leq \frac{C}{|\alpha n|^\frac13} \| f_n \|_{L^2(\Omega)}\,.\label{est.prop.mS_alpha.1'}
\end{align}
The decay order in each estimate of \eqref{est.prop.mS_alpha.1'} and \eqref{est.prop.mS_alpha.6}-\eqref{est.prop.mS_alpha.8} is better than the order in~\eqref{est.prop.all.alpha.1}-\eqref{est.prop.all.alpha.3} for the solution subject to the noslip boundary condition. This faster decay is shown below by an energy estimate thanks to  the slip boundary condition in \eqref{mS_alpha}, which reduces the magnitude of the boundary layer arising from the fast rotation in low and middle frequencies. 
}
\end{remark}

\begin{proofprop}{prop.mS_alpha} For simplicity we set  $\omega_n = \omega_n (r) := (\mathcal{P}_n \omega) e^{-i n \theta}$.

\noindent (A priori estimates) We first show the a priori estimates stated in \eqref{est.prop.mS_alpha.1}-\eqref{est.prop.mS_alpha.8}. As in the proof of Proposition \ref{prop.all.alpha}, for each $n$ mode the energy computation for $v_n$ based on the integration by parts is verified for any solution $(v, \nabla q)$ to \eqref{mS_alpha} such that $v\in \mathcal{Q}_0 L^2_\sigma (\Omega)\cap  W^{1,2}(\Omega)^2 \cap W^{2,2}_{loc} (\overline{\Omega})$ and $q\in W^{1,2}_{loc} (\overline{\Omega})$ (see the argument at the beginning of Subsection~\ref{subsec.general}). We also note that the $n$ mode of the vorticity $\omega_n\in W^{1,2}_0 ((1,\infty); r \dd r)$ satisfies the following ordinary differential equations on $(1,\infty)$ in the weak sense:
\begin{align}\label{proof.prop.mS_alpha.1}
- \bigg(\frac{\dd^2}{\dd r^2} + \frac1r \frac{\dd}{\dd r}  -\frac{n^2}{r^2}+  i \alpha n \big (  1-\frac1{r^2}\big ) \bigg ) \omega_n = \frac1r \frac{\dd}{\dd r} (r f_{\theta,n}) - \frac{i n}{r} f_{r,n}\,, 
\end{align}
with the Dirichlet boundary condition $\omega_n (1) =0$. 
We can multiply both sides by $r\bar{\omega}_n$ and integrate over $[1,\infty)$, which gives the following identities:
\begin{align}
\int_1^\infty \big ( |\partial_r \omega_n|^2 + n^2 \frac{|\omega_n|^2}{r^2} \big ) r \dd r & = - \Re \bigg (  \int_1^\infty f_{\theta,n} \frac{\dd \bar\omega_n}{\dd r} r \dd r  + in \int_1^\infty f_{r,n} \bar\omega_n \dd r\bigg )\,,\label{proof.prop.mS_alpha.2}\\
\alpha n \int_1^\infty (1-\frac{1}{r^2}) |\omega_n|^2 r \dd r & =   \Im \bigg (  \int_1^\infty f_{\theta,n} \frac{\dd \bar\omega_n}{\dd r} r \dd r  + in \int_1^\infty f_{r,n} \bar\omega_n \dd r\bigg )\,.\label{proof.prop.mS_alpha.3}
\end{align}
The first identity \eqref{proof.prop.mS_alpha.2} gives the bound
\begin{align*}
\| \nabla \mathcal{P}_n \omega\|_{L^2 (\Omega)} \leq \| f_n \|_{L^2 (\Omega)}\,,
\end{align*}
and then the second identity \eqref{proof.prop.mS_alpha.3} yields
\begin{align*}
\big\| \frac{\sqrt{|x|^2-1}}{|x|} \mathcal{P}_n \omega \big\|_{L^2(\Omega)}^2 \leq \frac{1}{|\alpha n|} \| f_n \|_{L^2 (\Omega)}\| \nabla \mathcal{P}_n \omega\|_{L^2 (\Omega)} \leq \frac{1}{|\alpha n|} \| f_n \|_{L^2 (\Omega)}^2\,.
\end{align*}
Thus \eqref{est.prop.mS_alpha.2} and \eqref{est.prop.mS_alpha.3} hold. The estimate \eqref{est.prop.mS_alpha.1} follows from \eqref{est.prop.mS_alpha.2} and \eqref{est.prop.mS_alpha.3} by the interpolation inequality \eqref{proof.apriori.middle.9} and $|\alpha|\geq 1$. 
The estimates on the velocity $v_n$ are obtained from the estimates for the streamfunction $\psi_n$ which is the solution to the Poisson equation~\eqref{eq.stream}. Let us first estimate $\displaystyle \frac{\dd ^2\psi_n}{\dd r^2}$ and $\displaystyle\frac{\psi_n}{r^2}\cdotp$
We write $\psi_n'=\displaystyle\frac{\dd \psi_n}{\dd r}$ and $\psi_n''=\displaystyle\frac{\dd ^2\psi_n}{\dd r^2}$ for simplicity.
We observe that 
\begin{align*}
\int_1^\infty \frac1r \psi_n' \bar{\psi}_n'' r \dd r & = \int_1^\infty \Big(\frac{\psi_n}{r}\Big)' \bar{\psi}_n'' r \dd r + \int_1^\infty \frac{\psi_n}{r^2} \bar{\psi}_n'' r \dd r \,.
\end{align*}
A similar computation yields
\begin{align*}
- \Re \int_1^\infty \frac{\psi_n}{r^2} \bar{\psi}_n'' r\dd r= \Re \int_1^\infty \Big(\frac{\psi_n}{r}\Big)' \bar{\psi}_n' \dd r = \int_1^\infty \Big|\Big (\frac{\psi_n}{r}\Big)'\Big|^2 r \dd r
\end{align*}
and 
\begin{align*}
\Re \int_1^\infty \frac1r \psi_n' \frac{\bar{\psi}_n}{r^2} r \dd r   =  \int_1^\infty\Big |\frac{\psi_n}{r^2}\Big|^2 r \dd r\,.
\end{align*}
Hence we have from \eqref{eq.stream}, by multiplying both sides by $r\bar{\psi}_n''$ and by integrating over $[1,\infty)$,
\begin{align*}
\int_1^\infty  |\psi_n''|^2 r \dd r + \Re \int_1^\infty\Big (\frac{\psi_n}{r}\Big)' \bar{\psi}_n'' r \dd r + (n^2-1) \int_1^\infty \Big| \Big(\frac{\psi_n}{r}\Big)'\Big|^2 r \dd r   = -\Re \int_1^\infty \omega_n \bar{\psi}_n'' r \dd r\,,
\end{align*}
and similarly, by multiplying by $\frac{\bar{\psi}_n}{r}$ in \eqref{eq.stream} and integrating over $[1,\infty)$,
\begin{align*}
\int_1^\infty\Big |\Big (\frac{\psi_n}{r}\Big)'\Big|^2 r \dd r + (n^2 - 1 )  \int_1^\infty\Big |\frac{\psi_n}{r^2}\Big|^2 r \dd r = \Re \int_1^\infty \omega_n \frac{\bar{\psi}_n}{r^2} r \dd r\,.
\end{align*}
Combining these two identities gives the following bounds: 
when $|n|\geq 2$,
\begin{equation}\label{proof.prop.mS_alpha.-1}\begin{aligned}
\| \psi_n''\|_{L^2 (\Omega)}  + \big\| (\frac{ n\psi_n}{r})'\big\|_{L^2(\Omega)} +  n^2 \big \|\frac{\psi_n}{r^2} \big\|_{L^2 (\Omega)} \leq C \| \omega_n \|_{L^2(\Omega)}\,,
\end{aligned}
\end{equation}
while when $|n|=1$,
\begin{align}\label{proof.prop.mS_alpha.-2}
\begin{split}
\| (\frac{\psi_n}{r})'\|_{L^2(\Omega)}^2  & \leq \| \omega_n \|_{L^2(\Omega)} \big \| \frac{\psi_n}{r^2}\big \|_{L^2 (\Omega)}\,,\\
\| \psi_n''\|_{L^2 (\Omega)}^2  & \leq C \| \omega_n \|_{L^2(\Omega)}\big ( \| \omega_n \|_{L^2(\Omega)} + \big\| \frac{\psi_n}{r^2}\big \|_{L^2 (\Omega)}\big )\,,
\end{split}
\end{align}
where $C$ is independent of $n$. 
Next we observe that the identity \eqref{proof.apriori.middle.6} holds even under the slip boundary condition,
and thus, \eqref{proof.apriori.middle.7} and \eqref{proof.apriori.middle.7'} are valid also for solutions to \eqref{mS_alpha}.
Hence, by recalling the relation $v_{r,n}=\frac{in}{r}\psi_n$ and using the interpolation inequality \eqref{proof.apriori.middle.9}, we obtain for $|n|\geq 2$, thanks to~\eqref{proof.prop.mS_alpha.-1} and \eqref{proof.apriori.middle.7},
\begin{align*}
\| v_{r,n} \|_{L^2 (\Omega)}  & \leq C \| \partial_r v_{r,n}\|_{L^2 (\Omega)}^\frac13 \big\| \frac{\sqrt{r^2-1}}{r} v_{r,n} \big\|_{L^2 (\Omega)}^\frac23  + C\big\| \frac{\sqrt{r^2-1}}{r} v_{r,n}\big \|_{L^2 (\Omega)} \nonumber \\
& \leq C \| \omega_n \|_{L^2 (\Omega)}^\frac13 \big ( \frac{1}{|\alpha n|} \| f_n \|_{L^2(\Omega)} \| v_n \|_{L^2(\Omega)} \big )^\frac13  + C \big ( \frac{1}{|\alpha n|} \| f_n \|_{L^2(\Omega)} \| v_n \|_{L^2(\Omega)} \big )^\frac12  \,,
\end{align*}
while for $|n|=1$ we use \eqref{proof.prop.mS_alpha.-2} and also \eqref{proof.apriori.middle.7'}, which give
\begin{align*}
\| v_{r,n} \|_{L^2 (\Omega)} & \leq C \| \partial_r v_{r,n}\|_{L^2 (\Omega)}^\frac13\big \| \frac{\sqrt{r^2-1}}{r} v_{r,n}\big \|_{L^2 (\Omega)}^\frac23  + C\big\| \frac{\sqrt{r^2-1}}{r} v_{r,n} \big\|_{L^2 (\Omega)} \nonumber \\ 
& \leq C \| \omega_n \|_{L^2 (\Omega)}^\frac16 \big ( \frac{1}{|\alpha|} \| f_n \|_{L^2(\Omega)} \| v_n \|_{L^2(\Omega)} \big )^\frac{5}{12} + C \big ( \frac{1}{|\alpha |} \| f_n \|_{L^2(\Omega)} \| v_n \|_{L^2(\Omega)} \big )^\frac12 \,.
\end{align*}
Similarly, we have from $v_{\theta,n}=-\psi_n'$, \eqref{proof.apriori.middle.9},  \eqref{proof.prop.mS_alpha.-1}, and \eqref{proof.apriori.middle.7}, for $|n|\geq 2$,
\begin{align*}
\| v_{\theta,n} \|_{L^2 (\Omega)}  
& \leq  C \| \omega_n \|_{L^2 (\Omega)}^\frac13 \big ( \frac{1}{|\alpha n|} \| f_n \|_{L^2(\Omega)} \| v_n \|_{L^2(\Omega)} \big )^\frac13  + C \big ( \frac{1}{|\alpha n|} \| f_n \|_{L^2(\Omega)} \| v_n \|_{L^2(\Omega)} \big )^\frac12  \,,
\end{align*}
on the other hand, for $|n|=1$, by applying \eqref{proof.apriori.middle.9}, \eqref{proof.prop.mS_alpha.-2}, \eqref{proof.apriori.middle.7} and \eqref{proof.apriori.middle.7'}, 
\begin{align*}
\| v_{\theta,n} \|_{L^2 (\Omega)}  
& \leq  C \| \omega_n \|_{L^2 (\Omega)}^\frac13 \big ( \frac{1}{|\alpha|} \| f_n \|_{L^2(\Omega)} \| v_n \|_{L^2(\Omega)} \big )^\frac13  \\
& \quad + C \| \omega_n \|_{L^2 (\Omega)}^\frac16 \big ( \frac{1}{|\alpha|} \| f_n \|_{L^2(\Omega)} \| v_n \|_{L^2(\Omega)} \big )^\frac{5}{12} + C \big ( \frac{1}{|\alpha |} \| f_n \|_{L^2(\Omega)} \| v_n \|_{L^2(\Omega)} \big )^\frac12  \,.
\end{align*}
Hence we have from \eqref{est.prop.mS_alpha.1} for the estimate of $\|\omega_n \|_{L^2 (\Omega)}$ and $|\alpha|\geq 1$,
\begin{align}\label{proof.prop.mS_alpha.4}
\| v_n \|_{L^2 (\Omega)} & \leq \frac{C}{|\alpha n|^\frac23} \| f_n \|_{L^2(\Omega)}\,, \qquad |n|\geq 1\,.
\end{align}  
The estimate \eqref{est.prop.mS_alpha.8} follows from \eqref{proof.apriori.middle.7} and \eqref{proof.prop.mS_alpha.4}.
As for the $L^\infty$ estimate in the case~$|n|\geq 2$, we have from the interpolation inequality and \eqref{proof.prop.mS_alpha.-1},
\begin{align*}
\| v_{r,n} \|_{L^\infty(\Omega)}^2 \leq C \| \partial_r v_{r,n} \|_{L^2(\Omega)} \| v_{r,n} \|_{L^2(\Omega)} & \leq C \| \omega_n \|_{L^2 (\Omega)} \| v_{r,n} \|_{L^2(\Omega)} \\
& \leq \frac{C}{|\alpha n|} \| f_n \|_{L^2 (\Omega)}^2\,,
\end{align*}
and similarly, 
$$
\| v_{\theta,n} \|_{L^\infty(\Omega)}^2\leq \frac{C}{|\alpha n|} \| f_n \|_{L^2(\Omega)}^2\cdotp 
$$
The case $|n|=1$ is handled in the same manner, and in this case we use \eqref{proof.prop.mS_alpha.-2} instead of \eqref{proof.prop.mS_alpha.-1} to estimate $\| \partial_r v_{r,n} \|_{L^2 (\Omega)}$ and $\| \partial_r v_{\theta,n}\|_{L^2 (\Omega)}$. This modification does not produce any change in the final estimate 
$$
\| v_n \|_{L^\infty (\Omega)}\leq \frac C{|\alpha|^{\frac12}} \| f_n \|_{L^2 (\Omega)}\,.
$$ The details are omitted here. 
Thus \eqref{est.prop.mS_alpha.1}-\eqref{est.prop.mS_alpha.8} hold.

\smallskip

\noindent (Existence and uniqueness) The uniqueness follows 
from the a priori estimates. The existence is shown by the same argument as in the proof of Proposition \ref{prop.all.alpha} in Subsection \ref{subsubsec.all.alpha} (the proof for the statement (ii)), so we omit the details here. The proof of Proposition \ref{prop.mS_alpha} is complete.
\end{proofprop}

Let us return to the proof of Proposition \ref{decompositionvelocity}. 
We denote by $v^{\rm slip}$ the solution to \eqref{mS_alpha}
with the external force $f$ in Proposition \ref{prop.mS_alpha}.
Note that $v^{\rm slip}$ satisfies the desired estimate~\eqref{est.decompositionvelocity}, 
while $v^{\rm slip}$ is not necessarily subject to the noslip boundary condition. 
Hence, starting from the perfect-slip solution $v^{\rm slip}$, 
our next task is to recover the noslip boundary condition  by the boundary layer analysis. 
To this end we consider the following equations related to the stream function:
\begin{equation}\label{eq.fast.slow}
\left\{
\begin{aligned}
& \bigg(\frac{\dd^2}{\dd r^2} + \frac1r \frac{\dd}{\dd r}  -\frac{n^2}{r^2}+  i \alpha n \big (  1-\frac1{r^2}\big ) \bigg ) \big ( \frac{\dd^2 }{\dd r^2} + \frac1r\frac{\dd}{\dd r} - \frac{n^2}{r^2} \big ) \varphi_n  = 0\,, \quad r>1\,,\\
&| \varphi_n (1) | \leq 1\,, \qquad  \lim_{r\rightarrow \infty} \varphi_n =0\,.
\end{aligned}\right.
\end{equation}
Indeed, the first equation in \eqref{eq.fast.slow} is nothing but the equation for the $n$ mode of the stream function in our problem (in polar coordinates), while here we do not need to impose the exact boundary value on the boundary $r=1$, and the key point is that there exists a solution to \eqref{eq.fast.slow} which possess a boundary layer structure due to~$|\alpha n|\gg 1$. We shall decompose it into a boundary layer part and a remainder:
 \begin{equation}\label{defphitilden}
\varphi_{n} = \varphi _{n,{ \rm BL}} + \widetilde \varphi_n  \,,
 \end{equation}
where~$\varphi _{n,{ \rm BL}} $ is a function on the boundary layer variable $|\beta| (r-1)$ and~$ \widetilde \varphi_n$ vanishes at~$r=1$ and is smaller than $\varphi_{n,\rm BL}$ up to its derivatives.  
The precise construction of $\varphi _{n,{ \rm BL}}$ and $\widetilde \varphi_n$ will be stated later.

\noindent Once $\varphi_n$ is constructed as in \eqref{defphitilden}, the $n$ mode of the noslip solution $v_n$ to \eqref{eq.fast.slow}, or equivalently its stream function $\psi_n$, is obtained by the following argument.
Noticing that \eqref{eq.fast.slow} has an explicit solution given by~$r\mapsto  r^{-|n|} $ (corresponding to a vorticity-free solution), we look for the stream function $\psi_n$ under the following form:
\begin{align}\label{def.psi^noslip}
\psi_n  (r)= a_n r^{-|n|} + b_n\varphi_{n}  (r) + \psi_n^{\rm slip} (r)
\end{align}
where   $\psi_n^{\rm slip}$ is the $n$ mode of the streamfunction for $v_n^{\rm slip}$, namely,
\begin{align}\label{def.psi^slip}
\psi^{\rm slip}_n = \frac{r}{in} v_{r,n}^{\rm slip}\,, \qquad n\ne 0\,,
\end{align}
and the coefficients~$a_n$ and~$b_n$ are determined from the  prescription that~$\psi_n$ and~$\psi_n'$ vanish at the boundary~$r=1$:
\begin{equation*}
\begin{aligned}
a_n + b_n\varphi_n(1) &= 0 \\
-|n|a_n +  b_n   \frac{\dd \varphi_n }{\dd r} (1)  &= - \frac{\dd  \psi_n^{\rm slip}}{\dd r}  (1)  \, \big (= v_{\theta,n}^{\rm slip} (1) \big )\,.
\end{aligned}
\end{equation*}
In other words there holds
 \begin{equation}\label{defbn}
b_n \Big(  \frac{\dd \varphi_n }{\dd r} (1)  + |n|\varphi_n(1) \Big) = v_{\theta,n}^{\rm slip} (1)\,, \qquad a_n = - b_n \varphi_n (1) \,.
\end{equation}
The key point is that  for~$|n| \leq \kappa |\alpha|^\frac12$ for some small enough~$\kappa$, the term $\displaystyle \frac{\dd \varphi_n }{\dd r} (1)  + |n|\varphi_n(1) $ will be shown to be nonzero and actually large, due to a specific boundary layer structure of~$\varphi_n $.  
  Combining \eqref{defphitilden} with \eqref{def.psi^noslip} gives the formula 
$$\psi_n (r) = \psi_n^{\rm slip} (r) + a_n r^{-|n|} + b_n \varphi_{n,\rm BL} (r) + b_n \widetilde \varphi_n (r)\,.$$
Thus, with the notation of Proposition~\ref{decompositionvelocity}, the remainder velocity $\widetilde v_n$ is given in terms of the stream function $\widetilde \psi_n$ defined as
$$
 \widetilde \psi_n (r):=  b_n\widetilde\varphi_n(r)\,. $$
\noindent We now focus on the construction of $\varphi_{n}$ and its associate velocity field.  To estimate a possible boundary layer thickness we observe from the first equation of \eqref{eq.fast.slow} that there is a natural scale balance between $-\frac{\dd^2 }{\dd r^2}$ and $-i \alpha n (1-\frac{1}{r^2})\approx -2i\alpha n  (r-1)$ near the boundary~$r=1$, which formally implies that the thickness of the boundary layer is $|2\alpha n|^{-\frac13} = |\beta|^{-1}$ where~$\beta$ is defined in~(\ref{def.beta}).  Before stating the result leading to the construction of the boundary layer term, let us recall the notation introduced in~(\ref{eq.stream}):
 \begin{equation}\label{defHn}
H_n:=
-\frac {\dd ^2}{\dd r^2} - \frac1r \frac \dd {\dd r}+\frac{n^2}{r^2} 
\,,
\end{equation}
 and let us denote
 \begin{equation}\label{defAn}
A_n  :=H_n-i\alpha n \Big(
1-\frac1{r^2}
\Big)
\end{equation}
 so that the first equation in~(\ref{eq.fast.slow}) translates into~$ A_nH_n  \varphi_{n} = 0$.
The next proposition is the construction of $\varphi_{n,\rm BL}$ in \eqref{defphitilden}, which describes the leading part of the boundary layer.
\begin{proposition}\label{propboundarylayer}
There exist $\kappa\in (0,1)$ and $C>0$ such that the following statement holds.
  If $|n| \leq \kappa |\alpha|^\frac12$ and if $|\beta|=(2|\alpha n|)^\frac13 $ is large enough, then there exist smooth functions~$\varphi_{ n,\rm BL}$ and $g_{n,\rm BL}$ on $[1,\infty)$ such that 
\begin{align}\label{eqonphibl}
\begin{split}
&A_n  H_n  \varphi_{n,\rm BL}= \frac1r \frac{\dd}{\dd r} (r g_{n,{\rm BL}})\,,  \quad \mbox{with} \quad \|g_{n,{\rm BL}}\|_{L^2(\Omega)} \le C |\beta|^{\frac32} \, ,   |\varphi_{ n,\rm BL}(1) |\leq 1\,, 
\end{split}
\end{align}
and
\begin{itemize}
\item there holds
\begin{align}\label{lowerbounddphidr}
\Big|  \frac{\dd \varphi_{ n,\rm BL} }{\dd r} (1)  + |n|\varphi_{ n,\rm BL}(1) \Big|  \geq \frac\kappa2|\beta|\,,
\end{align}
\item there is  a smooth function~$   G_{n,\alpha} $ decaying exponentially at infinity uniformly in~$n$ and~$\alpha$ such that
\begin{align}\label{def.varphi.BL}
\varphi_{n,{ \rm BL}}(r)=     G_{n,\alpha} \big( |\beta| (r-1) \big) \,, \qquad r\geq 1\,.
\end{align}
\end{itemize}
\end{proposition}
More precisely, $G_{n,\alpha}$ in \eqref{def.varphi.BL} satisfies the estimate \eqref{nondegenerate.6} stated below.
Let us postpone the proof of Proposition \ref{propboundarylayer} and conclude the proof of Proposition~\ref{decompositionvelocity}. We construct a couple~$(\varphi_{ n,\rm BL},g_{ n,\rm BL})$ as in Proposition~\ref{propboundarylayer}, which produces a boundary layer vector field
\begin{align}\label{def.v_nBL}
v_{r,n,{\rm BL}}(r) :=\frac{inb_n}r  G_{n,\alpha} \big( |\beta| (r-1) \big) \, , \quad v_{\theta,n,{\rm BL}}(r) :=-|\beta| b_n  G_{n,\alpha}' \big( |\beta| (r-1) \big)\,.
\end{align}
Next  we fix~$ \varphi_{n} (1) = \varphi_{n,{ \rm BL}} (1)$, so that the remainder $\widetilde  \varphi_{n}$ as in \eqref{defphitilden} vanishes at the boundary.
 Moreover,  from the requirement~$ A_nH_n  \varphi_{n} = 0$ and thanks to \eqref{eqonphibl}, $\widetilde\varphi_{n }$ is obtained as the solution to 
$$
A_nH_n\widetilde\varphi_{n } = -\frac1r\frac{\dd }{\dd r}  (r g_{n,{\rm BL}})  \qquad r>1 \,, \qquad (H_n \widetilde \varphi_n ) (1) = \widetilde \varphi_n (1) =0  \,.
$$
More precisely, the construction of $\widetilde \varphi_n$ is as follows:
we first construct $\widetilde \omega_n (r) e^{in \theta}$ as the solution to \eqref{eqw} with right-hand side~$-{\rm rot}\, \big (g_{n,{\rm BL}} e^{in\theta} {\bf e}_\theta\big )$. Then $\widetilde \varphi_n$ is obtained as the solution to $H_n \widetilde \varphi_n =\widetilde \omega_n$ under the boundary condition $\widetilde \varphi_n (1)=0$. Let us denote by $\widetilde w_n$ the velocity field whose stream function is $\widetilde \varphi_n$.
Then the estimates of $\widetilde w_n$ and $\widetilde \omega_n$ follow from Proposition \ref{prop.mS_alpha} and the fact that~$\| g_{n, \rm BL}\|_{L^2 (\Omega)} \leq C |\beta|^{\frac32}$ as stated in \eqref{eqonphibl}.
In particular,~\eqref{est.prop.mS_alpha.7} in Proposition~\ref{prop.mS_alpha}  produces  from~$|\beta|=|2\alpha n|^\frac13$ that
$$
\Big|\frac{\dd \widetilde\varphi_{n } }{\dd r}(1)\Big| \leq \| \widetilde w_{\theta,n}\|_{L^\infty (\Omega)} \lesssim |\beta|^{-\frac32} \| g_{n,\rm BL}\|_{L^2(\Omega)}  \leq C\,.
$$
Together with \eqref{defphitilden} and (\ref{lowerbounddphidr}), we find that
$$
\begin{aligned}
\Big| \displaystyle\ \frac{\dd \varphi_n }{\dd r} (1)  + |n|\varphi_n(1) 
\Big| & \geq 
\Big| \displaystyle\ \frac{\dd \varphi_{n,{ \rm BL}} }{\dd r} (1)  + |n|\varphi_{n,{ \rm BL}}(1) 
\Big| - \Big| \displaystyle\ \frac{\dd \widetilde \varphi_{n } }{\dd r} (1)  
\Big|\\
 & \geq 
\frac\kappa2 |\beta| - C
\\
 & \geq 
\frac\kappa4 |\beta| 
\end{aligned}
$$
for~$|\beta|\gg \kappa^{-1}$. 
Let us estimate the coefficients $a_n$ and $b_n$, which are defined in \eqref{defbn}. Since
$$
|v_{\theta,n}^{\rm slip} (1)| \leq \|v_{\theta,n}^{\rm slip} \|_{L^\infty (\Omega)} \lesssim |\alpha n|^{-\frac12} \|f_{n}\|_{L^2(\Omega)}
$$
thanks to Proposition~\ref{prop.mS_alpha} we infer that the parameter~$b_n$ satisfies
$$
|b_n| = \big |\frac{v_{\theta,n}^{\rm slip} (1)}{\frac{\dd \varphi_n }{\dd r} (1)  + |n|\varphi_n(1)} \big |\lesssim |\alpha n|^{-\frac56}\|f_{n}\|_{L^2(\Omega)}
$$
while since~$|\varphi_n(1) | \leq 1$,
$$
|a_n| = |b_n\varphi_n(1)| \lesssim |\alpha n|^{-\frac56}\|f_{n}\|_{L^2(\Omega)}\,.
$$
  Thus, the estimates of the boundary layer velocity $v_{n,\rm BL}$ easily follow from its definition in~\eqref{def.v_nBL}. The estimates of the remainder velocity $\widetilde v_n$ follow from Proposition \ref{prop.mS_alpha} and the estimates of $b_n$ since $\widetilde \psi_n = b_n \widetilde \varphi_n$ that implies $\widetilde v_n = b_n \widetilde w_n$. This concludes the proof of Proposition~\ref{decompositionvelocity}.
\end{proofprop}

\begin{proofprop}{propboundarylayer} Without loss of generality we assume from now on that~$\alpha>0$.
As already mentioned, we formally estimate the thickness of the boundary layer to be of the order $|\beta|=|2\alpha n |^\frac13$. One important remark here is the size of~$\frac{n^2}{r^2}$ in the operators $H_n$ and~$A_n$ defined in~(\ref{defHn}) and~(\ref{defAn}). Recall that we are interested in the regime $|n|\leq O (\alpha^\frac12)$. If~$|n|=O(\alpha^\frac12)$ then we observe that $|\beta|=|2\alpha n|^\frac13=O((\alpha^\frac32)^\frac13) = O(|n|)$, and thus, the term $\frac{n^2}{r^2}$ has the same size near the boundary as $\partial_r^2$ and $\alpha n (r-1)$.
Hence, in the construction of the boundary layer we also need to take into account the term $\frac{n^2}{r^2}$, for this term is no longer small in the regime $|n|=O(\alpha^\frac12)$.
With this remark in mind let us rewrite $H_n$ and $A_n$ in more convenient forms: we define
$$
\widetilde H_n  =
-\frac {\dd^2}{\dd r^2} +n^2  \quad \mbox{and} \quad 
\widetilde A_n =\widetilde H_n - 2i\alpha n (r-1) =  -\frac {\dd^2}{\dd r^2}  +\beta^3 \big ( r-1 +\frac{i n}{2\alpha} \big ) \,,
$$
so that
\begin{align*}
H_n =  \widetilde H_n   - \frac1r\frac{\dd }{\dd r} + n^2 (\frac{1}{r^2}-1 ) 
\end{align*}
and by using $1-\frac{1}{r^2} = 2(r-1) + (r-1)^2 \frac{1+2r}{r^2}$,
\begin{align*}
A_n & =  H_n - 2i \alpha n (r-1) - i\alpha n (r-1)^2 \frac{1+2 r}{r^2} \\
& = \widetilde H_n - 2i \alpha n(r-1) - \frac1r \frac{\dd }{\dd r} + n^2 (\frac{1}{r^2}-1 )  - i\alpha n (r-1)^2 \frac{1+2 r}{r^2} \\
& = \widetilde A_n - \frac1r\frac{\dd }{\dd r} + n^2 (\frac{1}{r^2}-1 )  - 2i\alpha n (r-1)^2 \frac{1+2 r}{r^2} \, \cdotp
\end{align*}
That is, $\widetilde H_n$ and $\widetilde A_n$ have the leading size of $H_n$ and $A_n$, respectively, for the boundary layer functions.
Then we write $A_n H_n$ as 
\begin{align}\label{eq.A_nH_n}
A_n H_n = \widetilde A_n \widetilde H_n + \big (A_n- \widetilde A_n \big )H_n  +\widetilde A_n \big ( H_n - \widetilde H_n\big )\,,
\end{align}
and we claim that
\begin{align}\label{eq.remainder.rot}
\big (A_n- \widetilde A_n \big )H_n  +\widetilde A_n \big ( H_n - \widetilde H_n\big ) = {\rm rot}_n\, R_n 
\end{align}
with a suitable operator $R_n$, where ${\rm rot}_n$ is defined in   polar coordinates with the $n$ mode for the angular variable. Note that  the leading operator $\widetilde A_n \widetilde H_n$, when applied to a boundary layer function of the form $h(|\beta| (r-1))$, formally has size~$O(|\beta|^4)$,  the term $\big (A_n- \widetilde A_n \big )H_n  + \widetilde A_n \big ( H_n - \widetilde H_n\big )$ has size  $O(|\beta|^3)$, and  then $R_n$ is of size~$O(|\beta|^2)$.

\smallskip
\noindent
Let us look for the boundary layer $\varphi_{n,{ \rm BL}}$ as a solution to
\begin{align}\label{eqboundarylayer}
\widetilde A_n \widetilde H_n \varphi_{n,{ \rm BL}} = 0\,.
\end{align}
Since $\widetilde H_n$ is easily inverted for $|n|\geq 1$, we start by considering the homogeneous problem~$\widetilde A_n \phi=0$.
By its very definition the operator $\widetilde A_n$ is nothing but the Airy operator with a complex coefficient.
Hence we introduce the Airy function $\Ai (z)$ which is a solution to 
$$
\frac{\dd^2 \Ai}{\dd z^2}  - z \Ai =0
$$ in $\C$; for details, see   Appendix~\ref{appendixairy}.
Then we define 
\begin{align}\label{def.G_beta}
\widetilde G_{n,\alpha} (\rho) := \Ai \big (c_- (\rho +\frac{in |\beta|}{2\alpha}) \big ) \,, \qquad \rho>0\,,
\end{align} 
which satisfies from $c_-^3=-i$,
\begin{equation}\label{eqAirymod}
\bigg ( \frac{\dd^2}{\dd \rho^2}   +i  (\rho + \frac{in}{2\alpha} |\beta|)  \bigg)\widetilde    G_{n,\alpha} = 0\,, \qquad \rho>0\,.
\end{equation}
Next we set 
\begin{align}\label{def.tildeG_beta}
  G_{n,\alpha} (\rho): =-  \int_\rho^\infty e^{-\frac{|n|}{|\beta|} (\rho -\tau )} \int_\tau^\infty e^{-\frac{|n|}{|\beta|}(\sigma-\tau)}\widetilde  G_{n,\alpha} ( \sigma) \dd \sigma \dd\tau\,,
\end{align}
which satisfies 
\begin{align}
- \frac{\dd ^2G_{n,\alpha}}{\dd \rho^2}   + \frac{n^2}{|\beta|^2}   G_{n,\alpha} =  \widetilde G_{n,\alpha}\,, \qquad \rho>0\,.
\end{align}
Finally we define
\begin{align}\label{def.C_0.n}
C_{0,n,\alpha}  =
\begin{cases}
&\displaystyle\frac1{G_{n,\alpha}(0)} \quad \text{ if }\quad  |G_{n,\alpha}(0)| \geq 1 \,,\\
& 1\qquad \text{ otherwise}
\end{cases}
\end{align}
and we
set 
\begin{align}\label{defvarphinBL}
\varphi_{n,{\rm BL}} (r) := C_{0,n,\alpha}  G_{n,\alpha} \big(|\beta| (r-1)\big)\,,
\end{align}
which satisfies from $-i |\beta|^3 = \beta^3 $,
\begin{align*}
\widetilde A_n \widetilde H_n \varphi_{n,{\rm BL}} =0\,, \qquad r>1\,,  \end{align*}
as desired.
Notice also that
$$
|\varphi_{n,{\rm BL}} (1) | \leq 1\,.
$$
The key quantity is  \begin{align*} 
\frac{\dd \varphi_{n,{\rm BL}}}{\dd r}|_{r=1} & = C_{0,n,\alpha} |\beta| \, \frac{\dd   G_{n,\alpha}}{\dd \rho} |_{\rho=0}  \\
& =C_{0,n,\alpha} |\beta| \Big ( \int_0^\infty e^{-\frac{|n|}{|\beta|}\sigma} \widetilde G_{n,\alpha} (\sigma) \dd \sigma  - \frac{|n|}{|\beta|} G_{n,\alpha} (0) \Big )\nonumber \\
& = C_{0,n,\alpha} |\beta| \Big ( \int_0^\infty e^{-\frac{|n|}{|\beta|}\sigma} \Ai \big (c_- (\sigma   + \frac{i n |\beta|}{2\alpha}) \big ) \dd \sigma  - \frac{|n|}{|\beta|} G_{n,\alpha} (0) \Big )\nonumber  \\
& = C_{0,n,\alpha} |\beta| \Big ( \frac{1}{c_-} \int_0^\infty e^{-\lambda s} \Ai  (s+ \frac{i  n |\beta| c_-}{2\alpha})  \dd s - \frac{|n|}{|\beta|} G_{n,\alpha} (0) \Big )\,, 
\end{align*}
with
$$
\lambda=\lambda_{n,\beta} := \frac{|n|}{|\beta|c_-}=\frac{|n|c_+}{|\beta|} \,\cdotp
$$
 We find from $\beta^3=|\beta|^3 c_-^3$ that 
\begin{align}\label{deflambda2}
\frac{in|\beta|c_-}{2\alpha} = \frac{n^2|\beta| c_-}{-2i\alpha n} = \frac{n^2|\beta|c_-}{\beta^3} = \frac{n^2|\beta|}{|\beta|^3 c_-^2} = \frac{n^2 c_+^2}{|\beta|^2} = \lambda^2
\end{align}
so
\begin{equation}\label{keyquantity}
\frac{\dd \varphi_{n,{\rm BL}}}{\dd r}|_{r=1} = C_{0,n,\alpha} |\beta| \Big ( \frac{1}{c_-} \int_0^\infty e^{-\lambda s}  \Ai (s + \lambda^2 ) \dd s - \frac{|n|}{|\beta|} G_{n,\alpha} (0) \Big )\,.
\end{equation}
Note that $ |\lambda|^2=4^{-1} |n|^{\frac43} \alpha^{-\frac23}$ is small when $|n|\ll \alpha^\frac12$.  
The proof of the following lemma is postponed to   Appendix~\ref{appendixairy}.
\begin{lemma}\label{constructionGbeta}
There holds
\begin{align}\label{nondegenerate.5} 
\widetilde C_0 : = \inf \Big \{  ~ |C_{0,n,\alpha}| ~~ | ~~\alpha \geq 1\,, ~ 1\leq |n|\leq \alpha^\frac12 \Big \} >0  \,,
\end{align}
and there is a constant~$\varepsilon\in (0,1)$ such that defining
$$
 \Sigma_\varepsilon:= \Big \{\mu \in \C \, | \,  \arg \, \mu = \frac\pi6 \, , \, 0 \leq |\mu| \leq\varepsilon \Big\}\,,
$$
then  
\begin{align}\label{nondegenerate.4} 
\widetilde \kappa_\varepsilon : = \inf_{\mu \in \Sigma_\varepsilon}\Big |\int_0^\infty e^{-\mu s} \Ai (s + \mu^2 ) \dd s\Big | >0  \,.
\end{align}
Moreover, the function~$G_{n,\alpha}$ defined in~{\rm(\ref{def.tildeG_beta})} 
satisfies for~$R$ large enough
\begin{align}\label{nondegenerate.6} 
\sup_{\alpha\geq 1} \, \sup_{1\leq |n|\leq \alpha^{1/2} } \, \sup_{\rho \geq R} \, e^{\rho} \, |  \frac{\dd ^k G_{n,\alpha}}{\dd \rho^k} (\rho)| <\infty\,, \qquad k=0,1,2,3\,.
\end{align}
\end{lemma}

Now let us return to the proof of Proposition~\ref{propboundarylayer}. We define~$\varphi_{n,{\rm BL}} $ as in~(\ref{defvarphinBL}). Note that thanks   to~(\ref{keyquantity}), \eqref{nondegenerate.4} and~\eqref{nondegenerate.5}  there holds as long as~$\frac{|n|}{|\beta|}=|\lambda| \leq\varepsilon $,
$$
\begin{aligned}
|\frac{\dd \varphi_{n,{\rm BL}}}{\dd r}  (1)| & \geq \widetilde C_0  |\beta| \Big ( \widetilde \kappa_\varepsilon  -  \frac{|n|}{|\beta|} |G_{n,\alpha} (0)| \Big ) \\
& \geq  \widetilde C_0  |\beta| \big ( \widetilde\kappa_\varepsilon  - \frac{\varepsilon}{\widetilde C_0} \big ) \,.
\end{aligned} 
$$
Here we have also used \eqref{def.C_0.n}. Hence
\begin{align}
|\frac{\dd \varphi_{n,{\rm BL}}}{\dd r}  (1)|  \geq \frac{ \widetilde C_0 \widetilde  \kappa_\varepsilon}2|\beta|
\end{align} 
  as long as~$2\varepsilon \leq \widetilde C_0 \widetilde  \kappa_\varepsilon$, which is possible since $\widetilde \kappa_\varepsilon$ is nonincreasing in $\varepsilon>0$.
This proves~(\ref{lowerbounddphidr}) since
$$
\begin{aligned}
\Big| \displaystyle\ \frac{\dd \varphi_{n,{\rm BL}} }{\dd r} (1)  + |n|\varphi_{n,{\rm BL}}(1) 
\Big| & \geq 
\Big| \displaystyle\ \frac{\dd \varphi_{n,{ \rm BL}} }{\dd r} (1)  \Big| -  |n|| \varphi_{n,{ \rm BL}}(1)   | \\
 & \geq   
  \frac{ \widetilde C_0 \widetilde  \kappa_\varepsilon}2  |\beta| - |n| \\
 & \geq  
 \frac {\check\kappa}2  |\beta| 
\end{aligned}
$$
with~$\check\kappa := \widetilde C_0 \widetilde  \kappa_\varepsilon/2$, and the last inequality  holds as long as~$|n|\leq \min(\varepsilon, \check\kappa/2)|\beta|$.
It then suffices to choose~$\kappa \leq \min(\sqrt 2 \varepsilon^\frac32, \check \kappa^\frac32/2)$ which ensures that
$$
|n|\leq \kappa \alpha^\frac12 \Rightarrow |n|\leq  \min ( \varepsilon, \frac{\check\kappa}{2}) \, |\beta|   \, .
$$
The result~(\ref{def.varphi.BL}) is an obvious consequence of the previous construction so  it remains to prove that~(\ref{eqonphibl}) is  satisfied for a suitable $g_{n,{\rm BL}}$. 
Let us recall~(\ref{eq.A_nH_n}). Then we define $g_{n,\rm{BL}}$ as 
\begin{equation}\label{defRntheta}
g_{n,\rm{BL}}(r) :=- \frac1r \int_r^\infty s \Big( \widetilde A_n \big ( H_n - \widetilde H_n\big )  +\big (A_n- \widetilde A_n \big )H_n \Big)\varphi_{n,{\rm BL}}\, \dd s\,,
\end{equation}
which then satisfies 
$$
\Big( \big (A_n- \widetilde A_n \big )H_n  + \widetilde A_n \big ( H_n - \widetilde H_n\big ) \Big)\varphi_{n,{\rm BL}}= \frac1r \frac{\dd }{\dd r} (r \, g_{n,\rm{BL}} )\,.
$$
Let us consider the estimate of $g_{n,{\rm BL}}$.
We compute the highest order terms in~$ \widetilde A_n \big ( H_n - \widetilde H_n\big )$ and~$\widetilde A_n \big ( H_n - \widetilde H_n\big )$, which are of order~$O(|\beta|^3)$ when~$n = O(\alpha^\frac12)$:
$$
\begin{aligned}
\frac1{ |\beta|^3} \widetilde A_n \big ( H_n - \widetilde H_n\big )\varphi_{n,{\rm BL}}& = \frac1r    G_{n,\alpha}''' +
\frac{n^2}{|\beta|}   \frac{r^2-1}{r^2}   G_{n,\alpha}''- \frac{n^2}{r|\beta|^2} 
  G_{n,\alpha}'
-\frac{n^4}{|\beta|^3 r^2} (r^2-1)   G_{n,\alpha} \\
&\quad + \frac{2i\alpha n}{|\beta|^2 r}(r-1)
  G_{n,\alpha}' +  \frac{2i\alpha n^3}{|\beta|^3 r^2}(r-1)^2 (r+1)
   G_{n,\alpha}
 + \rm{l.o.t}
 \end{aligned}
$$
where all the~$G_{n,\alpha}$ are computed at~$|\beta|(r-1)$ and the lower order terms are to be understood in terms of~$|\beta|$ for~$n = O(\alpha^\frac12)$. It is clear from this formula that as long as~$n^2 \lesssim \alpha$ then there is a function~$R_{n,\alpha,I} (r,\rho)$, uniformly bounded in~$r$ and exponentially decaying at infinity in~$\rho$,
such that
$$
\widetilde A_n \big ( H_n - \widetilde H_n\big )\varphi_{n,{\rm BL}} =  |\beta|^3   R_{n,\alpha,I} \big(r,|\beta|(r-1)\big)\,.
$$
Similarly
one has
$$
\begin{aligned}
\frac1{ |\beta|^3}\big (A_n- \widetilde A_n \big )H_n \varphi_{n,{\rm BL}}& = 
\frac1r  G_{n,\alpha}''' - \frac{n^2}{|\beta|r^2} (1-r^2)
  G_{n,\alpha}''+\frac{2 i\alpha n }{|\beta|r^2} (r-1)^2(1+2r)   G_{n,\alpha}''\\&\quad 
- \frac{n^2}{r|\beta|^2}  G_{n,\alpha}'+\frac{2i\alpha n^3}{|\beta|^2 r^2} (r-1)^2(1+2r)   G_{n,\alpha}
+\rm{l.o.t}
 \end{aligned}
$$
so
again
$$
\big (A_n- \widetilde A_n \big )H_n \varphi_{n,{\rm BL}}= |\beta|^3  R_{n,\alpha,II} \big(r,|\beta|(r-1)\big)
$$
with some $R_{n,\alpha,II} (r,\rho)$ which is uniformly bounded in $r$ and exponentially decaying as $\rho\rightarrow \infty$. This implies that~$g_{n,\rm{BL}}$ defined in~(\ref{defRntheta}) satisfies
$$
\|g_{n,\rm{BL}}\|_{L^2(\Omega)} \lesssim |\beta|^\frac32\,.
$$
The result~(\ref{eqonphibl}) follows.
\medskip
\noindent
Proposition~\ref{propboundarylayer} is proved.
\end{proofprop}

\section{The nonlinear problem}\label{sec.nonlinear}
In this section we construct the solution to the nonlinear problem \eqref{tildeNS_alpha}
\begin{equation}\tag{$\widetilde {\rm NS}_{\alpha}$}\label{tildeNS_alpha}
  \left\{
\begin{aligned}
  -\Delta v - \alpha ( x^\bot \cdot \nabla v - v^\bot ) + \nabla q + \alpha U^\bot {\rm rot}\, v & \,=\, 
  - v^\bot {\rm rot}\, v   + f \,,   \qquad x\in \Omega \,,\\
{\rm div}\, v &  \,=\, 0\,,   \qquad \qquad  x \in \Omega\,, \\
  v & \,=\, 0 \,, \qquad \qquad  x \in \partial \Omega\,, \\
\end{aligned}\right.
\end{equation}
in the class $v=(\mathcal{P}_{0} v, \mathcal{Q}_{0} v) \in X = \mathcal{P}_0 W_0^{1,\infty} (\Omega)^2 \times \mathcal{Q}_0 W^{1,2}_0 (\Omega)^2$ with a suitable pressure~$q\in W^{1,1}_{loc}(\overline{\Omega})$. As   already noted in the introduction, for the solvability of \eqref{tildeNS_alpha}, the key observation is the  decomposition of the nonlinear term $v^\bot {\rm rot}\, v$: we have for $v \in X$,
\begin{equation*}
\begin{split}
v^\bot {\rm rot}\,v
&\,=\,  (\mathcal{P}_{0} v)^\bot {\rm rot}\,\mathcal{Q}_{0}v
+ (\mathcal{Q}_{0}v )^\bot {\rm rot}\,\mathcal{P}_{0} v
+ (\mathcal{Q}_{0} v)^\bot {\rm rot}\, \mathcal{Q}_{0} v 
+ (\mathcal{P}_{0} v)^\bot {\rm rot}\,\mathcal{P}_{0} v\\
 & \, = \, \mathcal{Q}_{0}\Big((\mathcal{P}_{0} v)^\bot {\rm rot}\,\mathcal{Q}_{0} v
+ (\mathcal{Q}_{0}v)^\bot {\rm rot}\,\mathcal{P}_{0} v \Big )
+ (\mathcal{Q}_{0} v)^\bot {\rm rot}\, \mathcal{Q}_{0} v 
+ (\mathcal{P}_{0} v)^\bot {\rm rot}\,\mathcal{P}_{0} v\,.
\end{split}
\end{equation*}
  Here we have used $\mathcal{P}_0 \big ( (\mathcal{P}_{0} v)^\bot {\rm rot}\,\mathcal{Q}_{0}v
+ (\mathcal{Q}_{0}v )^\bot {\rm rot}\,\mathcal{P}_{0} v \big ) =0$.
Furthermore, since the last term on the right-hand side can be written in a gradient form, the problem~\eqref{tildeNS_alpha} is in fact reduced to the next system
\begin{equation}\tag{$\widehat {\rm NS}_{\alpha}$}\label{hatNS_alpha}
  \left\{
\begin{aligned}
  -\Delta v - \alpha ( x^\bot \cdot \nabla v - v^\bot ) + \nabla \tilde{q} + \alpha U^\bot {\rm rot}\, v & \,=\, 
  G(v) + f \,,   \qquad x\in \Omega \,,\\
{\rm div}\, v &  \,=\, 0\,,   \qquad \qquad  x \in \Omega\,, \\
  v & \,=\, 0 \,, \qquad \qquad  x \in \partial \Omega\,. 
\end{aligned}\right.
\end{equation}
Here we have set
\begin{equation}\label{proof.nonlinear.nonlinearity.1}
G(v) 
\,:=\,  
-\mathcal{Q}_{0}\Big((\mathcal{P}_{0} v)^\bot {\rm rot}\,\mathcal{Q}_{0} v
+ (\mathcal{Q}_{0}v)^\bot {\rm rot}\,\mathcal{P}_{0} v \Big )
- (\mathcal{Q}_{0} v)^\bot {\rm rot}\, \mathcal{Q}_{0} v\,,
\end{equation}
and $\tilde{q} := q + Q$, where $Q=Q(|x|)$ is a radial function satisfying $\nabla Q = - (\mathcal{P}_{0} v)^\bot {\rm rot}\,\mathcal{P}_{0} v$.

Our aim is to prove the  existence and uniqueness of solutions $(v, q)$ for \eqref{hatNS_alpha}  in a suitable subset of $X$, under some conditions on the external force $f =(\mathcal{P}_{0} f, \mathcal{Q}_{0} f) $ in~$ Y= \mathcal{P}_{0} L^{1} (\Omega)^2 \times \mathcal{Q}_{0} L^2 (\Omega)^2$. The proofs in Subsections \ref{subsec.nonlinear.general.alpha} and \ref{subsec.nonlinear.large.alpha} rely on the standard Banach fixed point argument, where the estimate of the nonlinearity $G(v)$ in the space $Y$  is important. Thanks to the identity
$$
 \mathcal{P}_{0} G (v) = - \mathcal{P}_0 \Big ( (\mathcal{Q}_{0} v)^\bot {\rm rot}\, \mathcal{Q}_{0} v \Big )  \, , 
$$
we see that $\mathcal{P}_{0} G(v)$ belongs to $L^1 (\Omega)^{2}$, which is the same summability as the space $Y$. In order to control the $L^{2}$-norm of $\mathcal{Q}_{0} G(v)$ in the iteration scheme, we introduce the closed subspace $X_{0}$ of $X$ equipped with the norm $\|\cdot\|_{X_{0}}$:
\begin{equation}\label{proof.nonlinear.space}
\begin{split}
X_{0}
\,:=\,
\bigg\{ &v \in X ~\bigg|~~~~
{\rm div}\, \mathcal{P}_0 v \,=\,{\rm div}\, \mathcal{Q}_0 v \,=\,0 ~~\text{in}~\Omega\,,  \\
& \|v\|_{X_{0}} \,:=\,
\|\mathcal{P}_{0} v\|_{L^{\infty}(\Omega)}
+ \|\nabla \mathcal{P}_{0} v\|_{L^{\infty}(\Omega)} \\
& \qquad\quad
+ \|\mathcal{Q}_{0} v\|_{L^{2}(\Omega)} 
+ \|\nabla \mathcal{Q}_{0} v\|_{L^{2}(\Omega)} 
+ \sum_{|n| \ge 1}
\|\mathcal{P}_{n} v\|_{L^{\infty}(\Omega)} 
< \infty
\bigg\}\,.
\end{split}
\end{equation}
Indeed, we can easily establish an a priori estimate for 
$$
 \mathcal{Q}_{0} G(v)= - \mathcal{Q}_{0}\Big( (\mathcal{P}_{0} v)^\bot {\rm rot}\,\mathcal{Q}_{0} v
+ (\mathcal{Q}_{0}v)^\bot {\rm rot}\,\mathcal{P}_{0} v + (\mathcal{Q}_{0} v)^\bot {\rm rot}\, \mathcal{Q}_{0} v\Big)$$
in~$L^{2}(\Omega)^{2}$ for $v \in X_{0}$; see \eqref{proof.nonlinear.3} below for the estimate of $\big\| \mathcal{Q}_{0} \big( (\mathcal{Q}_{0} v)^\bot {\rm rot}\, \mathcal{Q}_{0} v\big)\big\|_{L^2(\Omega)}$.

 After proving Theorem \ref{thm.main.1} and \ref{thm.main.2}, we revisit Theorem \ref{thm.main.2} in order to study the qualitative behavior of solutions. By fixing the external force $f \in Y$, we consider the fast rotation limit $|\alpha| \to \infty$ for the solution $(v, q)=(v^{(\alpha)}, q^{(\alpha)})$ to~\eqref{tildeNS_alpha}. The results are summarized in Theorem \ref{thm.main.3}, which will be proved in Subsection \ref{subsec.nonlinear.thm.1.3}.

\subsection{Useful estimates on \eqref{S_alpha}}
Before we give the proofs of Theorems \ref{thm.main.1} and   \ref{thm.main.2}, 
let us restate the main estimates for the linearized problem \eqref{S_alpha} for general $\alpha$, which will be used throughout this section, as well as some specific estimates corresponding to~$|\alpha| \gg 1$. 

\subsubsection{The case of general $\alpha$}
Let $\alpha \in \R\setminus\{0\}$ and $f \in Y$, and let~$(v,q) \in X \times W^{1,1}_{loc}(\overline{\Omega})$ be the unique solution to \eqref{S_alpha} given in Proposition \ref{prop.all.alpha}. For   notational simplicity, we denote by $v^{{\rm low}}$ and $v^{{\rm high}}$, respectively, the low and high frequency parts  of $\mathcal{Q}_{0} v$:
\begin{equation}\label{defvlowvhigh}
\begin{split}
v^{{\rm low}} \,:=\, \sum_{1 \le |n| <   1+ \sqrt{2|\alpha|}} \mathcal{P}_n v\,, 
\qquad\quad
v^{{\rm high}} \,:=\, \sum_{|n| \ge 1+ \sqrt{2|\alpha|}} \mathcal{P}_n v\,.
\end{split}
\end{equation}
From the estimates \eqref{est.prop.all.alpha.1} and \eqref{est.prop.all.alpha.4} in Proposition \ref{prop.all.alpha}, we have 
\begin{align}
\|v^{{\rm low}}\|_{L^{2}(\Omega)}
&\le
\frac{C}{|\alpha|^{\frac12}} 
\big(1 + \frac{1}{|\alpha|^{\frac12}} \big)
\|\mathcal{Q}_{0} f\|_{L^{2}(\Omega)}\,, 
\label{est.proof.nonlinear.1} \\
\|v^{{\rm high}}\|_{L^{2}(\Omega)}
&\le
\frac{C}{|\alpha|} 
\|\mathcal{Q}_{0} f\|_{L^{2}(\Omega)}\,, 
\label{est.proof.nonlinear.2}
\end{align}
and from \eqref{est.prop.all.alpha.3} and \eqref{est.prop.all.alpha.6} in the same proposition, we have the estimates for the derivatives.
\begin{align}
\|\nabla v^{{\rm low}}\|_{L^{2}(\Omega)}
& \le
C \big(1 + \frac{1}{|\alpha|^{\frac12}} \big)
\|\mathcal{Q}_{0} f\|_{L^{2}(\Omega)}\,,
\label{est.proof.nonlinear.3} \\
\|\nabla v^{{\rm high}}\|_{L^{2}(\Omega)}
&\le
\frac{C}{|\alpha|^{\frac12}}
\|\mathcal{Q}_{0} f\|_{L^{2}(\Omega)}\,.
\label{est.proof.nonlinear.4}
\end{align}
For the fixed point argument in the space $X_{0}$, the estimate of $\sum_{|n| \ge 1} 
\|\mathcal{P}_{n} v\|_{L^{\infty}(\Omega)}$ is also needed. In the case $0<|\alpha|<1$, \eqref{est.prop.all.alpha.1'} and \eqref{est.prop.all.alpha.4'} in Proposition \ref{prop.all.alpha}, and the H\"{o}lder inequality for sequences lead to
\begin{equation*}
\begin{split}
\sum_{|n| \ge 1} 
\|\mathcal{P}_{n} v\|_{L^{\infty}(\Omega)} 
& \le
\frac{C}{|\alpha|^{\frac34}} 
\sum_{|n| \ge 1} 
\frac{1}{|n|^{\frac34}}
\|\mathcal{P}_{n} f\|_{L^{2}(\Omega)} \\
& \le
\frac{C}{|\alpha|^{\frac34}}
\big( \sum_{|n| \ge 1}  \frac{1}{|n|^{\frac32}} \big)^{\frac12}\,
\|\mathcal{Q}_{0} f\|_{L^{2}(\Omega)} \,,
\end{split}
\end{equation*}
which implies 
\begin{align}
\sum_{|n| \ge 1} 
\|\mathcal{P}_{n} v\|_{L^{\infty}(\Omega)} 
\le
\frac{C}{|\alpha|^{\frac34}}
\|\mathcal{Q}_{0} f\|_{L^{2}(\Omega)} \,.
\label{est.proof.nonlinear.5}
\end{align}
In the case $|\alpha| \ge 1$ we have from \eqref{est.prop.all.alpha.1'}, 
\begin{align}
\sum_{1 \le |n| < 1+ \sqrt{2|\alpha|}} 
\|\mathcal{P}_{n} v\|_{L^{\infty}(\Omega)} 
\le
\frac{C}{|\alpha|^{\frac14}} 
\|\mathcal{Q}_{0} f\|_{L^{2}(\Omega)} \,,
\label{est.proof.nonlinear.6}
\end{align}
and from \eqref{est.prop.all.alpha.4'} with the H\"{o}lder inequality,
\begin{align}
\sum_{|n| \ge 1+ \sqrt{2|\alpha|}} 
\|\mathcal{P}_{n} v\|_{L^{\infty}(\Omega)} 
\le
\frac{C}{|\alpha|^{\frac12}} 
\|\mathcal{Q}_{0} f\|_{L^{2}(\Omega)} \,.
\label{est.proof.nonlinear.7}
\end{align}

\subsubsection{The case of large $|\alpha|$}
Now let us consider   the special case when~$|\alpha|\gg1$. Thanks to the boundary layer analysis in Proposition \ref{decompositionvelocity}, we have a better estimate for the linearized problem \eqref{S_alpha} in terms of   decay in the parameter $|\alpha n|$. Given the parameter~$\kappa$ defined in Proposition~\ref{decompositionvelocity} we assume that~$|\alpha|$ is large enough, as required in  Proposition \ref{decompositionvelocity},  and let us truncate frequencies, similarly to~(\ref{defvlowvhigh}), as follows:
\begin{equation}\label{defvlowvhigkappah}
\begin{split}
v^{{\rm low},\kappa} \,:=\, \sum_{ 1\leq  |n| \leq \kappa |\alpha|^\frac12} \mathcal{P}_n v\,, 
\qquad\quad
v^{{\rm high},\kappa} \,:=\, \sum_{|n|  > \kappa |\alpha|^\frac12 } \mathcal{P}_n v\,.
\end{split}
\end{equation}
Then thanks to~(\ref{est.decompositionvelocity}) we know that
\begin{align}
\|v^{{\rm low},\kappa}\|_{L^{2}(\Omega)}
&\le
\frac{C}{|\alpha|^{\frac23}} 
\|\mathcal{Q}_{0} f\|_{L^{2}(\Omega)}\,,
\label{est.proof.nonlinear.1''} \\
\|\nabla v^{{\rm low},\kappa}\|_{L^{2}(\Omega)}
& \le
\frac{C}{|\alpha|^{\frac13}} 
\|\mathcal{Q}_{0} f\|_{L^{2}(\Omega)}\,,
\label{est.proof.nonlinear.3'} 
\end{align}
and
\begin{align}
\sum_{1 \le |n| \le  \kappa |\alpha|^\frac12} 
\|\mathcal{P}_{n} v\|_{L^{\infty}(\Omega)} 
& \le
\frac{C}{|\alpha|^{\frac12}} 
\big( \sum_{1 \le |n| \le \kappa |\alpha|^\frac12}
\frac{1}{|n|} \big)^{\frac12}
\|\mathcal{Q}_{0} f\|_{L^{2}(\Omega)} \nonumber \\
& \le
\frac{C(\log |\alpha|)^{\frac12}}{|\alpha|^{\frac12}}
\|\mathcal{Q}_{0} f\|_{L^{2}(\Omega)}\,.
\label{est.proof.nonlinear.6'} 
\end{align}
Note that the constant~$C $ depends on~$\kappa$, which is fixed. 
One can see that the decay in terms of~$|\alpha|$ compared with~\eqref{est.proof.nonlinear.1}, \eqref{est.proof.nonlinear.3} and \eqref{est.proof.nonlinear.6} are improved. 

\medskip

\noindent Similarly thanks to~(\ref{est.prop.all.alpha.4}),~(\ref{est.prop.all.alpha.6}),~(\ref{est.prop.large.alpha.1}) and~(\ref{est.prop.large.alpha.3}) there holds
\begin{align}
\|v^{{\rm high},\kappa}\|_{L^{2}(\Omega)}
&\le
\frac{C}{|\alpha| } 
\|\mathcal{Q}_{0} f\|_{L^{2}(\Omega)}\,,
\label{est.proof.nonlinear.1'} \\
\|\nabla v^{{\rm high},\kappa}\|_{L^{2}(\Omega)}
& \le
\frac{C}{|\alpha|^{\frac12}} 
\|\mathcal{Q}_{0} f\|_{L^{2}(\Omega)}\,.
\label{est.proof.nonlinear.3'} 
\end{align}
Finally from~(\ref{est.prop.all.alpha.4'}),~(\ref{est.prop.large.alpha.2}) and the  H\"{o}lder inequality, we derive
\begin{align}
\sum_{|n| \ge \kappa|\alpha|^\frac12} 
\|\mathcal{P}_{n} v\|_{L^{\infty}(\Omega)} 
\le
\frac{C}{|\alpha|^{\frac12}} 
\|\mathcal{Q}_{0} f\|_{L^{2}(\Omega)} \,.
\label{est.proof.nonlinear.7'}
\end{align}

%
\subsection{Proof for general $\alpha$}\label{subsec.nonlinear.general.alpha}
%
In this subsection we prove Theorem \ref{thm.main.1}, by means of a fixed point argument. The solutions to \eqref{hatNS_alpha} will be found in the closed convex set $\mathcal{B}_{\vec{\delta},\epsilon}$ of $X_{0}$ defined as follows:
\begin{equation*}
\begin{split}
&\mathcal{B}_{\vec{\delta},\epsilon} 
\,:=\, \mathcal{B}_{(\delta_{1}, \delta_{2}, \delta_{3}, \delta_{4}), \epsilon} \\
&\,:=\,
\bigg\{ h \in X_{0} ~\bigg|~ 
\|\mathcal{P}_{0} h\|_{L^{\infty}(\Omega)}
+ \|\nabla \mathcal{P}_{0} h\|_{L^{\infty}(\Omega)}
\le \epsilon |\alpha|^{\delta_{1}} \,,\\
& \quad\quad\quad
\|\mathcal{Q}_{0} h\|_{L^{2}(\Omega)} 
\le \epsilon |\alpha|^{\delta_{2}}\,,\quad
\|\nabla \mathcal{Q}_{0} h\|_{L^{2}(\Omega)}
\le \epsilon |\alpha|^{\delta_{3}}\,, \quad
\sum_{|n| \ge 1}
\|\mathcal{P}_{n} h\|_{L^{\infty}(\Omega)} 
\le \epsilon |\alpha|^{\delta_{4}}
\bigg\}\,,
\end{split}
\end{equation*}
where we have set $\vec{\delta}:=(\delta_{1}, \delta_{2}, \delta_{3} ,\delta_{4})$, and the numbers $\delta_{1},\ldots,\delta_{4}$ and the positive number~$\epsilon$ will be chosen later. For any $h \in \mathcal{B}_{\vec{\delta},\epsilon}$, let $(v_{h}, q_{h})$ be the unique solution constructed  in Proposition \ref{prop.all.alpha} to the linear system
\begin{equation}\label{proof.nonlinear.linear}
  \left\{
\begin{aligned}
-\Delta v_{h} - \alpha (x^\bot \cdot \nabla v_{h} - v_{h}^\bot ) 
+ \nabla q_{h} 
+ \alpha U^\bot {\rm rot}\,v_{h} 
& \,=\, 
G(h) + f \,, \qquad x\in \Omega \,,\\
{\rm div}\, v_{h} &  \,=\, 0\,,   \qquad \qquad  x \in \Omega\,, \\
v_{h} & \,=\, 0 \,, \qquad \qquad  x \in \partial \Omega\,, \\
\end{aligned}\right.
\end{equation}
where the function $G$ is defined in \eqref{proof.nonlinear.nonlinearity.1}.

Let us start with the estimate of $G(h)$ in the space $Y$. The first two terms in the right-hand side of \eqref{proof.nonlinear.nonlinearity.1} with $v$ replaced by $h$ can be estimated as
\begin{equation}\label{proof.nonlinear.1}
\begin{split}
&~~~\|\mathcal{Q}_{0}\big((\mathcal{P}_{0} h)^\bot {\rm rot}\,\mathcal{Q}_{0} h\big) \|_{L^{2}(\Omega)}
+ \|\mathcal{Q}_{0}\big ((\mathcal{Q}_{0}h)^\bot {\rm rot}\,\mathcal{P}_{0} h \big ) \|_{L^{2}(\Omega)} \\
&\,=\,
\big(
\sum_{|n| \ge 1}
\|(\mathcal{P}_{0} h)^\bot {\rm rot}\,\mathcal{P}_{n} h\|_{L^{2}(\Omega)}^{2}
\big)^{\frac12} 
+ \big( \sum_{|n| \ge 1}
\|(\mathcal{P}_{n} h)^\bot {\rm rot}\,\mathcal{P}_{0} h\|_{L^{2}(\Omega)}^{2} 
\big)^{\frac12} \\
& \le
\|\mathcal{P}_{0} h\|_{L^{\infty}(\Omega)} 
\|\nabla \mathcal{Q}_{0} h\|_{L^{2}(\Omega)}
+ \|\nabla \mathcal{P}_{0} h\|_{L^{\infty}(\Omega)} 
\|\mathcal{Q}_{0} h\|_{L^{2}(\Omega)}\,.
\end{split}
\end{equation}
For the last term in the right-hand side of  \eqref{proof.nonlinear.nonlinearity.1} with $v$ replaced by $h$, we observe that 
\begin{equation*}
\mathcal{P}_{n} 
\big(
(\mathcal{Q}_{0} h)^\bot {\rm rot}\, \mathcal{Q}_{0} h 
\big)
\,=\,
\sum_{k\in \Z\setminus\{ 0,n\}} 
(\mathcal{P}_{k} h)^\bot {\rm rot}\,\mathcal{P}_{n-k} h\,.
\end{equation*}
Then, applying the H\"{o}lder inequality we have
\begin{equation}\label{proof.nonlinear.2}
\begin{split}
\|\mathcal{P}_{0}
\big( (\mathcal{Q}_{0} h)^\bot {\rm rot}\,\mathcal{Q}_{0} h\big)
\|_{L^{1}(\Omega)}
& \le
\sum_{k \in \Z\setminus \{ 0\}} 
\|\mathcal{P}_{k} h\|_{L^{2}(\Omega)}
\|\mathcal{P}_{-k} \nabla h\|_{L^{2}(\Omega)} \\
& \le
\|\mathcal{Q}_{0} h\|_{L^{2}(\Omega)} 
\|\nabla \mathcal{Q}_{0} h\|_{L^{2}(\Omega)}\,,
\end{split}
\end{equation}
and the Young inequality for sequences implies that
\begin{equation}\label{proof.nonlinear.3}
\begin{split}
\|\mathcal{Q}_{0}
\big( (\mathcal{Q}_{0} h)^\bot {\rm rot}\,\mathcal{Q}_{0} h \big)
\|_{L^{2}(\Omega)}
& \,=\,
\big(
\sum_{|n| \ge 1} 
\|\mathcal{P}_{n} (\mathcal{Q}_{0} h)^\bot {\rm rot}\,\mathcal{Q}_{0} h \|_{L^{2}(\Omega)}^{2} 
\big)^{\frac12} \\
& \le
\Big(
\sum_{|n| \ge 1} 
\big(
\sum_{k\in Z \setminus\{0,n\}} 
\|\mathcal{P}_{k} h\|_{L^{\infty}(\Omega)}
\|\nabla \mathcal{P}_{n-k} h\|_{L^{2}(\Omega)}
\big)^{2}
\Big)^{\frac12} \\
& \le
\big(
\sum_{|n| \ge 1} 
\|\mathcal{P}_{n} h\|_{L^{\infty}(\Omega)}
\big)\,
\|\nabla \mathcal{Q}_{0} h\|_{L^{2}(\Omega)}\,.
\end{split}
\end{equation}
Next we estimate the difference $G(h^{(1)})-G(h^{(2)})$ for $h^{(1)}, h^{(2)} \in X$. 
Setting ${\bf h}=(h^{(1)}, h^{(2)})$ for simplicity, we define the function $H({\bf h})$ on $\mathcal{B}_{\vec{\delta},\epsilon} \times \mathcal{B}_{\vec{\delta},\epsilon}$ as
\begin{equation}\label{proof.nonlinear.nonlinearity.2}
\begin{split}
&~~~H({\bf h}) 
\,:=\, G(h^{(1)}) - G(h^{(2)}) \\
&\qquad
\,=\,
\mathcal{Q}_{0}\Big ((\mathcal{P}_{0} h^{(1)} - \mathcal{P}_{0} h^{(2)})^\bot {\rm rot}\,\mathcal{Q}_{0} h^{(1)} 
+ (\mathcal{P}_{0} h^{(2)})^\bot {\rm rot}\, (\mathcal{Q}_{0} h^{(1)} - \mathcal{Q}_{0} h^{(2)})  \\
&\qquad\quad
\quad + (\mathcal{Q}_{0}h^{(1)} - \mathcal{Q}_{0} h^{(2)})^\bot {\rm rot}\,\mathcal{P}_{0} h^{(1)}  
+ (\mathcal{Q}_{0}h^{(2)})^\bot {\rm rot}\,(\mathcal{P}_{0} h^{(1)} - \mathcal{P}_{0} h^{(2)}) \Big) \\
&\qquad\quad
+ ( \mathcal{Q}_{0}h^{(1)} - \mathcal{Q}_{0}h^{(2)} )^\bot {\rm rot}\,\mathcal{Q}_{0} h^{(1)}  
+ ( \mathcal{Q}_{0} h^{(2)})^\bot {\rm rot}\, (\mathcal{Q}_{0} h^{(1)} - \mathcal{Q}_{0} h^{(2)})\,.
\end{split}
\end{equation}
The following estimates on $H({\bf h})$ are obtained in the same way as \eqref{proof.nonlinear.1}-\eqref{proof.nonlinear.3}:%
\begin{align}
\|\mathcal{P}_{0} H({\bf h}) \|_{L^{1}(\Omega)} 
& \le
\big(
\|\nabla \mathcal{Q}_{0} h^{(1)}\|_{L^{2}(\Omega)} 
+ \|\mathcal{Q}_{0} h^{(2)}\|_{L^{2}(\Omega)} 
\big)\,
\|h^{(1)} - h^{(2)} \|_{X_{0}}\,,
\label{proof.nonlinear.4} \\
\|\mathcal{Q}_{0} H({\bf h}) \|_{L^{2}(\Omega)} 
& \le
\Big(
\|\nabla \mathcal{Q}_{0} h^{(1)}\|_{L^{2}(\Omega)} 
+ \|\mathcal{Q}_{0} h^{(2)}\|_{L^{2}(\Omega)} 
+ \|\nabla \mathcal{P}_{0} h^{(1)}\|_{L^{\infty}(\Omega)} \nonumber \\
&  
+ \|\mathcal{P}_{0} h^{(2)}\|_{L^{\infty}(\Omega)} 
+ \sum_{|n| \ge 1} \| \mathcal{P}_{n} h^{(2)}\|_{L^{\infty}(\Omega)}
\Big)\,
\|h^{(1)} - h^{(2)} \|_{X_{0}}\,.
\label{proof.nonlinear.5}
\end{align}
From \eqref{proof.nonlinear.1}-\eqref{proof.nonlinear.3},  \eqref{proof.nonlinear.4}-\eqref{proof.nonlinear.5}, and the definition of $\mathcal{B}_{\vec{\delta},\epsilon}$, we obtain the following estimates on~$G(h)$ and $H({\bf h})$ in $Y$.
\begin{align}
\|\mathcal{P}_{0} G(h) \|_{L^{1}(\Omega)}
& \le
\epsilon^{2} |\alpha|^{\delta_{2} + \delta_{3}}\,,
\label{proof.nonlinear.G.0} \\
\|\mathcal{Q}_{0} G(h) \|_{L^{2}(\Omega)}
& \le
\epsilon^{2} 
(|\alpha|^{\delta_{1}+\delta_{2}}
+ |\alpha|^{\delta_{1}+\delta_{3}}
+ |\alpha|^{\delta_{3}+\delta_{4}})\,, 
\label{proof.nonlinear.G.1} \\
\|\mathcal{P}_{0} H({\bf h}) \|_{L^{1}(\Omega)} 
& \le
\epsilon (|\alpha|^{\delta_{2}}+|\alpha|^{\delta_{3}})\,
\|h^{(1)} - h^{(2)} \|_{X_{0}} \,,
\label{proof.nonlinear.H.0} \\
\|\mathcal{Q}_{0} H({\bf h}) \|_{L^{2}(\Omega)} 
& \le
\epsilon
(|\alpha|^{\delta_{1}} + |\alpha|^{\delta_{2}} 
+ |\alpha|^{\delta_{3}} + |\alpha|^{\delta_{4}} )\,
\|h^{(1)} - h^{(2)} \|_{X_{0}}\,.
\label{proof.nonlinear.H.1}
\end{align}

Now let us define the mapping $\Phi : \mathcal{B}_{\vec{\delta},\epsilon} \rightarrow X_0$ by setting $\Phi [h]: =  v_{h}$, where $v_{h}$ is the unique solution to  \eqref{proof.nonlinear.linear}. Our aim is to show that 

\noindent 
(i) $\Phi$ is Lipschitz continuous on $\mathcal{B}_{\vec{\delta},\epsilon}$ in the topology of $X_0$. Namely, there exists $\tau\in (0,1)$ such that $\| \Phi[h^{(1)}] - \Phi[h^{(2)}]\|_{X_0} \leq \tau \|h^{(1)} - h^{(2)}\|_{X_0}$ for any $h^{(1)}, h^{(2)}\in \mathcal{B}_{\vec{\delta},\epsilon}$, 

\noindent 
(ii) $\Phi$ is a mapping from $\mathcal{B}_{\vec{\delta},\epsilon}$ into $\mathcal{B}_{\vec{\delta},\epsilon}$, if the pair $(\vec{\delta}, \epsilon)$
and the external force $f=(\mathcal{P}_{0} f, \mathcal{Q}_{0} f)\in Y$ satisfy a suitable condition.

\medskip 

\noindent
For convenience in the following proof, let $K_{1}$ and $K_{2}$ denote the largest constant~$C$ (larger than $1$ without loss of generality) appearing in \eqref{est.prop.all.alpha.0},  \eqref{est.proof.nonlinear.1}-\eqref{est.proof.nonlinear.5}, and \eqref{est.prop.all.alpha.0}, \eqref{est.proof.nonlinear.1}-\eqref{est.proof.nonlinear.4},  \eqref{est.proof.nonlinear.6}-\eqref{est.proof.nonlinear.7}, respectively.

We first show (i). For any ${\bf h}=(h^{(1)}, h^{(2)}) \in \mathcal{B}_{\vec{\delta},\epsilon} \times \mathcal{B}_{\vec{\delta},\epsilon}$,
we observe that the differences~$u_{{\bf h}}:=\Phi [h^{(1)}]-\Phi [h^{(2)}]$ and $p_{{\bf h}}:=q_{h^{(1)}}-q_{h^{(2)}}$ solve the following system:
\begin{equation*}
  \left\{
\begin{aligned}
-\Delta u_{{\bf h}}
- \alpha ( x^\bot \cdot \nabla u_{{\bf h}} - u_{{\bf h}}^\bot ) 
+ \nabla p_{{\bf h}}
+ \alpha U^\bot {\rm rot}\,u_{{\bf h}}
& \,=\, 
H({\bf h})\,,  \qquad \qquad x\in \Omega \,,\\
{\rm div}\,u_{{\bf h}} &\,=\,0\,,   \qquad \qquad  x \in \Omega\,, \\
u_{{\bf h}} &\,=\, 0 \,, \qquad \qquad  x \in \partial \Omega\,,
\end{aligned}\right.
\end{equation*}
where $H({\bf h})$ is defined in \eqref{proof.nonlinear.nonlinearity.2}. We consider the case $0<|\alpha|<1$. Then for any~$h^{(1)}, h^{(2)} $ in~$ X_{0}$, by  \eqref{est.prop.all.alpha.0} and \eqref{est.proof.nonlinear.1}-\eqref{est.proof.nonlinear.5} combined with \eqref{proof.nonlinear.H.0}-\eqref{proof.nonlinear.H.1}, we see that
\begin{equation*}
\begin{split}
&~~~\|\Phi [h^{(1)}] - \Phi [h^{(2)}]\|_{X_{0}} \\
& \,=\,
\|\mathcal{P}_{0} u_{{\bf h}}\|_{L^{\infty}(\Omega)}
+ \|\nabla \mathcal{P}_{0} u_{{\bf h}}\|_{L^{\infty}(\Omega)} 
+ \|\mathcal{Q}_{0} u_{{\bf h}} \|_{L^{2}(\Omega)}
+ \|\nabla \mathcal{Q}_{0} u_{{\bf h}} \|_{L^{2}(\Omega)}
+ \sum_{|n| \ge 1} 
\|\mathcal{P}_{n} u_{{\bf h}}\|_{L^{\infty}(\Omega)} \\
& \le
K_{1} \|\mathcal{P}_{0} H({\bf h})\|_{L^{1}(\Omega)}
+ \frac{7 K_{1}}{|\alpha|} \|\mathcal{Q}_{0} H({\bf h})\|_{L^{2}(\Omega)} \\
& \le
8K_{1} \epsilon\,
(|\alpha|^{\delta_{1}-1} + |\alpha|^{\delta_{2}-1} 
+ |\alpha|^{\delta_{3}-1} + |\alpha|^{\delta_{4}-1} )\,
\|h^{(1)} - h^{(2)} \|_{X_{0}} \,.
\end{split}
\end{equation*}
Hence, if we choose the pair $(\vec{\delta}, \epsilon)$ to satisfy
\begin{equation}\label{proof.nonlinear.cond.small.alpha.1}
\delta_{j} \ge 1 \quad {\rm for} \quad j\,=\,1\ldots4\,,
\qquad 0 < \epsilon < \frac{1}{32K_{1}}\,,
\qquad {\rm when} \quad 0 < |\alpha| < 1\,,
\end{equation}
then the mapping $\Phi$ is Lipschitz continuous on the set $\mathcal{B}_{\vec{\delta},\epsilon}$. For the case $|\alpha| \ge 1$, from the estimates \eqref{est.prop.all.alpha.0},  \eqref{est.proof.nonlinear.1}-\eqref{est.proof.nonlinear.4}, 
and \eqref{est.proof.nonlinear.6}-\eqref{est.proof.nonlinear.7} combined with \eqref{proof.nonlinear.H.0}-\eqref{proof.nonlinear.H.1}, we have
\begin{align}
\|\Phi [h^{(1)}] - \Phi [h^{(2)}]\|_{X_{0}} 
& \le
K_{2} \|\mathcal{P}_{0} H({\bf h})\|_{L^{1}(\Omega)}
+ 8 K_{2}
\|\mathcal{Q}_{0} H({\bf h})\|_{L^{2}(\Omega)} \nonumber \\
& \le
9K_{1} \epsilon\,
(|\alpha|^{\delta_{1}} + |\alpha|^{\delta_{2}}
+ |\alpha|^{\delta_{3}} + |\alpha|^{\delta_{4}} )\,
\|h^{(1)} - h^{(2)} \|_{X_{0}}\,.
\label{proof.nonlinear.middle.alpha.difference}
\end{align}
Then  we obtain the next condition 
\begin{equation}\label{proof.nonlinear.cond.middle.alpha.1}
\delta_{j} \le 0 \quad {\rm for} \quad j\,=\,1\ldots4\,,
\qquad 0 < \epsilon < \frac{1}{36K_{2}}\,,
\qquad {\rm when} \quad  |\alpha| \ge 1\,,
\end{equation}
for the Lipschitz continuity of $\Phi$ on $\mathcal{B}_{\vec{\delta},\epsilon}$. We have   shown (i) provided the pair $(\vec{\delta}, \epsilon)$ satisfies the conditions \eqref{proof.nonlinear.cond.small.alpha.1} or \eqref{proof.nonlinear.cond.middle.alpha.1}. 

Next we prove  (ii). In the case $0<|\alpha|<1$, the estimates \eqref{est.prop.all.alpha.0} and \eqref{proof.nonlinear.G.0} imply
\begin{align}
\|\mathcal{P}_{0} \Phi [h]\|_{L^{\infty}(\Omega)}
+ \|\nabla \mathcal{P}_{0} \Phi [h]\|_{L^{\infty}(\Omega)} \nonumber 
& \le K_{1}
\big(
\|(\mathcal{P}_{0} G(h))_{\theta}\|_{L^{1}(\Omega)}
+ \|(\mathcal{P}_{0} f)_{\theta}\|_{L^{1}(\Omega)}
\big) \nonumber \\
& \le K_{1}
\big(
\epsilon^{2} |\alpha|^{\delta_{2}+\delta_{3}}
+ \|(\mathcal{P}_{0} f)_{\theta}\|_{L^{1}(\Omega)}
\big) \,,
\label{proof.nonlinear.6}
\end{align}
for any $h \in X_{0}$. From \eqref{est.proof.nonlinear.1}-\eqref{est.proof.nonlinear.2}, 
\eqref{est.proof.nonlinear.3}-\eqref{est.proof.nonlinear.4}, 
and \eqref{proof.nonlinear.G.1} we have,
\begin{align*}
& \|\mathcal{Q}_{0} \Phi [h] \|_{L^{2}(\Omega)} 
\le
\frac{3K_{1}}{|\alpha|} 
\big(
\|\mathcal{Q}_{0} G(h)\|_{L^{2}(\Omega)} + \|\mathcal{Q}_{0} f\|_{L^{2}(\Omega)}
\big) \\
& \qquad
\le
3K_{1} \Big(
\epsilon^{2} 
(|\alpha|^{\delta_{1}+\delta_{2}-1}
+ |\alpha|^{\delta_{1}+\delta_{3}-1}
+ |\alpha|^{\delta_{3}+\delta_{4}-1} )
+ |\alpha|^{-1} \|\mathcal{Q}_{0} f\|_{L^{2}(\Omega)}
\Big)\,, \\
& \|\nabla \mathcal{Q}_{0} \Phi [h] \|_{L^{2}(\Omega)}
\le
\frac{3K_{1}}{|\alpha|^{\frac12}} 
\big(
\|\mathcal{Q}_{0} G(h)\|_{L^{2}(\Omega)} + \|\mathcal{Q}_{0} f\|_{L^{2}(\Omega)}
\big) \\
& \qquad
\le
3K_{1} \Big(
\epsilon^{2} 
\big(|\alpha|^{\delta_{1}+\delta_{2}-\frac12}
+ |\alpha|^{\delta_{1}+\delta_{3}-\frac12}
+ |\alpha|^{\delta_{3}+\delta_{4}-\frac12} \big)
+ |\alpha|^{-\frac12} \|\mathcal{Q}_{0} f\|_{L^{2}(\Omega)}
\Big)\,, 
\end{align*}
and
\begin{align*}
\sum_{|n| \ge 1} 
& \|\mathcal{P}_{n} \Phi [h]\|_{L^{\infty}(\Omega)} 
\le
\frac{K_{1}}{|\alpha|^{\frac34}} 
\big(
\|\mathcal{Q}_{0} G(h)\|_{L^{2}(\Omega)} + \|\mathcal{Q}_{0} f\|_{L^{2}(\Omega)}
\big) \nonumber \\
& \le
K_{1} \Big(
\epsilon^{2} 
\big(|\alpha|^{\delta_{1}+\delta_{2}-\frac34}
+ |\alpha|^{\delta_{1}+\delta_{3}-\frac34}
+ |\alpha|^{\delta_{3}+\delta_{4}-\frac34} \big)
+ |\alpha|^{-\frac34} \|\mathcal{Q}_{0} f\|_{L^{2}(\Omega)}
\Big)\,.
\end{align*}
Hence, recalling the condition \eqref{proof.nonlinear.cond.small.alpha.1}, if we choose the pair $({\vec \delta}, \epsilon)$ and $f \in Y$ to satisfy
\begin{align}
\delta_{1} \,=\, \delta_{2} \,=\, \delta_{3} \,=\, \delta_{4} \,=\,1\,,
\qquad 0 < \epsilon < \frac{1}{32K_{1}}\,,
\qquad {\rm when} \quad 0 < |\alpha| < 1\,,
\label{proof.nonlinear.choice.small.alpha.1}
\end{align}
and
\begin{align}
\|(\mathcal{P}_{0} f)_{\theta}\|_{L^{1}(\Omega)} 
\le \frac{\epsilon}{2K_{1}} |\alpha|\,,
\quad
\|\mathcal{Q}_{0} f\|_{L^{2}(\Omega)} 
\le \frac{\epsilon}{6K_{1}} |\alpha|^{2}\,,
\quad {\rm when} \quad 0 < |\alpha| < 1\,,
\label{proof.nonlinear.choice.small.alpha.2}
\end{align}
then $\Phi$ defines a mapping from $\mathcal{B}_{\vec{\delta},\epsilon}$ into itself. 

\noindent
In the case $|\alpha| \ge1$, we have \eqref{proof.nonlinear.6} with $K_{1}$ replaced by $K_{2}$, and moreover, from \eqref{est.proof.nonlinear.1}-\eqref{est.proof.nonlinear.2}, 
\eqref{est.proof.nonlinear.3}-\eqref{est.proof.nonlinear.4}, 
and \eqref{est.proof.nonlinear.6}-\eqref{est.proof.nonlinear.7}, along with~(\ref{proof.nonlinear.G.1}),
we have
\begin{align}
& \|\mathcal{Q}_{0} \Phi [h]\|_{L^{2}(\Omega)}
\le
\frac{3K_{2}}{|\alpha|^{\frac12}} 
\big(
\|\mathcal{Q}_{0} G(h)\|_{L^{2}(\Omega)}
+ \|\mathcal{Q}_{0} f\|_{L^{2}(\Omega)}
\big) \nonumber \\
& \qquad
\le
3K_{2} \Big(
\epsilon^{2} 
(|\alpha|^{\delta_{1}+\delta_{2}-\frac12}
+ |\alpha|^{\delta_{1}+\delta_{3}-\frac12}
+ |\alpha|^{\delta_{3}+\delta_{4}-\frac12} )
+ |\alpha|^{-\frac12} \|\mathcal{Q}_{0} f\|_{L^{2}(\Omega)}
\Big)\,, \label{proof.nonlinear.7} \\
& \|\nabla \mathcal{Q}_{0} \Phi [h]\|_{L^{2}(\Omega)}
\le
3 K_{2}
\big(
\|\mathcal{Q}_{0} G(h)\|_{L^{2}(\Omega)}
+ \|\mathcal{Q}_{0} f\|_{L^{2}(\Omega)}
\big) \nonumber \\
& \qquad
\le
3 K_{2} \Big(
\epsilon^{2} 
(|\alpha|^{\delta_{1}+\delta_{2}}
+ |\alpha|^{\delta_{1}+\delta_{3}}
+ |\alpha|^{\delta_{3}+\delta_{4}} )
+ \|\mathcal{Q}_{0} f\|_{L^{2}(\Omega)}
\Big)\,, \label{proof.nonlinear.8} 
\end{align}
and
\begin{align}
& \sum_{|n| \ge 1} 
\|\mathcal{P}_{n} \Phi [h]\|_{L^{\infty}(\Omega)} 
\le
\frac{2K_{2}}{|\alpha|^{\frac14}} 
\big(
\|\mathcal{Q}_{0} G(h)\|_{L^{2}(\Omega)}
+ \|\mathcal{Q}_{0} f\|_{L^{2}(\Omega)}
\big) \nonumber \\
& \qquad
\le
2K_{2} \Big(
\epsilon^{2} 
(|\alpha|^{\delta_{1}+\delta_{2}-\frac14}
+ |\alpha|^{\delta_{1}+\delta_{3}-\frac14}
+ |\alpha|^{\delta_{3}+\delta_{4}-\frac14} )
+ |\alpha|^{-\frac14} \|\mathcal{Q}_{0} f\|_{L^{2}(\Omega)}
\Big)\,.\label{proof.nonlinear.9}
\end{align}
Then recalling \eqref{proof.nonlinear.cond.middle.alpha.1}, if the pair $({\vec \delta}, \epsilon)$ and $f \in Y$ satisfy
\begin{equation}\label{proof.nonlinear.choice.middle.alpha.1}
\delta_{1} \,=\, \delta_{2} \,=\, \delta_{3} \,=\, \delta_{4} \,=\, 0\,,
\qquad 0 < \epsilon < \frac{1}{36K_{2}}\,,
\qquad {\rm when} \quad |\alpha| \ge 1\,,
\end{equation}
and
\begin{equation}\label{proof.nonlinear.choice.middle.alpha.2}
\|(\mathcal{P}_{0} f)_{\theta}\|_{L^{1}(\Omega)} 
\le \frac{\epsilon}{2K_{2}}\,,
\qquad
\|\mathcal{Q}_{0} f\|_{L^{2}(\Omega)} 
\le \frac{\epsilon}{6K_{2}}\,,
\qquad {\rm when} \quad |\alpha| \ge 1\,,
\end{equation}
then we see that $\Phi$ defines a mapping $\Phi :\mathcal{B}_{\vec{\delta},\epsilon} \to \mathcal{B}_{\vec{\delta},\epsilon}$. 

Now we have shown that the mapping $\Phi$ defines a contraction on $\mathcal{B}_{\vec{\delta},\epsilon}$ under the conditions \eqref{proof.nonlinear.choice.small.alpha.1} and \eqref{proof.nonlinear.choice.small.alpha.2} when $0<|\alpha|<1$, and under \eqref{proof.nonlinear.choice.middle.alpha.1} and \eqref{proof.nonlinear.choice.middle.alpha.2} when $|\alpha|\ge1$. Then there is a unique fixed point of $\Phi$ in the closed convex set $\mathcal{B}_{\vec{\delta},\epsilon}$. Hence we finally obtain a unique solution to the nonlinear problem \eqref{hatNS_alpha} in $\mathcal{B}_{\vec{\delta},\epsilon}$.

The estimates in Theorem \ref{thm.main.1} are obtained as follows. We consider only the case $0<|\alpha|<1$, in particular the estimates \eqref{est.thm.main.1.1}-\eqref{est.thm.main.1.3}, since the case $|\alpha|\ge1$ can be handled similarly. Let $v$ denote the unique fixed point of $\Phi$ in $\mathcal{B}_{\vec{\delta},\epsilon}$. Applying the linear estimates~\eqref{est.proof.nonlinear.1} and \eqref{est.proof.nonlinear.2} to \eqref{hatNS_alpha} along with~(\ref{proof.nonlinear.1})-(\ref{proof.nonlinear.3}) we see that
\begin{align*}
\|\mathcal{Q}_{0} v\|_{L^{2}(\Omega)}
& \le
\frac{3K_{1}}{|\alpha|} 
\big(
\|\mathcal{Q}_{0} G(v)\|_{L^{2}(\Omega)} 
+ \|\mathcal{Q}_{0} f\|_{L^{2}(\Omega)}
\big) \\
& \le
\frac{3K_{1}}{|\alpha|} 
\Big(
\|\mathcal{P}_{0} v\|_{L^{\infty}(\Omega)} 
\|\nabla \mathcal{Q}_{0} v\|_{L^{2}(\Omega)} 
+ \|\nabla \mathcal{P}_{0} v\|_{L^{\infty}(\Omega)} 
\|\mathcal{Q}_{0} v\|_{L^{2}(\Omega)} \\
& \qquad\qquad
+ \big(\sum_{|n| \ge 1} \|\mathcal{P}_{n} v\|_{L^{\infty}(\Omega)} \big)
\|\mathcal{Q}_{0} v\|_{L^{2}(\Omega)}
+ \|\mathcal{Q}_{0} f\|_{L^{2}(\Omega)}
\Big) \\
& \le
3K_{1} \epsilon 
\big( 
\|\nabla \mathcal{Q}_{0} v\|_{L^{2}(\Omega)} 
+ 2 \|\mathcal{Q}_{0} v\|_{L^{2}(\Omega)} 
\big)
+ \frac{3K_{1}}{|\alpha|} 
\|\mathcal{Q}_{0} f\|_{L^{2}(\Omega)}\,.
\end{align*}
In the last inequality we have applied the bounds derived from the assumption $v \in \mathcal{B}_{\vec{\delta},\epsilon}$ with~$\delta_j = 1$ for all~$j$. Since $0<6K_{1} \epsilon<\frac15$ under the choice of $\epsilon$ in \eqref{proof.nonlinear.choice.small.alpha.1}, we obtain
\begin{align}
\|\mathcal{Q}_{0} v\|_{L^{2}(\Omega)}
& \le
\frac{15}{4}
K_{1} \epsilon
\|\nabla \mathcal{Q}_{0} v\|_{L^{2}(\Omega)} 
+ 
\frac{15}{4}
\frac{K_{1}}{|\alpha|} 
\|\mathcal{Q}_{0} f\|_{L^{2}(\Omega)}\,.
\label{proof.nonlinear.10}
\end{align}
On the other hand, if we apply \eqref{est.proof.nonlinear.3} and \eqref{est.proof.nonlinear.4} to \eqref{hatNS_alpha}, then we have by the same argument
\begin{align*}
\|\nabla \mathcal{Q}_{0} v\|_{L^{2}(\Omega)}
& \le
\frac{3K_{1}}{|\alpha|^{\frac12}} 
\big(
\|\mathcal{Q}_{0} G(v)\|_{L^{2}(\Omega)} 
+ \|\mathcal{Q}_{0} f\|_{L^{2}(\Omega)}
\big) \\
& \le
6 K_{1} \epsilon |\alpha|^{\frac12} \|\mathcal{Q}_{0} v\|_{L^{2}(\Omega)}
+ 3K_{1} \epsilon |\alpha|^{\frac12} \|\nabla \mathcal{Q}_{0} v\|_{L^{2}(\Omega)} 
+ \frac{3K_{1}}{|\alpha|^{\frac12}} 
\|\mathcal{Q}_{0} f\|_{L^{2}(\Omega)}\,.
\end{align*}
Inserting \eqref{proof.nonlinear.10} in the above inequality, and by the smallness of $K_{1} \epsilon$ again, we see that
\begin{align*}
\|\nabla \mathcal{Q}_{0} v\|_{L^{2}(\Omega)}
& \le
\frac{C}{|\alpha|^{\frac12}} 
\|\mathcal{Q}_{0} f\|_{L^{2}(\Omega)}\,,
\end{align*}
which implies the estimate \eqref{est.thm.main.1.3}. We can obtain \eqref{est.thm.main.1.2} by \eqref{proof.nonlinear.10} combined with 
\eqref{est.thm.main.1.3} with the condition $0<|\alpha|<1$. 
  For the estimate \eqref{est.thm.main.1.2'}, by \eqref{est.proof.nonlinear.5}, (\ref{proof.nonlinear.1}) and (\ref{proof.nonlinear.3}), and the condition $v \in \mathcal{B}_{\vec{\delta},\epsilon}$ with~$\delta_j = 1$ for all~$j$ we have
\begin{align*}
\sum_{|n| \ge 1} \|\mathcal{P}_{n} v\|_{L^{\infty}(\Omega)}
& \le
\frac{K_{1}}{|\alpha|^{\frac34}} 
\big(
\|\mathcal{Q}_{0} G(v)\|_{L^{2}(\Omega)} 
+ \|\mathcal{Q}_{0} f\|_{L^{2}(\Omega)}
\big) \\
& \le
K_{1} \epsilon |\alpha|^{\frac14}
\big( 
\|\nabla \mathcal{Q}_{0} v\|_{L^{2}(\Omega)} 
+ 2 \|\mathcal{Q}_{0} v\|_{L^{2}(\Omega)} 
\big)
+ \frac{K_{1}}{|\alpha|^{\frac34}} 
\|\mathcal{Q}_{0} f\|_{L^{2}(\Omega)}\,,
\end{align*}
which combined with \eqref{est.thm.main.1.2} and \eqref{est.thm.main.1.3} leads to \eqref{est.thm.main.1.2'}. 
The estimate \eqref{est.thm.main.1.1} follows from \eqref{est.prop.all.alpha.0} 
and the nonlinear estimate \eqref{proof.nonlinear.2} with $h$ replaced by $v$. This completes the proof of Theorem \ref{thm.main.1}.
\BOX
\subsection{Proof for large $|\alpha|  $}\label{subsec.nonlinear.large.alpha}
In this subsection we prove Theorem \ref{thm.main.2}. 
We adopt the same notation as in the proof of Theorem \ref{thm.main.1} of the previous subsection. Let $K_{3}$ denote the largest constant of $C$ (larger than~$1$ without loss of generality) appearing in \eqref{est.prop.all.alpha.0} and \eqref{est.proof.nonlinear.1''}-\eqref{est.proof.nonlinear.7'} for convenience. 

We show that the mapping $\Phi$ is a contraction on the set $\mathcal{B}_{\vec{\delta},\epsilon}$. For the Lipschitz continuity of $\Phi$ on $\mathcal{B}_{\vec{\delta},\epsilon}$, we observe that the estimate \eqref{proof.nonlinear.middle.alpha.difference} is improved in terms of the decay in $|\alpha|$ into
\begin{align*}
& \|\Phi [h^{(1)}] - \Phi [h^{(2)}]\|_{X_{0}} 
\le
K_{3} \|(\mathcal{P}_{0} H({\bf h}))_{\theta}\|_{L^{1}(\Omega)}
+ \frac{6 K_{3}}{|\alpha|^{\frac13}} 
\|\mathcal{Q}_{0} H({\bf h})\|_{L^{2}(\Omega)} \nonumber \\
& \qquad
\le
7K_{3} \epsilon\,
( |\alpha|^{\delta_{1}-\frac13}+ |\alpha|^{\delta_{2}} 
+ |\alpha|^{\delta_{3}} + |\alpha|^{\delta_{4}-\frac13} )\,
\|h^{(1)} - h^{(2)} \|_{X_{0}}\,.
\end{align*}
Hence, if we choose the pair $(\vec{\delta}, \epsilon)$ to satisfy
\begin{equation}\label{proof.nonlinear.cond.large.alpha.1}
\delta_{1} \le \frac13\,, \quad
\delta_{2} \le 0\,, \quad 
\delta_{3} \le 0\,, \quad 
\delta_{4} \le \frac13\,, 
\quad 0 < \epsilon < \frac{1}{28K_{3}}\,,
\end{equation}
then the mapping $\Phi$ is Lipschitz continuous on $\mathcal{B}_{\vec{\delta},\epsilon}$. 

Next we check that $\Phi$ is a mapping from $\mathcal{B}_{\vec{\delta},\epsilon}$ into $\mathcal{B}_{\vec{\delta},\epsilon}$. We have \eqref{proof.nonlinear.6} with $K_{1}$ replaced by $K_{3}$, and we see that \eqref{proof.nonlinear.7}-\eqref{proof.nonlinear.9} are improved in terms of decay in $|\alpha|$ respectively to
\begin{align*}
& \|\mathcal{Q}_{0} \Phi [h]\|_{L^{2}(\Omega)}
\le
\frac{2K_{3}}{|\alpha|^{\frac23}} 
\big(
\|\mathcal{Q}_{0} G(h)\|_{L^{2}(\Omega)}
+ \|\mathcal{Q}_{0} f\|_{L^{2}(\Omega)}
\big) \\
& \qquad
\le
2K_{3}\Big(
\epsilon^{2} 
\big(|\alpha|^{\delta_{1}+\delta_{2}-\frac23}
+ |\alpha|^{\delta_{1}+\delta_{3}-\frac23}
+ |\alpha|^{\delta_{3}+\delta_{4}-\frac23} \big)
+ |\alpha|^{-\frac23} \|\mathcal{Q}_{0} f\|_{L^{2}(\Omega)}
\Big)\,, \\
& \|\nabla \mathcal{Q}_{0} \Phi [h]\|_{L^{2}(\Omega)}
\le
\frac{2K_{3}}{|\alpha|^{\frac13}} 
\big(
\|\mathcal{Q}_{0} G(h)\|_{L^{2}(\Omega)}
+ \|\mathcal{Q}_{0} f\|_{L^{2}(\Omega)}
\big) \\
& \qquad
\le
2K_{3} \Big(
\epsilon^{2} 
\big(|\alpha|^{\delta_{1}+\delta_{2}-\frac13}
+ |\alpha|^{\delta_{1}+\delta_{3}-\frac13}
+ |\alpha|^{\delta_{3}+\delta_{4}-\frac13}\big )
+ |\alpha|^{-\frac13} \|\mathcal{Q}_{0} f\|_{L^{2}(\Omega)}
\Big)\,,  
\end{align*}
and
\begin{align*}
& \sum_{|n| \ge 1} 
\|\mathcal{P}_{n} \Phi [h]\|_{L^{\infty}(\Omega)} 
\le
\frac{2K_{3} (\log |\alpha|)^{\frac12}}{|\alpha|^{\frac12}} 
\big(
\|\mathcal{Q}_{0} G(h)\|_{L^{2}(\Omega)}
+ \|\mathcal{Q}_{0} f\|_{L^{2}(\Omega)}
\big) \nonumber \\
& \qquad
\le
2K_{3}\Big(\epsilon^{2} 
(\log |\alpha|)^{\frac12}
\big(|\alpha|^{\delta_{1}+\delta_{2}-\frac12}
+ |\alpha|^{\delta_{1}+\delta_{3}-\frac12}
+ |\alpha|^{\delta_{3}+\delta_{4}-\frac12}\big ) \nonumber \\
& \qquad\qquad
+ (\log |\alpha|)^{\frac12} |\alpha|^{-\frac12} 
\|\mathcal{Q}_{0} f\|_{L^{2}(\Omega)}
\Big)\,.
\end{align*}
Then, under the condition \eqref{proof.nonlinear.cond.large.alpha.1}, 
if we choose the pair $({\vec \delta}, \epsilon)$ and $f \in Y$ to satisfy
\begin{equation}\label{proof.nonlinear.choice.large.alpha.1}
\delta_{1} \,=\, \delta_{4} \,=\, \frac13\,,\quad
\delta_{2} \,=\, \delta_{3} \,=\, 0\,,
\qquad 0 < \epsilon < \frac{1}{28K_{3}}\,,
\end{equation}
and
\begin{equation}\label{proof.nonlinear.choice.large.alpha.2}
\|(\mathcal{P}_{0} f)_{\theta}\|_{L^{1}(\Omega)} 
\le \frac{\epsilon}{2 K_{3}} |\alpha|^{\frac13}\,,
\quad
\|\mathcal{Q}_{0} f\|_{L^{2}(\Omega)} 
\le \frac{\epsilon}{6 K_{2}} |\alpha|^{\frac13}\,,
\end{equation}
then $\Phi$ defines a mapping $\Phi :\mathcal{B}_{\vec{\delta},\epsilon} \to \mathcal{B}_{\vec{\delta},\epsilon}$. Now we have shown that $\Phi$ is a contraction on~$\mathcal{B}_{\vec{\delta},\epsilon}$ if we assume the conditions
\eqref{proof.nonlinear.choice.large.alpha.1} and \eqref{proof.nonlinear.choice.large.alpha.2}. Hence there exists a unique fixed point $v$ of $\Phi$ in $\mathcal{B}_{\vec{\delta},\epsilon}$. The estimates \eqref{est.thm.main.2.1}-\eqref{est.thm.main.2.3} can be obtained in the same way as in the proof of Theorem \ref{thm.main.1}. This completes the proof of Theorem \ref{thm.main.2}.
\BOX

\subsection{Proof of Theorem \ref{thm.main.3}}\label{subsec.nonlinear.thm.1.3}
This subsection is devoted to the proof of Theorem \ref{thm.main.3}. For a given $f \in Y$, let us take $\alpha \in \R$ large enough to satisfy both the condition \eqref{assumption.thm.main.2} in Theorem \ref{thm.main.2} and the assumption in Proposition \ref{decompositionvelocity}. Then, from Theorem \ref{thm.main.2} we see that there exists a solution $(v^{(\alpha)}, q^{(\alpha)}) $ in~$ X \times W^{1,1}_{loc}(\overline{\Omega})$ to \eqref{tildeNS_alpha} satisfying the estimates  \eqref{est.thm.main.2.1}-\eqref{est.thm.main.2.3}. Hence, the proof will be completed  as soon as we show all the estimates in Theorem \ref{thm.main.3} and the decomposition~\eqref{decom.thm.main.3} of $v^{(\alpha)}$. Note that $v^{(\alpha)}$ also solves the system \eqref{hatNS_alpha}, which is introduced in the beginning of this section, with a suitable new pressure $\tilde q^{(\alpha)} \in W^{1,1}_{loc}(\overline{\Omega})$. 

\noindent
We deduce the estimate \eqref{est.thm.main.3.1} from the   triangle inequality
\begin{align}
\| v^{(\alpha)} -  v_0^{{\rm linear}}\|_{L^\infty (\Omega)} 
& \le 
\| \mathcal{Q}_{0} v^{(\alpha)} \|_{L^\infty (\Omega)}
+ \| \mathcal{P}_{0} v^{(\alpha)} -  v_0^{{\rm linear}}\|_{L^\infty (\Omega)} 
\nonumber \\
& \le 
\sum_{|n|\geq 1} \| \mathcal{P}_n v^{(\alpha)} \|_{L^\infty (\Omega)} 
+ \| \mathcal{P}_{0} v^{(\alpha)} -  v_0^{{\rm linear}}\|_{L^\infty (\Omega)}
\,. \label{proof.nonlinear.11}
\end{align}
Since $v:=\mathcal{P}_{0} v^{(\alpha)} - v_0^{{\rm linear}} \in W_0^{1,\infty} (\Omega)^2$ is a solution to the next system
\begin{equation*}
  \left\{
\begin{aligned}
  -\Delta v - \alpha ( x^\bot \cdot \nabla v - v^\bot ) + \nabla q + \alpha U^\bot {\rm rot}\, v & \,=\, 
- \mathcal{P}_{0} \Big( (\mathcal{Q}_{0} v)^\bot {\rm rot}\, \mathcal{Q}_{0} v \Big)\,,   \qquad x\in \Omega \,,\\
{\rm div}\, v &  \,=\, 0\,,   \qquad \qquad  x \in \Omega\,, \\
  v & \,=\, 0 \,, \qquad \qquad  x \in \partial \Omega\,, 
\end{aligned}\right.
\end{equation*}
with some pressure $q \in W^{1,1}_{loc}(\overline{\Omega})$, we have from \eqref{est.prop.all.alpha.0} and \eqref{proof.nonlinear.2} with $h$ replaced by $v^{(\alpha)}$, 
\begin{align}
\| \mathcal{P}_{0} v^{(\alpha)} -  v_0^{{\rm linear}}\|_{L^\infty (\Omega)}
& \le 
C\| \mathcal{Q}_{0} v^{(\alpha)}\|_{L^2 (\Omega)}
\|\nabla \mathcal{Q}_{0} v^{(\alpha)}\|_{L^2 (\Omega)}
\,. \label{proof.nonlinear.12}
\end{align}
Then, by \eqref{proof.nonlinear.11}-\eqref{proof.nonlinear.12} along with \eqref{est.thm.main.2.2}-\eqref{est.thm.main.2.3} we obtain the estimate \eqref{est.thm.main.3.1}. 

\noindent
The decomposition \eqref{decom.thm.main.3} and the related estimates follow from the results in Proposition~\ref{decompositionvelocity}. Indeed, we know that $(v^{(\alpha)}, \tilde q^{(\alpha)})$ solves \eqref{hatNS_alpha}, and we have the next estimate of~$\|\mathcal{P}_{n} G(v^{(\alpha)})\|_{L^2 (\Omega)}$, which is uniform both in $|n| \ge 1$ and $|\alpha|$, combined with estimates~\eqref{est.thm.main.2.1}-\eqref{est.thm.main.2.3}:
\begin{align*}
&\|\mathcal{P}_{n} G(v^{(\alpha)})\|_{L^2 (\Omega)} 
\le \|\mathcal{Q}_{0} G(v^{(\alpha)})\|_{L^2 (\Omega)} \\
& \le
\big\|\mathcal{Q}_{0}\big((\mathcal{P}_{0} v^{(\alpha)})^\bot {\rm rot}\,\mathcal{Q}_{0} v^{(\alpha)} \big) \big\|_{L^{2}(\Omega)}
+ \big\|\mathcal{Q}_{0}\big ((\mathcal{Q}_{0} v^{(\alpha)})^\bot {\rm rot}\,\mathcal{P}_{0} v^{(\alpha)} \big ) \big\|_{L^{2}(\Omega)} \\
& \quad
+ \big\| \mathcal{Q}_{0}
\big(
(\mathcal{Q}_{0} v^{(\alpha)})^\bot {\rm rot}\, \mathcal{Q}_{0} v^{(\alpha)} 
\big) \big\|_{L^2 (\Omega)} \\
& \le
\|\mathcal{P}_{0} v^{(\alpha)}\|_{L^{\infty}(\Omega)} 
\|\nabla \mathcal{Q}_{0} v^{(\alpha)}\|_{L^{2}(\Omega)}
+ \|\nabla \mathcal{P}_{0} v^{(\alpha)}\|_{L^{\infty}(\Omega)} 
\|\mathcal{Q}_{0} v^{(\alpha)}\|_{L^{2}(\Omega)} \\
& \quad
+ \Big(
\sum_{|n| \ge 1} 
\|\mathcal{P}_{n} v^{(\alpha)}\|_{L^{\infty}(\Omega)}
\Big)\,
\|\nabla \mathcal{Q}_{0} v^{(\alpha)}\|_{L^{2}(\Omega)}\,,
\end{align*}
where the nonlinear estimates \eqref{proof.nonlinear.1} and \eqref{proof.nonlinear.3} with $h$ replaced by $v^{(\alpha)}$ are applied. The proof of Theorem \ref{thm.main.3} is complete. 
\BOX

%

\appendix

\section{Solving the boundary layer equation: proof of Lemma~\ref{constructionGbeta}}\label{appendixairy}
In this section we prove Lemma~\ref{constructionGbeta}.  All the results concerning the Airy function can be found for instance in~\cite{Abramowitz}, Chapter 10.
Let us first recall that a solution to the Airy equation
$$
 \frac{\dd^2 f(\rho)}{\dd \rho^2}  - \rho f (\rho)= 0 \quad \mbox{in}\, \, \R 
$$  is given by the Airy function~$\Ai$, which can be extended as an entire analytic function on $\C$ satisfying
  \begin{equation}\label{Airyeq}
 \frac{\dd^2 f(z)}{\dd z^2}  - z f (z)= 0 \quad \mbox{in}\, \, \C\,.
\end{equation}
 It is    the inverse Fourier transform of
$$
\xi \mapsto \exp \big( \frac{i\xi^3}{3}\big) \,
$$
and satisfies
$$
 \Ai(0) = \frac{1}{ 3^{\frac23}\Gamma\left(\frac23\right)}\, , \quad  \Ai'(0) =-  \frac{1} {3^{\frac13} \Gamma\left(\frac13\right)}\,,$$
where~$\Gamma$ is the Gamma function. Moreover
\begin{equation}\label{decayAiry}
\Ai(z) \sim_{|z|\to \infty}  z^{-\frac14 } \exp \big(
-\frac23 z^\frac32
\big) \, , \quad |\mbox{arg} \, z | <\pi-\epsilon \,,  \quad \epsilon>0\,.
\end{equation}
The results~(\ref{nondegenerate.5}) and~(\ref{nondegenerate.6})  are easy consequences of~(\ref{decayAiry}). 
Indeed we can write 
$$
\begin{aligned}
|C_{0,n,\alpha}|^{-1} = \Big | 
\int_0^\infty e^{\frac{2|n|}{|\beta|} t} \int_t^\infty e^{-\frac{|n|}{|\beta|} s} \Ai (c_- s + \lambda^2) \, \dd s\dd t
\Big | \,,
\end{aligned}
$$
where $\lambda=\frac{|n|c_+}{|\beta|}$ with $|\beta|=(2|\alpha n|)^\frac13$ and $c_\pm =\frac{\sqrt{3}\pm i}{2}$ (hence $\frac{in|\beta|c_-}{2\alpha}=\lambda^2$),
and therefore $|C_{0,n,\alpha}|^{-1}$ is bounded uniformly in the set $\{(n,\alpha)\in \Z\times \R ~|~|\alpha|\geq 1\,, ~ 1\leq |n|\leq |\alpha|^\frac12\}$ thanks to~(\ref{decayAiry}). The result~(\ref{nondegenerate.6})  is obtained in the same way.
As to~(\ref{nondegenerate.4}), it is  known (see for instance~\cite{Abramowitz} page 449)  that 
$$
\int_0^\infty   \Ai (s  ) \dd s = \frac13\,,
$$
so the result follows by continuity of the map~$\mu \mapsto\displaystyle \int_0^\infty e^{-\mu s} \Ai (s + \mu^2 ) \dd s$. This
concludes the proof of Lemma~\ref{constructionGbeta}.    \BOX

\section{Proof of the interpolation inequality \eqref{proof.apriori.middle.9}}\label{appendix.interpolation}

We may assume that $g\in W^{1,2}((1,\infty); r\dd r)$ is nontrivial. 
Let $\delta\in (0,1]$ be a fixed number which will be determined later.
Then we have 
\begin{align}\label{proof.appendix.interpolation.1}
\|g\|_{L^2(\Omega)}^2 & = 2\pi \int_1^\infty |g|^2 r\dd r \nonumber \\
& \leq  C \int_1^{1+\delta} r\dd r \| g\|_{L^\infty( (1,\infty))}^2 + \frac{C}{\delta} \int_{1+\delta}^2 \frac{r^2-1}{r^2} |g|^2 r\dd r  + C \int_2^\infty \frac{r^2-1}{r^2} |g|^2 r\dd r \nonumber \\
& \leq C \delta \| g\|_{L^\infty ((1,\infty))}^2 + \frac{C}{\delta} \| \frac{\sqrt{r^2-1}}{r} g\|_{L^2(\Omega)}^2 + C \| \frac{\sqrt{r^2-1}}{r} g\|_{L^2 (\Omega)}^2\,.
\end{align}
Let us take 
\begin{align*}
\delta= \frac{\| \frac{\sqrt{r^2-1}}{r}g\|_{L^2(\Omega)}}{\| g\|_{L^\infty((1,\infty))} + \| \frac{\sqrt{r^2-1}}{r}g\|_{L^2(\Omega)}}\,.
\end{align*}
Then 
\begin{align}\label{proof.appendix.interpolation.2}
\|g\|_{L^2(\Omega)}^2 
& \leq C \| \frac{\sqrt{r^2-1}}{r}g\|_{L^2(\Omega)} \| g\|_{L^\infty ((1,\infty))} + C \| \frac{\sqrt{r^2-1}}{r} g\|_{L^2 (\Omega)}^2\,.
\end{align}
The estimate \eqref{proof.appendix.interpolation.2} combined with the standard interpolation inequality
$$ \|g\|_{L^\infty((1,\infty))}\leq C \| \partial_r g \|_{L^2((1,\infty))}^\frac12 \| g\|_{L^2((1,\infty))}^\frac12 \leq C \| \partial_r g\|_{L^2 (\Omega)}^\frac12 \| g\|_{L^2 (\Omega)}^\frac12$$
yields
\begin{align*}
\| g\|_{L^2 (\Omega)}^2 \leq C \| \frac{\sqrt{r^2-1}}{r}g\|_{L^2(\Omega)}^\frac43 \|\partial_r  g\|_{L^2 (\Omega)}^\frac23 + C \| \frac{\sqrt{r^2-1}}{r} g\|_{L^2 (\Omega)}^2\,.
\end{align*}
The proof of \eqref{proof.apriori.middle.9} is complete.

\section*{Acknowledgement}
The third author is partially supported by JSPS Program for Advancing Strategic International Networks
to Accelerate the Circulation of Talented Researchers, 'Development of Concentrated Mathematical Center Linking to Wisdom of  the Next Generation', which is organized by Mathematical Institute of Tohoku University.
The second and third authors are grateful to Universit{\'e} Paris Diderot for their kind hospitality during their stay in the spring semester of 2017.

\end{document}